\documentclass[12pt]{article}
\usepackage{amsthm}
\usepackage{latexsym}
\usepackage{fleqn}
\usepackage[latin1]{inputenc}
\usepackage{amssymb}
\usepackage{amsfonts}

\newcommand{\bR}{{\mathbb{R}}}

\newcommand{\tg}{\tilde{\gamma}}

\newcommand{\cF}{{\cal F}}
\newcommand{\cS}{{\cal S}}
\newcommand{\cB}{{\cal B}}

\newcommand{\capp}{{\rm cap}}

\newcommand{\arc}{{\rm arc}}

\newcommand{\cone}{{\rm cone}}
\newcommand{\wk}{\mbox{$\,<$\hspace{-5pt}\footnotesize )$\,$}}
\usepackage[T1]{fontenc}           
\usepackage{textcomp}              
\usepackage{epsfig}        

\theoremstyle{definition} 
\begin{document}

\begin{center}
\textbf{\Large Ball polytopes and the V\'{a}zsonyi problem}\\[0.3cm]
\emph{Yaakov S. Kupitz\footnote{Partially supported by the Landau
Center at the Mathematics Institute of the Hebrew University of
Jerusalem Minerva Foundation, Germany, and by Deutsche
Forschungsgemeinschaft.}, Horst Martini, Micha A. Perles}
\end{center}

\begin{abstract}
Let V be a finite set of points in Euclidean $d$-space $(d \ge 2)$. The
intersection of all unit balls $B(v,1)$ centered at $v$, where $v$ ranges
over $V$, henceforth denoted by $\cB(V)$
is the \emph{ball polytope} associated with $V$. Note that $\cB(V)$ is non-empty iff the
circumradius of $V$ is $\le 1$. After some preparatory discussion on spherical convexity
and spindle convexity, the paper focuses on two central themes\footnote{For the sake of priority it should
be noticed that, in the present form (except for a few chances in the footnotes and one additional reference),
this paper was submitted to the journal ``Discrete and Computational Geometry'' in February 1, 2008.}.
\begin{enumerate}
\item[a)] Define the boundary complex of $\cB(V)$ (assuming it is
non-empty, of course), i.e., define its vertices, edges and facets in
dimension $3$ (in dimension $2$ this complex is just a circuit), and
investigate its basic properties\footnote{This
face structure turns out to be quite intricate or, quoting \cite{Bezdek}, p. 202, lines
22-23: ``the face structure of these objects [ball polytopes] is not
obvious at all''.}.
\item[b)] Apply results of this investigation to characterize finite sets of diameter $1$ in
(Euclidean) $3$-space for which the diameter is attained a maximal
number of times as a segment (of length $1$) with both endpoints in $V$. A
basic result for such a characterization goes back to Gr\"{u}nbaum, Heppes and
Straszewicz, who proved independently of each other, in the late 1950's by
means of ball polytopes, that the diameter of $V$ is attained at most $2|V|-2$
times, thus affirming a conjecture of V\'{a}zsonyi from circa 1935. Call $V$
\emph{extremal} if its diameter is attained
this maximal number $(2|V|-2)$ of times. We extend the aforementioned basic
result by showing that $V$ is extremal iff $V$ coincides with the set of
vertices of its ball polytope $\cB(V)$ and show that in this
case the boundary complex of $\cB(V)$ is self-dual in some
strong sense. The problem of constructing new types of extremal
configurations will be addressed in a subsequent paper, but we do present
here some such new types.
\end{enumerate}
\end{abstract}

\textbf{\small MSC Numbers (2000):} 52C10, 52A30, 52B10, 52B05, 05C62, 52A40


\section{Introduction}

Let $V$ be a set of $n\,(\ge 2)$ points in the Euclidean
$d$-space. The \emph{diameter} of $V$, diam$V$, is the maximum of
the (Euclidean) distances between points of $V$. Denote by $e(V)$
the number of pairs $\{x,y\} \subset V$ such that $\|x-y\| = {\rm
diam} V$, and define $e(d,n)$ to be the maximum of $e(V)$ over all
sets $V$ of $n$ points in ${\mathbb R}^d$.
We call $V \,(V \subset {\mathbb R}^d, \# V = n)$ an
\emph{extremal configuration} if $n > d$ and $e (V) = e(d,n)$. The
case $d=1$ is not interesting $(e(1,n) = 1$ for $n \ge 2$).

For $d=2$, it is well known (see \cite{Pa-Ag}, pp. 213-214) that
$e(2,n) = n \,\, (n \ge 3)$, and the extremal configurations are
fully understood (cf. \cite{Woo} and \cite{Ku 1}). Namely, $V
\subset {\mathbb R}^2$ is extremal if and only if ${\rm vert}P
\subseteq V \subseteq {\rm bd}P$ for some Reuleaux polygon $P$.
When $P$ is a Reuleaux $k$-gon $(3 \le k \le n, \, k$ odd), the
diameters of $V$ form a self-intersecting $k$-circuit (the main
diagonals of $P$) plus $n-k$ ``dangling edges''. Each dangling
edge connects a vertex of $P$ with an interior point of the
opposite arc. See Figure 1, where $k=5$ (a Reuleaux pentagon) and
$n=9$. \vspace*{0.4cm}

\begin{picture}(280,130)
\put(65,0){\includegraphics[scale=0.5]{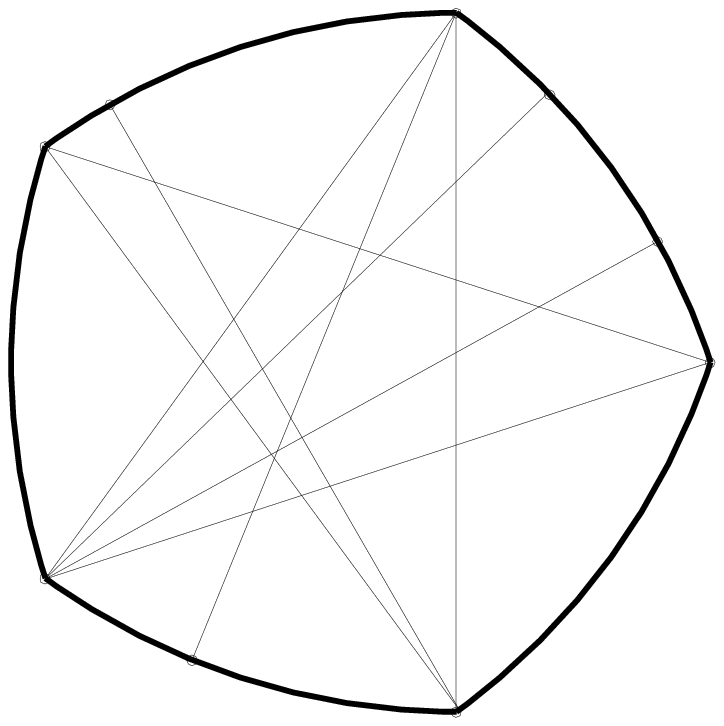}} 
\end{picture}
\vspace{-1cm}
\begin{center}
-- Figure 1 --
\end{center}

If $d=3$, then $e(3,n) = 2n-2$ for $n \ge 4$. This has been
conjectured by Andrew V\'{a}zsonyi (ne\'{e} Endree Weiszfeld\footnote{Also known as the originator of the ``Weiszfeld
algorithm'', to approximate successively the point of minimum sum
of distances to a finite set of points in ${\mathbb R}^3$ (called
the generalized Fermat-Torricelli point).}, see \cite{Er 1}) and
proved in the late 1950's independently by B. Gr\"{u}nbaum, A.
Heppes and S. Straszewicz (see \cite{Pa-Ag}). So we refer to the
equality $e(3,n) = 2n-2$, $n \ge 4$, as the GHS-\emph{Theorem}.

The extremal configurations for $d=3$, however, still form a
mostly uncharted territory. Very few examples have been described
in the literature\footnote{See \cite{Neaderhouser}, Examples 3.1,
3.2, 3.3. Example 3.1 is generalized in our Example 1.1 below;
Example 3.2 is the ``suspended $(2k-1)$-gon'' described in Example
1.2 below. Their Example 3.3 of $7$ points is brilliant, and as we
shall see in another occasion it is obtained from their Example
3.1 for $n=4$ ($4$ vertices of a tetrahedron) by an operation
called  ``ball truncation'' described in \cite{Ku-Ma-Pe}; cf. also Example
9.2 and Remark 9.1 below for the flavor of this operation.}. The
purpose of this paper is to shed some light on the extremal
configurations in dimension 3. We call them the \emph{extremal
configurations} and the problem of describing them \emph{the
V\'{a}zsonyi problem}.

For $d \ge 4$ the quantitiy $e(d,n)$ grows quadratically with $n$.
More precisely, if $d=2m$ or $d=2m+1$ ($d$ fixed) and $n \to
\infty$, then $e(d,n)$ is asymptotically $(1-\frac{1}{m}) \cdot
\frac{n^2}{2}$; see \cite{Pa-Ag}.

It will be convenient to use the notion of diameter graph.

\textbf{Definition 1.1 (Diameter graph):} For $V \subset {\mathbb
R}^d$ the geometric graph $D(V) = \langle V,E \rangle$ whose set
of vertices is $V$ and whose edges are the pairs $\{x,y\} \subset
V$ with $\|x-y\| = {\rm diam} V$ is the \emph{diameter graph} of
$V$.

Let us start with a description of the simplest example of a set
$V \subset {\mathbb R}^3$ with $\# V = n \ge 4$ and $e(V) = 2n-2$.

\textbf{Example 1.1:} When $n=4$, let $V = \{p_0,p_1,p_2,p_3\}$ be the set of vertices
of a regular tetrahedron of edge length 1, say. The diameter graph
of $V$ is $K_4$, and the diameter $1$ is attained exactly six $(=
2 \cdot 4 - 2)$ times.
For $\{i,j,k,l\} = \{0,1,2,3\}$ denote by $C_{kl}$ the circle (of
radius $\frac{1}{2} \sqrt{3}$) with center $\frac{1}{2} (p_k +
p_l)$ that passes through $p_i$ and $p_j$, and define $A_{ij}$ to
be the closed short circular arc of $C_{kl}$ with endpoints $p_i$
and $p_j$. Define also $A^\circ_{ij} = A_{ij} \setminus \{p_i,
p_j\}$. Then $A^\circ_{ij}$ is the locus of all points $x \in
{\mathbb R}^3$ that satisfy
\[
\| x-p_k\| = \| x - p_l \| = 1\,, \,\, \|x-p_i\| < 1\,, \,
\|x-p_j\| < 1\,.
\]
For $n > 4$, take $V$ to be the union of $\{p_0,p_1,p_2,p_3\}$ and
a set of $n-4$ additional points on some $A^\circ_{ij}$, say on
$A^\circ_{12}$. The diameter $1$ is attained six times at the
edges of the tetrahedron, and twice more $(\|x-p_3\| =
\|x-p_4\|=1)$ for each additional point $x$. So far this
configuration was known to Gr\"unbaum \cite{Gr}, Heppes \cite{He}
and Straszewicz \cite{Str} (see also \cite{Neaderhouser}, Example
3.1), and here is a generalization of it.

In \cite{Sch-Pe-Ma-Ku},
Lemma 4.2, we have shown in detail that if $x \in A^\circ_{ij}\,,
\, y \in A^\circ_{i',j'}$, then $\|x-y\| < 1$ provided $\{i,j\}
\cap \{i',j'\} \not= \emptyset$ (i.e., $x$ and $y$ lie on the same
arc or on adjacent arcs), whereas $\|x-y\|
> 1$ if $\{i,j\} \cap \{i',j'\} = \emptyset$, i.e., $x$ and $y$
lie on opposite arcs. (A more detailed form of this lemma with a geometric
proof is given in Lemma 8.1 below.)
Thus we obtain an extremal configuration by
taking the vertices $p_0,p_1,p_2,p_3$ of the tetrahedron, and
merely choosing $n-4$ additional points on three mutually adjacent arcs
$A^\circ_{ij}$. These three arcs either share a vertex (like
$A^0_{01}, A^0_{02}, A^0_{03}$), or miss a vertex (like $A^0_{12},
A^0_{23}, A^0_{31})$.

Another simple class of extremal configurations can be described
as follows.

\textbf{Example 1.2 (the suspended \boldmath$(2k-1)$\unboldmath-gon)}
(cf. \cite{Neaderhouser}, Example 3.2):
Let $P = P_{2k-1} \, (k \ge 2)$ be a regular convex
$(2k-1)$-gon of diameter $1$ located in the $(x,y)$-plane of
${\mathbb R}^3$ and with center at the origin $o$. The diameter of
$P$ is attained $2k-1$ times at the main diagonals of $P$. On the
$z$-axis mark a point $c$ at distance $1$ from the vertices of $P$
and put $V = \{c\} \cup {\rm vert} P$. Then $n =_{\rm def} \# V =
2k$, conv$V$ is a pyramid and the diameter of $V$ is attained
$4k-2 = 2n-2$ times: at the main diagonals of the base $P$, and at
the edges emanating from the apex $c$. Note that the diameter
graph $D(V)$ of $V$ is an $n$-wheel in the sense of graph theory.
(An $n$-wheel is a graph consisting of an $(n-1)$-circuit plus one
universal vertex, the center of the wheel.) For later reference we call
this example \emph{a suspended $(2k-1)$-gon}.
When $k \ge 3$, this
example allows some perturbation of the points of $V \setminus
\{c\}$ on the unit sphere centered at $c$. A detailed account of
this possibility will be given in a subsequent paper.

\textbf{Remark:} The suspended $(2k-1)$-gon of Example 1.2 serves in
Example 9.2 below to construct yet another family of extremal configurations
by a construction called ``ball truncation''; cf. also Remark 9.1,~2) below.

By these examples one sees that $e(3,n) \ge 2n-2$ for $n \ge 4$.

The existing proofs that $e(3,n) \le 2n-2$ for $n \ge 4$ proceed
by induction on $n$. The case $n=4$ is obvious. If $\#V =  n+1$
and a point $x \in V$ is incident with at most two diameters
(i.e., at most two points of $V$ are at distance diam$V$ from
$x$), then we apply the induction hypothesis to $V \setminus
\{x\}$ and find that $e(V) \le (2n-2)+2 = 2(n+1)-2$. Thus the
only case that needs some extra consideration is when each point
$x \in V$ is at distance diam$V$ from at least three other
points of $V$. This ``extra consideration'' consists of three
steps:
\begin{enumerate}
\item Construction of a convex body ${\cal B}(V)$ (``${\cal B}$''
stands for ``Ball'') related to $V$. In fact, ${\cal B}(V)$ is
defined as the intersection of all balls of radius diam$V$ with
centers in $V$.
\item Defining a ``spherical face structure'' on
the boundary of ${\cal B}(V)$ and establishing a relation between
this face structure and the diameter graph $D(V)$ of $V$.
\item Applying Euler's polyhedral formula $(v-e+f=2)$ to the face
structure of ${\cal B}(V)$, to conclude that $e(V) \le 2n-2$.
\end{enumerate}
We follow these steps carefully and conclude that if $e(V) =
2n-2$, then the face structure of ${\cal B}(V)$ is \emph{strongly
self-dual}, i.e., it admits a self-duality of order 2 (involution)
which is strong in the following sense: it is fixed-point free when acting as an automorphism of the
first barycentric subdivision of the boundary cell complex of
${\cal B}(V)$. Nevertheless, the face structure of ${\cal B}(V)$, i.e., the poset of its faces,
may differ considerably from that of an ordinary convex
$3$-polytope. The intersection of two facets may have more than
one connectivity component, and there may be digonal facets. Consequently the set of faces, partially
ordered by inclusion, need not be a lattice\footnote{So it is \emph{mistakenly} called a ``face
lattice'' systematically in \cite{Bez-Nasz}, p. 256, lines 21, 29, 36, p. 260,
line 15, and in \cite{Bezdek}, p. 211, lines 24, 30. In \cite{Bezdek},
Corollary 6.8 (p. 213), it is shown that, under a strong assumption on $V$ explained in footnote 15 below, $\cB(V)$ is a
lattice indeed, but this assumption takes the discussion far from the generic case.}, and its $1$-skeleton
forms a planar graph that is not necessarily $3$-connected. This
is true even if $V$ is a critical configuration (cf. Proposition and Definition
2.1 below). There is a critical configuration $V \subset \bR^3\,, \, \# V = 8$, for
which $\cB(V)$ is not combinatorially equivalent to a $3$-polytope. The
configuration will be described in a subsequent paper.
There we shall construct new examples of
extremal configurations that are modelled after some types of
strongly self-dual convex $3$-polytopes.

\textbf{Notation:} Our considerations take place in Euclidean
$d$-dimensional space ${\mathbb R}^d$, $d \ge 2$, with origin $o$, mainly
referring to the case $d=3$.
We abbreviate, as usual, aff, conv, vert, int, relint, bd,
relbd for \emph{affine hull}, \emph{convex hull, vertex set,
interior, relative interior, boundary}, and \emph{relative
boundary}, respectively.

For $p \in {\mathbb R}^3$ define by $B(p) =_{\rm def} \{x \in
{\mathbb R}^3: \|x-p\|\le 1\} = p + B(o)$ the \emph{solid ball
centered at} $p$, where $\|\cdot\|$ is the usual Euclidean norm in
${\mathbb R}^3$, and by $S(p) =_{\rm def} {\rm bd} B(p) = \{x \in
{\mathbb R}^3 : \|x-p\| = 1\}$ the \emph{unit sphere centered at}
$p$, usually abbreviated $S_p$.

The \emph{circumball} of a bounded set $W \subset {\mathbb R}^3$
is the (unique) closed ball of smallest radius that includes $W$.
We denote by cr$(W)$ (read: circumradius of $W$), resp. cc$(W)$
(read: circumcenter of $W$), the radius, resp. center, of the
circumball, and define cr$(W) = 0$ when $\# W=1$, and cr$(W) =
\infty$ when $W$ is unbounded. (When cr$(W) = \infty$, cc$(W)$ can
be taken as any point of ${\mathbb R}^3$.)

\section{Critical configurations}

Let $V$ be an extremal configuration for the V\'{a}zsonyi problem
with $n$ points, $n \ge 4$. If $n=4$, then the diameter  graph
$D(V)$ is just $K_4$. If $n \ge 4$, choose a point $v \in V$ of
minimal valence, say $\delta$, in $D(V)$. Then $e(V \setminus
\{v\}) = e (V) - \delta = 2 n- 2 - \delta$. Note that diam$(V
\setminus \{v\}) = {\rm diam} V$, since $\delta \le n-1 < 2 n-2$,
and thus the removal of $v$ does not obliterate all diameters of
$V$. Assuming the GHS-Theorem (an extended version of which is
proved in Theorem 7.1 below), $e(V \setminus \{v\}) \le 2 (n-1) -2$. It
follows that $\delta \ge 2$. Moreover, if $\delta = 2$, then $V
\setminus \{v\}$ is again extremal. This leads us to the following

\textbf{Proposition -- Definition 2.1 (Critical Configuration):}
\emph{Every point of an extremal configuration (for the
V\'{a}zsonyi problem) is incident with at least two diameters. If
every point of an extremal configuration is incident with at least
three diameters, then $V$ is \textbf{critical}}.

\textbf{Examples:} The suspended $(2k-1)$-gon (Example 2.1 above) yields
a critical configuration. The $7$-point set described in \cite{Neaderhouser}, Example 3.3 (mentioned above),
is critical as well.

\textbf{Corollary 2.1:} \emph{The diameter graph $D(V)$ of an
extremal configuration $V$ does not contain two adjacent
$2$-valent vertices}.

\textbf{Proof:} If $v$ and $w$ are adjacent $2$-valent vertices of
$D(V)$, the removal of $v$ leaves $w$ as a $1$-valent vertex of
$D(V \setminus\{v\})$, and thus $V \setminus \{v\}$ is not
extremal, by Proposition 2.1. So $D(V \setminus \{v\})$ has less
than $2 \# (V \setminus \{v\})-2$ diameters. Hence $D(V)$ has less
than $2 \# (V \setminus \{v\})-2 + 2 = 2 \# V-2$ diameters, i.e.,
$V$ is not extremal, a contradiction. \hfill \rule{2mm}{2mm}

If an extremal configuration $V$ (with $\# V \ge 5)$ is not critical,
then successive removal of $2$-valent vertices $v_1, v_2, \dots,
v_k$ (the vertex $v_i$ being $2$-valent in $D(V \setminus \{v_1,v_2,
\dots,v_{i-1}\}))$ will finally produce a critical
subconfiguration $V_k$ with $n-k$ points. This will happen after
at most $n-k$ steps. Note that the sequence $v_1, v_2, \dots, v_k$
of removed vertices is not necessarily unique. We shall see,
however, that the final critical subconfiguration is unique.
Denote the set $V_k$ by $C$ ($C$ stands for ``critical''). This
set $C$ has the following properties:
\begin{enumerate}
\item[1)] $C \subseteq V$, and $C$ is extremal,
\item[2)] $D(C)$ is a spanned subgraph
of $D(V)$.
\end{enumerate}
This follows from the fact, mentioned above, that removal of one
point from an extremal configuration does not affect its diameter,
and that removal of a $2$-valent vertex from an extremal
configuration leaves it extremal.
\begin{enumerate}
\item[3)] Every vertex of $D(C)$ has valence
$\ge 3$.
\item[4)] For any set $V',\, C \varsubsetneqq V'
\subseteq V$, the graph $D(V')$ has a vertex of valence $\le 2$.
\end{enumerate}
In fact, if $v_i$ is the first vertex in the sequence $v_1, v_2,
\dots, v_k$ of removed vertices that belongs to $V'$, then the
valence of $v_i$ in $D(V')$ is smaller than or equal to the valence
of $v_i$ in $D(V \setminus \{v_1, v_2, \dots, v_{i-1}\})$, which
is $2$.

Next we show that $C$ is uniquely determined by the properties 1) --
4) above in the following strong sense.

\textbf{Theorem 2.1:} \emph{If $C'$ is a critical subconfiguration
of $V$, then $C' \subseteq C$}.

\textbf{Proof:} Assume, by r.a.a., that $C'$ is a critical
subconfiguration of $V$ and $C' \nsubseteq C$. This means that $C
\cup C' \varsupsetneqq C$. Now consider property 4), with $V' = C
\cup C'$. There is a vertex $x \in C \cup C'$ whose valence in
$D(C \cup C')$ is $\le 2$. This is impossible, since if $x \in C$
or $x \in C'$, then the valence of $x$ in $D(C)$ or in $D(C')$,
respectively, is already $\ge 3$. \hfill \rule{2mm}{2mm}

\textbf{Corollary 2.2:} \emph{The set $C$ is the unique maximal
critical subconfiguration of $V$, where  ``maximal'' here is with
respect to inclusion as well as to cardinality}.

The question whether a critical configuration may still include an
extremal configuration as a proper subset (and therefore includes
a smaller critical configuration) turns out to be
delicate. The answer is (surprisingly) ``yes''.
The above mentioned example
$V \subset \bR^3, \# V=8,  {\rm diam} V=1$, for which $\cB(V)$ is not
$3$-polytopal, is critical and has an extremal subconfiguration of four points.
Meanwhile we just introduce the relevant technical term.

\textbf{Definition 2.2:} An extremal configuration $V$ is
\emph{strongly critical} if no proper subset of $V$ is extremal.

We conjecture that an extremal configuration
is strongly critical iff the spherical face structure of its Ball polytope
is isomorphic to that of an ordinary $3$-polytope; cf. Conjecture 9.2 below.

The only critical configurations we have seen so far are the
suspensions of regular $(2k-1)$-gons $(k \ge 2)$, whose diameter
graph is a $2k$-wheel. For $k=2$, these include the set of
vertices of a regular tetrahedron. In the planned subsequent paper we shall
construct a large variety of critical configurations, with an
arbitrarily large number of points. Most of the rest of this paper
is devoted to the study of Ball polytopes and their relation to
the extremal configurations. A detailed account of this confronts us
with the notions of spindle convexity and spherical convexity discussed in paragraphs 3 and 4 below.

\section{Spindle convexity, ball hulls and the Ball connection}

In this section we introduce the notion of spindle convexity and
the operations $S \to {\frak B} (S)$ (Ball hull) and $S \to \cB(S)$
(Ball set) for subsets $S$ of ${\mathbb R}^3$, and
establish the basic properties of these notions and their
interrelations.

\textbf{Definition and Notation 3.1.:} A set $S \subset {\mathbb
R}^3$ is \emph{spindle-convex} (abbr. $s$-\emph{convex}) if for
any two distinct points $a, b \in S$, and for any circle $C$ of
radius $r\ge 1$ that passes through $a$ and $b$, the short
circular arc of $C$ with endpoints $a,b$ is also included in $S$.
We denote by $sp(a,b)$ (read: spindle $a,b$) the union of these
arcs including the segment $[a,b]$.

\textbf{Remarks 3.1:}
\begin{enumerate}
\item If $S \subset {\mathbb R}^3$ is $s$-convex, then $S$ is strictly convex.
The proof is easy and left to the reader, with just one hint: If $S$ is $s$-convex,
$a,b \in S$ and
$y \in {\rm relint}[a,b]$, then $y$ is contained in a short
circular arc of radius $r \ge 1$ whose two endpoints, $x,z$ say,
belong to $sp(a,b) \setminus [a,b] \subset S$, and by the
$s$-convexity of $S$ we have $y \in sp(x,z) \setminus [x,z]
\subset S$.
\item If $\|a-b\|<2$, then the union of all short circular arcs of
radius $r \ge 1$ with endpoints $a$ and $b$, including the segment
$[a,b]$, forms a spindle with endpoints $a$ and $b$. When $\|a-b\|
\ge 2$, the union of these short circular arcs, including
semicircles with endpoints $a$ and $b$, forms a solid ball
centered at $c=\frac{1}{2}(a+b)$ with radius $\frac{1}{2}
\|a-b\|$\footnote{So it is \emph{not} $\bR^3$ as stated in \cite{Bezdek},
p. 203, line 11.}.
In both cases we denote the spindle with endpoints $a$
and $b$ by ${\rm sp}(a.b)$.
\item Spindle convexity was studied under the name ``\"Uberkonvexit\"at''  by
Mayer \cite{Ma}. We introduced this notion in \cite{Ku-Ma-Pe}; cf. also \cite{Bezdek},
p. 203.
\end{enumerate}

An immediate consequence of this definition is

\textbf{Proposition 3.1:} \emph{The intersection of any family of
$s$-convex sets in ${\mathbb R}^3$ is $s$-convex}.

\textbf{Lemma 3.1 ($s$-convexity of the unit ball):} \emph{A
closed ball of radius $1$ in ${\mathbb R}^3$ is $s$-convex}.

Lemma 3.1 will follow from a refined version, namely Lemma 3.1' below. We
will need this refined version in the sequel (e.g., in the proof of Proposition 5.1 (i)).

\textbf{Lemma 3.1' (simplified Sallee Lemma):}\emph{ Let $B = B(z,1)\subset {\mathbb
R}^3$ be a closed unit ball centered at $z$.  Let $a,b$ be two
distinct points in $B$. Let $C$ be a circle of radius $r \ge 1$
with center $c$ that passes through $a$ and $b$. Let $A \subset C$
be the open short circular arc with endpoints $a,b$. (When $a,b$
are opposite points of $C$, $A$ denotes any of the two open
semicircles of $C$ with endpoints $a,b$.) Then $A \subset {\rm
int} B$, with just one exception: If $a, b \in \partial B$ and $r
= 1$, then $c=z$ and} $C \subset
\partial B$\footnote{This lemma goes back essentially to Sallee in \cite{Sa},
Lemma 3.3, (p. 317). The lemma appears there in a rather distorted form, which is repeated in
\cite{Bez-Nasz}, Lemma 1.1 (SL1). This last formulation also lacks the delicate
exceptional case $a, b \in \partial B$ and $r=1$ appearing in our formulation. Some
overlaps between \cite{Sa}, \cite{Bezdek}, \cite{Bez-Nasz}, and the present paper do exist, and we will point out on
them occasionally in footnotes.}.

\textbf{Proof:} Put $H =_{\rm def}$ aff$C$ and $C' =_{\rm def} H
\cap \partial B, D' =_{\rm def} H \cap B = {\rm conv} C'$. $C'$ is
a circle of radius $r' \le 1$, with center $c'$. Note that $r'<1$,
unless $H$ passes through the center $z$ of $B$, in which case
$r'=1$ and $c'=z$. Now consider the circle $C$ (with radius $r \ge
1$ and center $c$) and the circle $C'$ (with radius $r' \le 1$ and
center $c'$). Both these circles lie in $H$.

If $a \in {\rm int} B$ or $b \in {\rm int} B$, then $C$ meets
relint$D'$. But $C \nsubseteq D'$, since $r \ge r'$, hence $C$
meets $C'$ at two points, which divide $C$ into two unequal open
arcs. The short arc is contained in ${\rm relint} D'$ and includes
$A (= \widehat{ab})$, and the long one is exterior to ${\rm
relint} D'$. The same holds true if both $a$ and $b$ are boundary
points of $B$ (i.e., $a,b \in C \cap C'$) provided $r > r'$, and
even when $r=r'(=1)$ provided $C\neq C'$ (i.e., $c\neq c'$)).
There remains the case where $a,b \in \partial B, r'=r=1$ (hence
$c'=z$), and $c=c'$. In this case $C'=C$, and thus $A \subset
\partial B$. \hfill \rule{2mm}{2mm}

\textbf{Corollary 3.1.:} \emph{The intersection of any family of
closed unit balls in ${\mathbb R}^3$ is $s$-convex}.

Under closedness assumption the converse statement of Corollary 3.1 holds as well;
see Theorem 3.1 below.

\textbf{Definition 3.2 (Ball hull):} For a set $S\subset {\mathbb
R}^3,$ the intersection of all closed unit balls in ${\mathbb
R}^3$ that include $S$, denoted by ${\frak B} (S)$, is
the {\it Ball hull} of $S$.

\textbf{Remarks:} It follows that ${\frak B}
(\emptyset)=\emptyset$ and $S \subset T \Longrightarrow {\frak
B}(S) \subset {\frak B} (T)$. If the circumradius cr$(S)$ of $S$
is $>1$, then $S$ is not included in any unit ball, and ${\frak
B}(S) = {\mathbb R}^3$ (= the intersection of the empty family of
unit balls in ${\mathbb R}^3$). If cr$(S) = 1$, then ${\frak
B}(S)$ is just the unique closed unit ball that includes $S$.

In all cases $S \subset {\frak B} (S)$ and ${\frak B} (S)$ is
$s$-convex. Note that $S = {\frak B}(S)$ iff $S$ is the
intersection of some family of closed unit balls in ${\mathbb
R}^3$. Next we show that every $s$--convex set in ${\mathbb R}^3$
is the intersection of a family of closed unit balls, i.e., $S$ is
$s$--convex iff $S={\frak B}(S)$.

\textbf{Theorem 3.1:} \emph{Suppose that $S$ is a closed subset of
${\mathbb R}^3$. Then $S$ is $s$-convex iff $S$ is the
intersection of a family of (closed) unit balls in} ${\mathbb
R}^3$\footnote{The easier ``if'' part of this theorem seems to be
overlooked until now, and the heavier ``only if'' part seems to appear already in \cite{Ma}.
We formulated the ``only if'' part in \cite{Ku-Ma-Pe}, Theorem 1, and it also appears in
\cite{Bezdek}, Lemma 3.1, with a small gap in the proof. All these appearences are independent of each other.}.

\textbf{Remarks 3.2:}
\begin{enumerate}
\item The family may be empty, in which case $S=\bR^3$. So
an alternative formulation of Theorem 3.1 is: Suppose $S \subset \bR^3$ is closed.
$S$ is $s$-convex iff either $S = \bR^3$ or $S$ is the intersection of
a nonempty family of unit balls in $\bR^3$.
\item This theorem is the $s$-convexity analogue of the theorem from
linear convexity saying that a closed set $S \subset \bR^3$ is convex iff
either $S=\bR^3$ or $S$ is the intersection of a nonempty family of closed
half-spaces.
\end{enumerate}

\textbf{Proof of Theorem 3.1:} The ``if'' part follows from Corollary 3.1 above.
To establish the ``only if'' part, it suffices to prove

\textbf{Proposition 3.3 (``separation out'' by a unit sphere):}
\emph{Suppose $S$ is a closed, $s$--convex proper subset of ${\mathbb
R}^3$. If $z \in {\mathbb R}^3\setminus S$, then there is a closed
unit ball $B = B_z \subset {\mathbb R}^3$ that includes $S$ and
misses} $z$.

Indeed, if Proposition 3.3 holds true, then $S = \bigcap
\{B_z:z\in {\mathbb R}^3 \setminus S\}$. (Note that the family $\{B_z: z \in \bR^3
\setminus S\}$ is empty iff $S=\bR^3$, which is compatible with the equality
$S=\cap\{B_z: z \in \emptyset\}$, since the intersection of an empty set of unit balls
is $\bR^3$.)

\textbf{Proof of Proposition 3.3:}
For $S = \emptyset$ the statement is clear. Assume, therefore, that $S \neq \emptyset$
and denote by $a$ the (unique) point of $S$ nearest to $z$.
Let $u = \frac{a-z}{\|a-z\|}$, and let $H \subset {\mathbb R}^3$ be
the plane through $a$ perpendicular to $u$, i.e., $H = \{ x \in
{\mathbb R}^3 : \langle u,x \rangle  = \langle u,a \rangle\}$. Let
$H^-$ be the closed half--space bounded by $H$ that contains $z$,
and let $H^+$ be the other closed half-space bounded by $H$. It
follows from the convexity of $S$ that $S\subset H^+$.

Note that $a \in S$ and $S$ is $s$--convex, hence strictly convex,
and therefore $S \cap H = \{a\}$, i.e., $S \subset \{a\} \cup {\rm
int} H^+$. The set $\{a\} \cup {\rm int}H^+$ is the union of all
spheres tangent to $H$ whose centers lie on the open ray
$\{a+\lambda u:\lambda > 0\}$. Define $B_z =_{\rm def} B(a+u)$,
the unit ball centered at $a+u$. Clearly $z\notin B_z$. If remains
to show that $S\subset B_z$.

Assume, to the contrary, that there is a point $b \in S \setminus
B_z$. Clearly $b \not= a$. There is a unique $\rho > 1$ such that
$b$ lies on the sphere $S_\rho =_{\rm def} S^2 (a + \rho u, \rho)$
of radius $\rho$ centered at $a + \rho \mu$. (An explicit
calculation shows that $\rho = \frac{\|a-b\|^2}{2 \langle b-a, u
\rangle}$.)

Let $A \subset S_\rho$ be a short (closed) circular arc of radius
$\rho$ with endpoints $a,b$.
\vspace{1cm}

\begin{picture}(280,130)
\put(50,0){\includegraphics[scale=0.5]{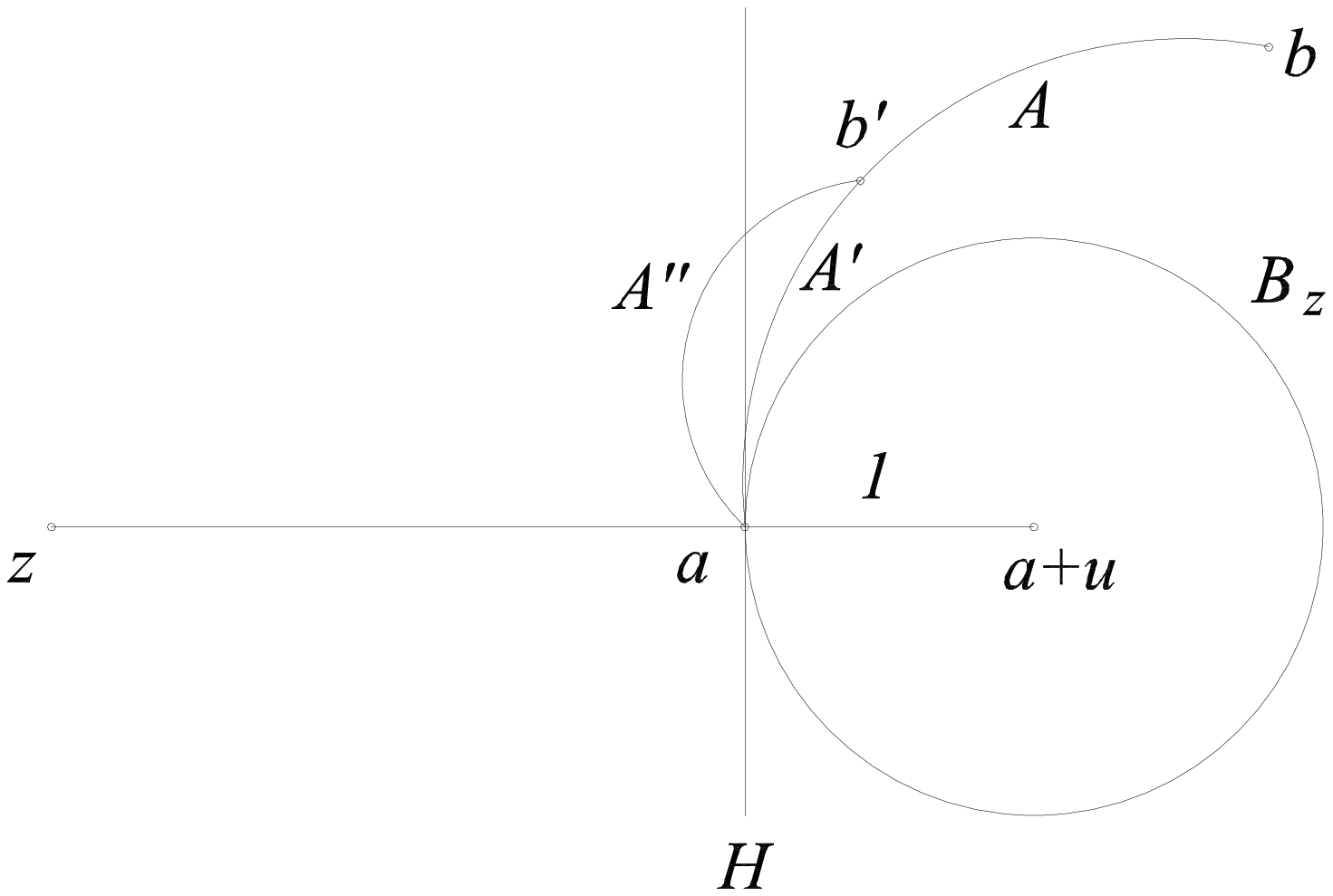}} 
\end{picture}

\begin{center}
-- Figure 2 --
\end{center}

($A$ is uniquely determined, unless $b = a+2 \rho u$, in which
case $A$ is half a great circle (a meridian) on $S_\rho$.) Let
$b'$ be a point on $A$ such that $0 < \| b'-a \| < 2$, and let
$A'$ be the part of $A$ between $a$ and $b'$. In the plane aff$A'$
draw a short circular arc $A''$ of radius $1$ with endpoints $a$
and $b'$, on the same side of the line aff$(a,b')$ as $A'$. Since
$A'\, (\subset H^+)$ is tangent to $H$ at $a$, and $A''$ is more
curved $(1 < \rho), A''$ protrudes into $H^-$ near $a$; see Figure
2. We have $A' \subset A \subset S$, hence $b' \in S$ and $A''
\subset S$, since $S$ is $s$-convex. But $S \subset H^+$, a
contradiction. \hfill \rule{2mm}{2mm}

\textbf{Remark:} Of course, the characterization of $s$-convex
sets given in Theorem 3.1 holds true in ${\mathbb R}^d$ for $d \ge
2$.

\textbf{Definition 3.3 (Ball set):}
For a subset $V$ of ${\mathbb R}^3$, the set ${\cal B}(V) =_{\rm
def} \{y \in {\mathbb R}^3 : (\forall x \in V) [\|x-y\|\le 1]\}$ is
the \emph{Ball set} (or \emph{Ball connection}) of $V$.
Here $\|\cdot\|$ denotes the usual Euclidean norm in $3$-space\footnote{Ball sets in general and Ball polytopes in particular
(cf. Definition 3.4 below) were studied by Sallee in \cite{Sa} in
connection with bodies of constant width (cf. footnote 20~(b) below) and, in a more general setting, in \cite{Ku-Ma-Pe},
\cite{Bez-Nasz}, and \cite{Bezdek}.}.

We can view ${\cal B}(V)$ as the intersection of all
(closed) unit balls centered at points of $V$ or, alternatively,
as the set of centers of all unit balls that include $V$.
Thus we have the following properties.
\begin{enumerate}
\item[(A)] ${\cal B}(V) = \bigcap \{v+B : v \in V\}$, and
\item[(B)] ${\cal B}(V) =\{z \in {\mathbb R}^3 : V \subset z +
B\}$, where $B=_{\rm def} \{x \in {\mathbb R}^3 : \|x\| \le 1\}$.
\end{enumerate}
Here are some more properties of the operation $V \to {\cal B}(V)$:
\begin{enumerate}
\item[(C)] ${\cal B}(\emptyset) = {\mathbb R}^3$ (= the
intersection of the empty set of unit balls is ${\mathbb R}^3$).
\item[(D)] $V \not= \emptyset \Rightarrow {\cal B}(V)$ is compact
and $s$-convex, hence strictly convex (see (A)).
\item[(E)] For $p \in
{\mathbb R}^3\,, \, {\cal B}(\{p\}) = p+B=B(p)$ is the unit ball
centered at $p$.
\item[(F)] cr$(V) = 1 \Rightarrow {\cal B}(V) =
{\rm cc}(V)$, where cc$(V)$ is the circumcenter of $V$ (i.e.,
${\cal B}(V)$ is $0$-dimensional).
\item[(G)]
\begin{enumerate}
\item[(i)] When cr$(V) = 1-\varepsilon < 1\,, \, {\cal B}(V)$
includes a ball of radius $\varepsilon\, (> 0)$ centered at cc$(V)$,
i.e., ${\cal B}(V)$ is full-dimensional.
\item[(ii)] cr$(V) > 1
\Rightarrow {\cal B}(V) = \emptyset$. Hence
\item[(iii)] ${\cal B}(V)$ is full dimensional $\Leftrightarrow$ cr$(V) < 1$\footnote{This simple condition,
equivalent to full dimensionality of ${\cal B}(V)$, seems to be new. Sallee, e.g., assumes in \cite{Sa}, p. 316,
that diam$V \le 1$ (in order to ensure full dimensionality of ${\cal B}(V)$), which implies cr$(V) < 1$; see
Theorem 5.2~a) below. Thus Sallee's condition is sufficient, and clearly not necessary for full
dimensionality of $\cB(V)$.}.
\end{enumerate}
\item[(H)] $V_1 \subset V_2 \Rightarrow {\cal B}(V_2) \subset
{\cal B}(V_1)$ (inclusion reversion, but strict inclusion does not
imply strict reversed inclusion). \item[(I)] ${\cal B}(V) =
{\cal B}({\rm cl~conv} V)$. In fact, ${\cal B}(V) = {\cal B}
({\frak B}(V))$, since $V$ and ${\frak B}(V)$ are included in
exactly the same unit balls. \item[(J)] ${\frak B}(V)$ is the
largest set $X \subset {\mathbb R}^3$ such that ${\cal B}(X) =
{\cal B}(V)$. In fact, if $Y \subset {\mathbb R}^3$ and $Y
\nsubseteq {\frak B}(V)$, then there is a point $y \in Y \setminus
{\frak B}(V)$, and there is a point $z$ such that $V \subset z +
B$ but $y \notin z+B$. Thus $z \in {\cal B}(V)$, but $z \notin
{\cal B}(Y)$, i.e., ${\cal B}(Y) \nsupseteqq {\cal B}(V)$.
\item[(K)] diam$V \le 1 \Leftrightarrow V \subset {\cal B}(V)$.
\item[(L)] diam$V \le 1 \Rightarrow {\cal B}(V) = \{x \in
{\mathbb R}^3: \mbox{ diam} (V \cup \{x\}) \le 1\}$.
\item[(M)]
\begin{enumerate}
\item[(i)] ${\cal B}({\cal B}(V)) = {\frak B}(V)$ and, in
particular,
\item[(ii)] ${\cal B}({\cal B}(V)) \supset V$.
\end{enumerate}
\end{enumerate}

\textbf{Proof:} (i) $x \in {\cal B}({\cal B}(V))$
\vspace*{-0.4cm}
$$
\begin{array}{lll}
& \Leftrightarrow & (\forall y \in {\mathbb R}^3)
\left[\left(\forall v \in V\right)\left(\|y-v\|\le 1\right)
\Rightarrow
\left(\|x-y\|\le 1\right)\right]\\
& \Leftrightarrow & \left(\forall y \in {\mathbb R}^3\right)
\left[\left(V \subset B(y)\right) \Rightarrow \left(x \in
B(y)\right)\right]
\Leftrightarrow  x \in {\frak B}(V)\,.
\end{array}
$$
(ii) Follows from $V \subset {\frak B}(V)$ and (i). \hfill \rule{2mm}{2mm}
\begin{enumerate}
\item[(N)] ${\cal B}({\cal B}({\cal B}(V))) = {\cal B}(V)$.
\end{enumerate}

To prove $(N)$ note that the inclusion ``$\supset$'' follows from (M) (ii)
above, where $V$ is replaced by ${\cal B}(V)$, and ``$\subset$''
follows again from (M) (ii) by using (H) (inclusion reversion).
(Alternatively, use (M) (i) and ${\cal B}(V) = {\cal B}({\frak
B}(V))$ proved in (I) above.)
From (N) we infer:

\textbf{Corollary:} \emph{The operation $V \to {\cal B}(V)$ is
injective (i.e., one to one) on its image set, that is}, ${\cal
B}(S) \not= {\cal B}(T) \Rightarrow {\cal B}({\cal B}(S)) \not=
{\cal B}({\cal B}(T))$.
\begin{enumerate}
\item[(O)] ${\cal B}(V) = V \Leftrightarrow V$ is a set of
constant width $1$. (Hence, a set $V \subset \bR^3$ of constant width $1$ is $s$-convex.)
\end{enumerate}
This characterization of bodies of constant width (through their Ball connection) is known as ``the spherical
intersection property'', proved, e.g., by Eggleston  in \cite{Egg}, Theorem 52 in p. 123,
and in \cite{Egg 2}, pp. 166-167. See also \cite{Cha-Gr}, pp. 62-63,
for a proof of this characterization. (A related notion is studied in \cite{Ku-Mar}.)

\textbf{Definition 3.4 (Ball polytope):} For $V \subset {\mathbb
R}^3$ finite and non-empty with cr$(V) < 1$, the Ball set ${\cal B}(V)$
is the \emph{Ball polytope} associated with $V$.

Note that cr$(V) < 1 \Leftrightarrow {\cal B}(V)$ is full-dimensional
(see (G) (iii) above).

Since a Ball polytope is a finite intersection of balls, it is strictly convex (provided cr$(V) < 1$), i.e.,
$\left(\forall x, y \in {\cal B}(V)\right) \left[x \not= y
\Rightarrow {\small \frac{1}{2}} (x+y) \in \mbox{ int } {\cal
B}(V)\right]\,.
$
The Ball polytope ${\cal B}(V)$ of a finite set $V \subset
{\mathbb R}^3$ (with, say, $\# V > 1$ and cr$(V) < 1$) has a
natural ``spherical'' face structure denoted by ${\cal S}{\cal F}
({\cal B}(V))$. We shall investigate this face structure
thoroughly, and show (in \S~8--9 below) that for extremal configurations $V$ of the
V\'azsonyi problem it
admits a canonical fixed-point free involutory self-duality. This
means an anti-isomorphism (antimorphism, for short) of ${\cal
S}{\cal F} ({\cal B}(V))$ of order 2 which, when applied to the
first barycentric subdivision of ${\cal S}{\cal F}({\cal B}(V))$,
acts as a fixed-point free automorphism. (Note that the first
barycentric subdivision of ${\cal S}{\cal F}({\cal B}(V))$ is a
purely $2$-dimensional simplicial complex which is
order-isomorphic to the collection of all chains (``flags'' in the
current terminology) in ${\cal S}{\cal F}({\cal B}(V))$, partially
ordered by inclusion.)

The exposition of the spherical face structure ${\cal S}{\cal
F}({\cal B}(V))$ of ${\cal B}(V)$ can be viewed as a detailed
account of the proof of V\'{a}zsonyi's conjecture mentioned in the
introduction. The insight gained by this exposition will lead to
the construction of many new types of critical configurations in
in the planned subsequent paper. In
their original proof of the V\'{a}zsonyi conjecture, Gr\"unbaum,
Heppes and Straszewicz (independently) defined the Ball polytope
${\cal B}(V)$ for critical $V$ only (they did not use the term
``Ball polytope'') and treated its face structure very sketchy.
Our more general definition of ${\cal B}(V)$ (for general (finite)
configurations $V \subset {\mathbb R}^3$) is aimed to gain a better
understanding of this face structure.
The next paragraph collects basic results on spherical convexity
and apexed cones, also needed in the sequel.

\section{Spherical convexity and apexed cones}

Let $S(o)$ be the unit sphere centered at the origin.

\textbf{Definition 4.1 (spherical convexity, strict spherical
convexity, spherical convex hull, small cap):} A set $F \subset
S(o)$ is \emph{spherically convex} [\emph{strictly spherically
convex}] if
for every two distinct non-antipodal points $a,b \in F$ the
geodetic line that connects $a$ and $b$ on $S(o)$ (i.e., the short
circular arc with endpoints $a,b$ of the unit circle centered at
$o$ which passes through $a$ and b) is contained in $F$ [is
contained in a ${\rm relint}F$, with possible exception of the
endpoints $a$ and $b$]\footnote{In \cite{Bezdek}, Definition 5.3,
spherical convexity is defined for sets contained in a hemisphere only.
This assumption is not necessary, as we see here, but it is harmless because most spherically convex
sets are contained in a hemisphere.}.

\textbf{Examples 4.1:}  Here is a list of some spherically convex subsets of $S(o)$.
\begin{enumerate}
\item[a)] The empty set and any one-point set.
\item[b)] A set of two antipodal points.
\item[c)] $S(o)$ is strictly spherically convex.
\item[d)] Let $a,b \in S(o)$, $a \not= -b$ (i.e., $a,b$ are not antipodal) and let
arc$(a,b)$ be the geodetic line between $a$ and $b$ on $S(o)$. Then arc$(a,b)$, relint arc$(a,b)$,
arc$(a,b) \setminus \{a\}$, and arc$(a,b) \setminus \{b\}$ are spherically convex. If $a,b \in S(o)$
are antipodal points, then any closed/open/half-open half great circle on $S(o)$ with endpoints
$a,b$ is spherically convex.
\item[e)] A great circle of $S(o)$.
\item[f)] A closed hemisphere of $S(o)$ is spherically convex, but not strictly spherically convex.
\item[g)] A small cap of $S(o)$ is strictly spherically convex (see below).
\item[e)] A spherical triangle on $S(o)$ is spherically convex but not strictly spherically convex,
since its edges are geodetic arcs.
\end{enumerate}
Note that all these examples are the intersection of $S(o)$ with some cone apexed at $o$. In Corollary 4.1 below
it is shown that this situation is typical.

Assume that $V \subset S(o)$. The intersection
of all spherically convex sets that contain $V$ (in $S(o)$) is the \emph{spherical convex hull} of $V$, henceforth denoted
by sph-conv$(V)$. Clearly, sph-conv$(V)$ is spherically convex and $V \subset W \Rightarrow$ sph-conv$(V) \subset$
sph-conv$(W)$.

Let $H$  be a plane that misses the origin $o$ such that dist$(o,H) \le 1$.
The intersection of $S(o)$ with the closed half-space bounded by
$H$ that misses $o$ is a \emph{small cap} of $S(o)$. A small cap is obviously strictly spherically convex.

\textbf{Definition 4.2 (apexed cone):} Let $V \subset {\mathbb
R}^3 \setminus \{o\}$ be a compact set. The cone cone$_o V$ is
\emph{apexed} (at $o$) if there exists a plane $H \subset \bR$ strictly
separating $V$ and $o$. Equivalently, cone$_o V$ is apexed (at
$o$) iff there exists a plane $H'$ strictly supporting cone$_oV$
at $o$.

\textbf{Proposition 4.1:} \emph{Assume $V \subset S(o)$. Then}
\begin{itemize}
\item[(i)] \emph{if} cone$_oV$ \emph{is apexed $($at} $o$), \emph{then} cone$_oV \cap S(o)$ \emph{is contained in a small cap
of} $S(o)$, \item[(ii)] cone$_oV \cap S(0)$ \emph{is spherically
convex in} $S(o)$, \item[(iii)] sph-conv$(V) = {\rm cone}_o V \cap
S(o)$.
\end{itemize}

\textbf{Proof:}
\begin{itemize}
\item[(i)] This follows readily from Definition 4.2.
\item[(ii)] The cases $V = \emptyset$ and $\# V=1$ are trivial. In the case $\# V=2$ and if the two points of
$V$ are antipodal, cone$_oV \cap S(o)=V$ and the set of two antipodal points is spherically convex (cf. Examples 4.1 b) above).
Let $a,b \in {\rm cone}_o V \cap S(o)$ be two non-antipodal points,
and let $z$ be a point on the short arc between $a$ and $b$ on the
unit circle centered at $o$ passing through $a$ and $b$. Then $z \in {\rm cone}_o (\{a,b\})$, i.e., $z = \lambda
a + \mu b$ for some $\lambda, \mu \ge 0$. The point $a$, resp. $b$, is a non-negative combination of a finite subset
$\{a_1, \dots, a_m\} \subset V$, resp. $\{b_1, \dots, b_l\} \subset V$, i.e., $a = \sum \limits^m_{i=1} \lambda_i a_i$,
resp. $b = \sum \limits^l_{i=1} \mu_i b_i$, for some $\lambda_i \ge 0 \, (1 \le i \le m)$, resp. $\mu_i \ge 0 \, (1 \le i
\le l)$. Thus
\[
z =  \lambda a + \mu b  =  \lambda \sum \limits^m_{i=1} \lambda_i
a_i + \mu \sum \limits^l_{i=1} \mu_i b_i
 =  \sum \limits^m_{i=1} \lambda \lambda_i a_i + \sum \limits^l_{i=1} \mu \mu_i b_i\,,
\]
which is a non-negative combination of $\{a_i : 1 \le i \le m\}
\cup \{b_i : 1 \le i \le l\} \subset V$, hence $z \in {\rm cone}_o
V$. Since $\|z\|=1 \,,\, z \in {\rm cone}_o V \cap S (o)$.
\item[(iii)] By (ii), $\mbox{sph-conv} V\subset {\rm cone}_o V \cap
S(o)$. To prove ``$\supset$'' assume that $F$ is a spherically
convex subset of $S(o)$ containing $V$. We show that cone$_o V
\cap S(o) \subset F$. Let $z \in {\rm cone}_o V \cap S(o)$.
Then $z$ is a non-negative combination of a finite subset $W
\subset V$. By Carath\'{e}odory's Theorem ($2$-dimensional case)
there are three points $w_1, w_2, w_3 \in W$ no two of which are antipodal such that
\begin{equation}\label{sternsternstern}
z = \lambda_1 w_1 + \lambda_2 w_2 + \lambda_3 w_3\,,
\end{equation}
where $\lambda_i \ge 0 \,\, (1 \le i \le 3)$, i.e., $z \in {\rm cone}_o(\{w_1,w_2,w_3\})$.

If there is a $\lambda_i \, (1 \le i \le 3)$ such that $\lambda_i = 0$, say $\lambda_3 = 0$, then $z \in
{\rm cone}_o(\{w_1,w_2\})$, i.e., $z$ lies on the short arc connecting $w_1$ and $w_2$ of the unit circle with
center $o$ passing through $w_1$ and $w_2$. Since $w_1,w_2 \in F$ and $F$ is closed to such short circles
$z \in F$.

Assume now that $\lambda_1,\lambda_2,\lambda_3 > 0$. The vector
$\lambda_1 w_1 + \lambda_2 w_2$ is a scalar product of the unit
vector $w =_{\rm def} \frac{\lambda_1w_1 +
\lambda_2w_2}{\|\lambda_1w_1 + \lambda_2w_2\|}$ by the positive
scalar $\lambda =_{\rm def} \|\lambda_1 w_1 + \lambda_2 w_2 \|$.
The unit vector $w$ lies on the short arc with endpoints $w_1,
w_2$ of the unit circle centered at $o$ and passing through
$w_1$ and $w_2$, hence $w \in F$. By (\ref{sternsternstern}) we
have $z = \lambda w + \lambda_3 w_3$, hence $z$ lies on the short
arc with endpoints $w$ and $w_3$ of the unit circle centered at
$o$ and passing through $w$ and $w_3$, hence $z \in F$. This
proves (iii). \rule{2mm}{2mm}
\end{itemize}

From Proposition 4.1 (iii) we conclude

\textbf{Corollary 4.1:} \emph{If $V \subset S(o)$ is spherically convex, then $V = \mbox{cone}_o V \cap S(o)$}.

Theorem 4.1 below is stated just for a general background without proof, since it is not needed in the
sequel.

\textbf{Theorem 4.1:} \emph{If $F \subset S(o)$ is a closed, spherically convex set
not containing a pair of antipodal points, then}
\begin{enumerate}
\item[a)] $F$ \emph{is contained
in a small cap of} $S(o)$,
\item[b)] \emph{$F$ is the intersection of a family of closed hemispheres of $S(o)$}.
\end{enumerate}

A close connection between spindle convexity and spherical convexity is enunciated in the following statement
(compare \cite{Bezdek}, Lemma 5.6).

\textbf{Proposition 4.2:} \emph{If $K \subset {\mathbb R}^3$ is
spindle convex, then $S(o) \cap K$ is spherically convex on
$S(o)$, and either $K = \bR^3$ or $K = B(o)$ or $S(o) \cap K$
is contained in small cap of $S(o)$}.

\textbf{Proof:} If $\# (S(o) \cap K \le 1$, there is nothing to prove. Assume
$\# (S(o) \cap K) \ge 2$,
let $a,b \in F =_{\rm def} S(o) \cap K, a \not= b$, and let arc$(a,b)$
be a short arc (possibly half a circle)
of radius 1 on $S(o)$ with endpoints $a,b$. Since $K$ is $s$-convex, arc$(a,b) \subset K$, hence
arc$(a,b) \subset F$ and $F$ is spherically convex. By Remark 3.2, 1) either $K = \bR^3$ or $K$ is
the intersection of a non-empty family of unit balls. In the second case $K \subset B(p)$ for some
$p \in \bR^3$ (the unit ball centered at $p$). If $K =B(o)$, we are done. Otherwise $K \subset B(p)$
for some $p \not= o$, hence $S(o) \cap K \subset S(o) \cap B(p)$, and since $p \not= o$, clearly
$S(o) \cap B(p)$ is a small cap of $S(o)$. \hfill \rule{2mm}{2mm}

\textbf{Remark:} The intersection $S(o) \cap K$ alluded to above need not be \emph{strictly}
spherically convex. Let $a,b,c \in S(o)$ be three points no two of which are antipodal such that $o \notin
{\rm aff} \{a,b,c\}$. Then $\triangle =_{\rm def}$ sph-conv$(\{a,b,c\}$) is a spherical
triangle which is not strictly spherically convex (the edges of $\triangle$ are geodetic arcs of $S(o)$).
Put $K=_{\rm def} \cap \{B \subset \bR^3:B$ is a unit ball and $\triangle \subset B\}$. $K$ is spindle convex
(by Theorem 3.1), and it can be shown that $S(o) \cap K = \triangle$, thus it is not strictly spherically convex.

\section{The face structure of Ball polytopes}

\textbf{Essential and inessential points:}

Let $V \not= \emptyset$ be a finite set of points in ${\mathbb
R}^3$. We assume that cr$(V) < 1$; ipso facto ${\cal B}(V)$ is
full-dimensional, and diam$V < 2$.

\textbf{Definition 5.1:} A point $v \in V$ is \emph{essential} if
${\cal B}(V) \varsubsetneqq {\cal B}(V \setminus \{v\})$,
otherwise (if ${\cal B}(V) = {\cal B}(V \setminus \{v\})$) it is \emph{inessential}.

Note that $v$ is essential iff there is a unit ball that includes
$V \setminus \{v\}$ and misses $v$. In particular, if $v$ is
essential, then $v$ is an extreme point of $V$. Moreover, this
characterization of essential points shows that if $v$ is
essential in $V$, then $v$ is essential in every subset of $V$
that contains $v$.
Denote by ess$(V)$ the set of essential points of $V$. Soon we
shall show that ${\cal B} ({\rm ess} (V)) = {\cal B}(V)$ (Theorem 5.1
below).

\textbf{Lemma 5.1:} \emph{We have $v \in {\rm ess} (V)$ iff there
is a point $x \in {\mathbb R}^3$ such that $\|x - v \| = 1$ and
$\|x-w\|< 1$ for $w \in V \setminus \{v\}$}\footnote{See \cite{Bez-Nasz}, Proposition 1.1,
for a similar statement.}.

\textbf{Proof:} ($\Rightarrow$) Assume $v \in {\rm ess} (V)$. Let
$c$ be the center of a unit ball that includes $V \setminus \{v\}$
but misses $v$, and let $b$ be an interior point of ${\cal B}(V)$. Then
$\|b-w\|< 1$ for $w \in V, \,\, \|c-v\|> 1$ and $\|c-w\|\le 1$ for
$w \in V \setminus \{v\}$. For $0 \le \vartheta \le 1$, define
$a(\vartheta)=_{\rm def} (1-\vartheta) b + \vartheta c$. Choose
$\vartheta, \, 0 < \vartheta < 1$, such that
$\|a(\vartheta)-v\|=1$, and define $x = a (\vartheta)$. Then
$\|x-w\|<1$ for $w \in V \setminus \{v\}$.

$(\Leftarrow)$ Assume $\|x-v\|=1$ and $\|x-w\|<1$ for $w \in V
\setminus \{v\}$. Choose a positive number $\varepsilon$ such that
$\|x-w\|\le 1-\varepsilon$ for $w \in V \setminus \{v\}$, and
define $x' =_{\rm def} x + \varepsilon (x-v)$. Then $\|x'-v\| = 1
+ \varepsilon > 1$, but $\|x'-w\|\le 1$ for $w \in V \setminus
\{v\}$, hence $v \in {\rm ess} (V)$. \hfill \rule{2mm}{2mm}

Next we show that the removal of an inessential point does not
change the status (essential/inessential) of the remaining points
with respect to the remaining set. We have seen already that if $w
\in V \setminus \{v\}$  is essential in $V$, then it is essential
in $V \setminus \{v\}$. The next lemma does the rest.

\textbf{Lemma 5.2:} \emph{Assume $V \subset {\mathbb R}^3$ is
finite, with} cr$(V) < 1$. \emph{If $v$ and $v'$ are inessential
in $V$, then $v'$ is inessential in $V \setminus \{v\}$}.

\textbf{Proof:} Assume, on the contrary, that $v' \in {\rm
ess}(V\setminus \{v\})$. By Lemma 5.1 there is a point $x \in
{\mathbb R}^3$ such that $\|x-v'\|=1$, and $\|x-w\|< 1$ for $w \in
V \setminus \{v,v'\}$. What about $\|x-v\|$? If $\|x-v\| < 1$,
then $v'$ is essential in $V$, contrary to our hypothesis. If
$\|x-v\|> 1$, define $x(\varepsilon)=_{\rm def} (1-\varepsilon) x
+ \varepsilon v'$. For $\varepsilon > 0$ sufficiently small, we
find that $\|x (\varepsilon) - v' \| < 1\,, \, \|x (\varepsilon) -
v \|> 1$, and $\| x (\varepsilon) - w \| < 1$ for all $w \in V
\setminus \{v,v'\}$. Thus $v \in {\rm ess} (V)$, contrary to our
hypothesis. If $\|x-v\| = 1$, then the points $x,v,v'$ form an
isosceles triangle with $\|x-v\| = \| x-v'\| = 1$. Move the point
$x$ slightly parallel to the base $[v,v']$ of this triangle. The
resulting point $x'$ satisfies: $\|x'-v\| > 1\,, \, \|x' - v'\| <
1$ (or vice versa), and $\|x'-w\| < 1$ for all $w \in V \setminus
\{v,v'\}$. Thus $v \in {\rm ess} (V)$ (or $v' \in {\rm ess} (V)$),
contrary to our hypothesis. \hfill \rule{2mm}{2mm}

\textbf{Theorem 5.1:} \emph{Assume that $V \subset {\mathbb R}^3$
is finite with} cr$(V) < 1$. \emph{Then ${\cal B}(V) = {\cal
B}({\rm ess}(V))$}.

\textbf{Proof:} Remove the inessential points of $V$ one by one,
and apply Lemma 5.2 each time to ensure that all inessential
points remain inessential. \hfill \rule{2mm}{2mm}

Clearly, for any proper subset $W \varsubsetneq {\rm ess}(V)$ one
has ${\cal B}(W) \supsetneqq {\cal B}(V)$ (=${\cal B} ({\rm ess}(V)$), hence, by
Theorem 5.1, ess$(V)$ is the smallest subset $X$ of $V$ satisfying
${\cal B}(X) = {\cal B}(V)$. This brings us to

\textbf{Definition 5.2:} A finite set $V \subset {\mathbb R}^3$
satisfying cr$(V) < 1$ and $V = {\rm ess}(V)$ is \emph{tight}\footnote{In \cite{Bezdek},
Definition 6.3, such a set is called \emph{reduced}, and in \cite{Bez-Nasz}, p. 258, it is called
\emph{non-redundant}. Sallee in \cite{Sa} had no name for this notion, though he clearly
had it in mind in his (confused) definition of \emph{a finite intersection of balls} in p. 316, lines
1-5.}.

If $V \subset {\mathbb R}^3$ is finite and not tight, then after removal of the inessential points from $V$
there remains a (tight) set whose ball polytope is the same
${\cal B}(V)$. Hence, in addition to the standard assumption
cr$(V)<1$, we will also assume
that $V$ is tight when dealing with the face structure of the Ball
polytope ${\cal B}(V)$ of $V$. Theorem 5.2 a), b) below guarantee these two properties in case that
$V \subset {\mathbb R}^3$ is extremal for the V\'{a}zsonyi
problem.

\textbf{Theorem 5.2:} \emph{Assume that $V \subset {\mathbb R}^3$
is finite and} diam$V=1$. \emph{Then}
\begin{enumerate}
\item[a)] cr$(V) < 1$.
\item[b)] \emph{If a point $v \in V$ is
incident with (at least) two diameters of $V$, then} $v \in {\rm
ess} (V)$.
\item[c)] \emph{ If $V$ is extremal for the V\'{a}zsonyi problem, then $V$ is tight}.
\end{enumerate}

\newpage

\textbf{Proof:}
\vspace*{-0.4cm}
\begin{enumerate}
\item[a)] In fact, we have cr$(V) \le {\sqrt{\frac{3}8}} \,  (<
1)$, where equality is attained iff $V$ is the set of vertices of
a tetrahedron of edge-length 1 (Jung's Theorem; see
\cite{Egg}, Theorem 49 in p. 111, or \cite{Jung}).
\item[b)] Assume $a,b \in V\,, \, \|a-v\| =
\|b-v\|=1$. Define $c =_{\rm def} \frac{1}{2} (a+b)$. Denote by
$H$ the plane through $c$ perpendicular to $[a,b]$ and by
$H_a,H_b$ the closed half-spaces determined by $H$ that contain
$a$ and $b$, respectively. Note that
\begin{equation}\label{stern}
V \subset {\cal B}(V) \subset B(a,1) \cap B(b,1) = (B(a,1) \cap
H_b) \cup (B(b,1)\cap H_a)\,.
\end{equation}
Denote by $C$ the circle with center $c$ that passes through $v$,
i.e., $C = \{x \in {\mathbb R}^3 : \| x-a \| = \| x-b \| = 1\}$.
Denote by $w$ the point at distance $1$ from $v$ on the ray
$\overrightarrow{vc}$. Note that $\|w-x\| < 1$ for every point $x$
on $C$ or inside the disc conv$C$, except $v$, i.e., conv$C
\setminus \{v\} \subset {\rm int} B(w,1)$. By spindle convexity of
$B (w,1)$ (Lemma 3.1) we have
\begin{equation}\label{sternstern}
(B(a,1) \cap H_b) \cup (B(b,1) \cap H_a) \subset \{v\} \cup
\mbox{int}B (w,1)\,.
\end{equation}
Thus $\|v-w\|=1$, and by (\ref{stern}) and (\ref{sternstern}) we
have $\|v'-w\| < 1$ for $v' \in V \setminus \{v\}$. This shows, by
Lemma 5.1, that $v$ is essential in $V$.
\item[c)] This follows from a), b) and Proposition-Definition 2.1 above.
\hfill \rule{2mm}{2mm}
\end{enumerate}

Now we are ready to introduce the spherical face structure ${\cal
S} {\cal F} ({\cal B} (V))$ of the Ball polytope ${\cal B} (V)$.
We assume that $V$ is finite, nonempty with cr$(V) < 1$, and tight $(V
= {\rm ess} (V))$. ${\cal S} {\cal F} ({\cal B}(V))$ will be a
finite $5$-tiered poset (not necessarily a lattice) consisting of
the empty set $\emptyset$, vertices, edges, facets, and ${\cal
B}(V)$ itself. We shall introduce and discuss these elements in
the following order: First facets, then vertices, then edges.

\textbf{Facets:}

\textbf{Definition 5.3 (facet of \boldmath${\cal B}(V)$\unboldmath):}
For a point $p \in V$ the set $
F_p =_{\rm def} \{x \in {\cal B}(V) : \| x-p \| = 1\}$
is a \emph{facet} of ${\cal B}(V)$.

By definition, $F_p$ is a compact subset of the unit $2$-sphere
centered at $p$, $S_p =_{\rm def} \{x \in {\mathbb R}^3 : \| x - p
\| = 1\}$. Moreover, since $V$ is tight, $p$ is essential in $V$,
and thus (by Lemma 5.1) there is a point $x \in S_p$ such that $\|
x - q \| < 1$ for $q \in V \setminus \{p\}$. The point $x$ is
interior to $F_p$, relative to $S_p$. Hence relint$F_p \not=
\emptyset$ and ${\cal B}(V)$ has exactly $\# V$ facets. For $V =
\{p\}$ we have, of course, $F_p = S_p$. When $V \supsetneqq
\{p\}$, we can write
$F_p : = S_p \cap \bigcap \limits_{q \in V \setminus \{p\}} \{x \in
{\mathbb R}^3 : \| x-q \| \le \|x-p \|\}$.
For $x \in {\mathbb R}^3$ and $q \in V \setminus \{p\}$ we have
$\|x-q \| \le \| x-p \| \Leftrightarrow \langle x-p, q-p \rangle
\ge \frac{1}{2} \| q-p \|^2$. (This can be verified by a
straightforward calculation.) The set
$H_{p,q} =_{\rm def} \{x \in {\mathbb R}^3 : \langle x-p, q-p
\rangle \ge {\tiny\frac{1}{2}} \|q - p \|^2 \}$
is a closed half-space of ${\mathbb R}^3$ that misses $p$. The
intersection of $S_p$ with such a half-space is a small cap of
$S_p$. Thus we have expressed $F_p$ as an intersection of $\# V -
1$ or less small caps of $S_p$:
$F_p = \bigcap \limits_{q \in V \setminus \{p\}} (S_p \cap H_{p,q})$.
A small cap $C =_{\rm def} S_p \cap H_{p,q}$ is strictly
spherically convex (in the sense of Definition 4.1
above), hence we have

\textbf{Corollary 5.2:} \emph{If $V \subset {\mathbb R}^3$ is
finite, tight, $\# V \ge 2$, and {\rm cr}$(V) < 1$, then for $p
\in V$ the facet $F_p$ is a strictly spherically convex subset of
$S_p$, contained in a small cap of $S_p$, compact with a non-empty relative interior}
$F_p = \{x \in S_p : \|x-q\| \le 1$ for $q \in V \setminus \{p\}\} = S_p \cap (\cB (V \setminus \{p\})$ and
relint$F_p = \{x \in S_p : \|x-q\|<1$
\emph{for} $ q \in V \setminus \{p\}\}$.

Next we show that $F_p$ is topologically a $2$-disc
(= $2$-cell).

\textbf{Theorem 5.3:} \emph{If $F \subset S(o)$ is compact,
strictly spherically convex and contained in a small cap of $S(o)$ with} relint$F \not= \emptyset$,
\emph{then there is a homeomorphism of $F$ to a closed hemisphere
of $S(o)$, which maps} relbd$F$ \emph{$($the relative boundary of
$F$$)$ to a great circle in} $S(o)$.

\textbf{Proof:} Assume, without loss of generality, that the north
pole $N(0,0,1)$ is relatively interior to $F$. We notice that for $(\varphi,
\vartheta) \in (-\pi,\pi] \times [-\frac{\pi}{2},\frac{\pi}{2}]$
the vector $(\cos \varphi \cos \vartheta, \sin \varphi
\cos \vartheta, \sin \vartheta)$ is a parametric representation of the meridian
line of $S(o)$ of fixed longitude $\varphi \in (-\pi,\pi]$ (spherical
coordinates).

Define
$\vartheta(\varphi) =_{\rm def} \min \left\{\vartheta: \left(\begin{array}{c}
\cos \varphi \cos \vartheta\\ \sin \varphi \cos \vartheta\\
\sin \vartheta \end{array}\right) \in F \right\}$
for $-\pi < \varphi \le \pi$.
Clearly, $-\frac{\pi}{2} < \vartheta (\varphi) < \frac{\pi}{2}$, and
one can easily check that $\vartheta(\varphi)$
is a continuous function of $\varphi$ (this follows from the compactness and the strict spherical convexity of $F$).
A general point of $F$ is of the form $x = (\cos \varphi \cos \vartheta, \sin \varphi \cos \vartheta,
\sin \vartheta)$, $\vartheta (\varphi) \le \vartheta \le \frac{\pi}{2}$, and
$\vartheta = (1-\lambda) \vartheta (\varphi) + \lambda \cdot
\frac{\pi}{2}\,,\,\, 0 \le \lambda \le 1$.
Map $x$ to
$f(x) =_{\rm def} (\cos \varphi \cos (\lambda \cdot \frac{\pi}{2}), \sin \varphi \cos (\lambda \cdot \frac{\pi}{2}),
\sin (\lambda \cdot \frac{\pi}{2}))$.

This defines a homeomorphism between $F$ and the northern
hemisphere of $S(o)$, which maps relbd$F$ to the equator. \hfill
\rule{2mm}{2mm}

A point $x \in {\cal B}(V)$ is a boundary point of ${\cal B}(V)$ if it
satisfies at least one of the inequalities $\|x-p\| \le 1 \,\, (p
\in V)$ as an equality. Thus
\begin{equation}\label{aa}
{\rm bd}{\cal B}(V) = \bigcup \{F_p: p \in V\}\,.
\end{equation}

\textbf{Theorem 5.4:} \emph{For} $p \in V$, relbd$F_p = \bigcup \{F_p \cap F_q : q \in V \setminus
\{p\}\}$.

\textbf{Proof:} Assume $x \in F_p$. Then $\|x-p\| = 1$ and
$\|x-q\| \le 1$ for $q \in V \setminus \{p\}$. If all these
inequalities are strict, then they hold for $x' \in S_p$ instead
of $x$, provided
$\|x'-x\| \le 1 - \|x-q\| \, \mbox{ for } \, q \in V \setminus
\{p\}$.

Thus $x$ is interior to $F_p$ relative to $S_p$. If, on the other
hand, $\|x-q \| = 1$ for some $q \in V \setminus \{p\}$, then the
points $x, p, q$ form an isosceles triangle $T$ with
$\|x-p\|=\|x-q\|=1$. (We have $\|p-q\|<2$, since otherwise cr$(V)
\ge 1$, contrary to our assumptions.) Moving $x$ slightly within
aff$T$ away from $q$ along a unit circle centered at $p$, we
obtain a point $x' \in S_p$ that satisfies $\|x'-q\| > 1$. Hence
$x' \notin {\cal B}(V)$ and, in particular, $x' \notin F_p$. Thus $x
\in
\partial F_p$. \hfill \rule{2mm}{2mm}

What can be said about the intersection $F_p \cap F_q$ for $p,q
\in V, p \not= q$? By definition
$F_p \cap F_q = \{x \in {\cal B}(V): \|x-p\| = \|x-q\| = 1\}$.

The set
$C_{pq} =_{\rm def} \{x \in {\mathbb R}^3: \|x-p\| = \|x-q\|=\}$
is the ``equatorial'' circle of the lentil shaped body ${\cal
B}(\{p,q\})$, of radius $r_{pq}$ satisfying $0 < r_{pq} =_{\rm
def} \sqrt{1-\frac{1}{4} \|p-q\|^2} < 1$ and centered at $c_{pq}
=_{\rm def} \frac{1}{2} (p+q)$, lying on the plane that bisects
the segments $[p,q]$ perpendicularly. The intersection of $C_{pq}$
with any unit ball $B(v,1)$, $v \in V \setminus \{p,q\}$,
is a-priori either the whole of $C_{pq}$, or a closed arc of
$C_{pq}$, or a singleton or empty. It follows that $F_p \cap
F_q = C_{pq} \cap {\rm bd}{\cal B}(V)$ has at most $n-2$
components (closed arcs or singletons).
Next we prove (part (i) of the following proposition) that
for $\# V \ge 3$ the intersection $F_p \cap F_q$ is not the whole of $C_{pq}$, not a
singleton, nor empty, which is a key for some later developments.

\textbf{Proposition 5.1:}
 \emph{Assume that $V \subset \bR^3$, $\# V
\ge 3$ is finite, tight, cr$(V) < 1$, and let $p,q \in V$. Then}
\begin{itemize}
\item[(i)] $C_{pq}$ \emph{is a circle of radius $0
< r_{pq} < 1$, and for $v \in V \setminus \{p,q\}$, $C_{pq}
\cap B(v,1)$ is a relatively closed arc of $C_{pq}$ $($with two
$($different$)$ endpoints; ipso facto it has a non-empty relative
interior$)$. Equivalently: ${\cal B}(v,1)$ misses a nonempty,
proper (relatively open) arc of $C_{pq}$}.
\item[(ii)] $F_p \cap F_q = C_{pq} \cap {\cal
B}(V)$ \emph{is either empty, or it is a finite union of at most
$\# V-2$ closed arcs $($connectivity components$)$ of $C_{pq}$, with
each arc having either two endpoints or being a single point}.
\item[(iii)] \emph{Each
endpoint of a connectivity component of $C_{pq} \cap {\cal B} (V)$ is
common to at least} three \emph{facets of ${\cal B}(V)$ $($including
$F_p$ and $F_q$, of course$)$}.
\item[(iv)] \emph{An isolated point of $C_{pq} \cap {\cal B}(V)$ is common to at least} four \emph{facets
of ${\cal B}(V)$}.
\end{itemize}

\textbf{Proof:}
\begin{itemize}
\item[(i)] Since $V$ is tight and cr$(V) < 1$, ${\cal B} (\{p,q\})
= {\cal B} (p) \cap {\cal B}(q)$ is a full-dimensional ``lentil''
shaped body whose ``equator'' $C_{pq}$ is a circle of radius
$r_{pq}$, $0 < r_{pq} < 1$. It remains
to show that, for $v \in V \setminus
\{p,q\}$, $C_{pq} \setminus B(v) \not= \emptyset$ and that
$B(v) \cap C_{pq}$ has more than one point.
(Note that $B(v,1)$ is closed, hence $C_{pq} \setminus B(v,1)$ is
relatively open in $C_{pq}$; so it cannot be a singleton.)

Assume, by r.a.a., that $C_{pq} \subset \cB(v)$. The boundary bd$(B(p) \cap
B(q))$ of $B(p) \cap B(q)$ consists of two small caps cap$(p) =_{\rm def} S(p) \cap B(q)$ of $S(p)$ and
cap$(q) =_{\rm def} S(q) \cap B(p)$ of $S(q)$. Let $x \in \mbox{ relint cap } (p)$ (= cap$(p) \setminus C_{pq}$) and
let $C \subset S(p)$ be a great circle of $S(p)$ (of radius $1$ and centered at $p$) passing through $x$.
The intersection $C \cap {\rm cap} (p)$ is a short arc of $C$, both of its endpoints lie on $C_{pq} \subset B(v)$.
By Lemma 3.1' (simplified Sallee Lemma) relint $A \subset {\rm int} B(v)$, unless $v=p$, which is impossible ($v \not= p$ by
assumption). Thus $x \in {\rm int} B(v)$ and cap$(p) \subset B(v)$. A similar argument shows that cap$(q) \subset B(v)$,
hence bd$(B(p) \cap B(q)) = {\rm cap}(p) \cup
{\rm cap} (q) \subset B(v)$. It follows that $B(p) \cap B(q) \subset B(v)$ (by convexity of $B(p) \cap B(q)$), a contradiction to
the assumption that $v$ is essential. Thus
\begin{equation}\label{sterna}
C_{pq} \nsubseteq B(v)\,,
\end{equation}
and similarly (interchanging the roles of $q$ and $v$) we get
\begin{equation}\label{sternb}
C_{pv} \nsubseteq B(q)\,.
\end{equation}
By (\ref{sterna}) $C_{pq} \cap B(v)$ is either a proper arc of $C_{pq}$ or a singleton or empty. In order to
disprove the latter two possibilities we consider two cases.

\textbf{Case I:} $v \in H =_{\rm def} {\rm aff} C_{pq}$.\\
If $\|v-\frac{1}{2} (p+q)\| \le r_{pq}\, (< 1)$, then the circle $B(v) \cap H$ of radius $1$ has its center
$v$ in conv$\,C{pq}$, it is not contained in conv$\,C_{pq}$ (since $r_{pq} \le 1$), and it does not contain $C_{pq}$ (by (\ref{sterna})),
hence it intersects $C_{pq}$ in two points. Assume now that $\|v-\frac{1}{2} (p+q)\| > r_{pq}$, i.e., $v$ lies out of the circle
$C_{pq}$ (in the plane $H$).
Then
\begin{equation}\label{6}
{\rm dist}(v,B (p) \cap B(q))= {\rm dist}(v,C_{pq})\,.
\end{equation}
Since cr$\{v,p,q\}< 1$, it follows from (\ref{6})
that dist$(v,C_{pq}) < 1$, and by (\ref{sterna}) the circle $B(v) \cap H$ intersects $C_{pq}$ in two points.

\textbf{Case II:} $v \notin H$.\\
Let $H^+, H^-$ be the two open half-spaces bounded by $H$ and assume, without loss of generality, that $\{v, q\} \subset H^-$ and $p \in H^+$.

If $C_{pq} \cap C_{pv}$ is  singleton or empty, then -- both circles lying on $S(p)$ -- either $C_{pv} \subset H \cup H^+$ or $C_{pv} \subset H \cup H^-$.
In the first case $C_{pq} \subset B(v)$, contradicting (\ref{sterna}), and in the second case $C_{pv} \subset B(q)$, contradicting
(\ref{sternb}). Thus $\# (C_{pq} \cap C_{pv}) \ge 2$, i.e., $C_{pq} \cap B(v)$ is a proper, non-singleton arc of $C_{pq}$.

\item[(ii)] By part (i), $F_p \cap F_q = C_{pq} \cap {\cal B}(V) = \bigcap \{C_{pq}
\cap B(v,1) : v \in V \setminus \{p,q\}\}$ is a finite
intersection of $\# V-2$ closed arcs on $C_{pq}$, which is either empty or
of the type described in the proposition.
\item[(iii)] Let $x$ be an endpoint of a connectivity component
$A$ of $C_{pq} \cap {\cal B}(V)$. The point $x$ is common to $F_p$
and $F_q$, and if it does not belong to any other facet $F_v$, $v
\in V \setminus \{p,q\}$, then $x \in {\rm int} {\cal B}(V
\setminus \{p,q\})$. It follows that there is a whole neighborhood
of $x$ on $C_{pq}$ included in ${\cal B} (V)$, contradicting the
assumption that $x$ is an endpoint of a connectivity component of
$C_{pq} \cap {\cal B}(V)$.
\item[(iv)] Follows from (i) and (iii). \hfill \rule{2mm}{2mm}
\end{itemize}

The following example shows that the intersection of two facets $F_p \cap F_q = C_{pq} \cap
{\cal B}(V)$ can be, indeed, \emph{any} finite disjoint union of closed (not necessarily short) arcs of $C_{pq}$,
including singletons.
In the next paragraph we define the $1$-skeleton of a ball polytope and refer to this as a basic example of a ball polytope whose
$1$-skeleton is not $3$-connected.

\textbf{Example 5.1:} Let $p = p (0,0,h), \, q = q (0,0,-h), \, 0
< h < 1$, and let $C =_{\rm def} S(p) \cap S(q) = C_{pq}$ be a circle in the
$(x,y)$-plane centered in the origin of radius $\sqrt{1-h^2}$, $0 <
\sqrt{1-h^2} < 1$. The set bd${\cal B}(\{p,q\})$ is the
union of two small caps of $S(p)$ and $S(q)$, respectively,
the relative boundary of both of which is $C$. Assume that $C$ is positively oriented
(counterclockwise in all our figures).

\textbf{Notation:} For $u,v \in C$, $u \not= v$, denote by orarc$(u,v)$ the
positively oriented arc of $C$ beginning in $u$ and terminating in
$v$. (orarc$(u,v)$ may be a short as well as a long arc of $C$, of
course).
\vspace{1cm}

\begin{picture}(280,130)
\put(50,0){\includegraphics[scale=0.5]{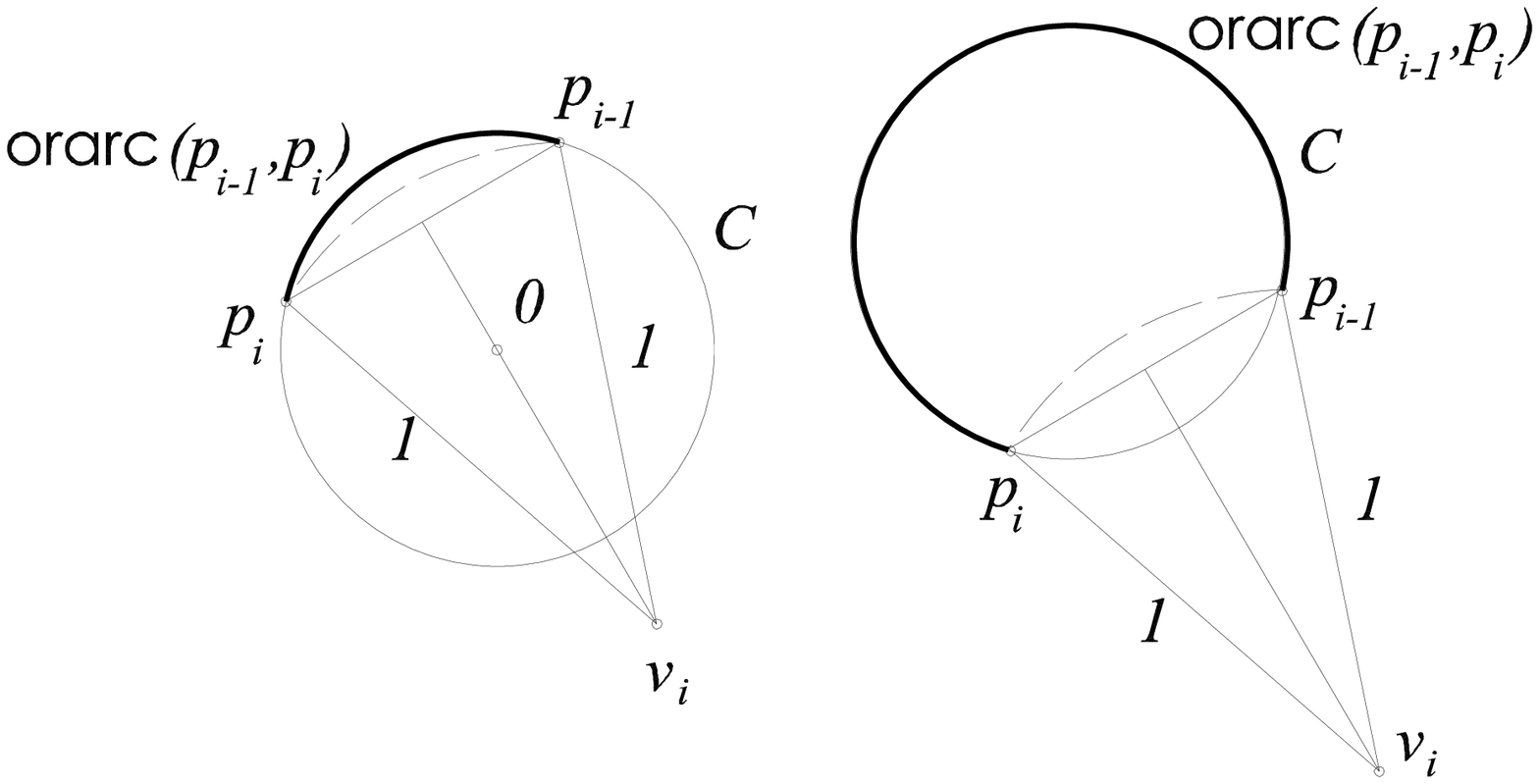}} 
\end{picture}
\vspace{-1cm}

\begin{center}
-- Figures 3 and 4 --
\end{center}

Mark succesively $n$ points $p_0, p_1, \dots, p_{n-1} \,, n \ge
2$, on $C$ so that the cyclic order of the indices $0,1,\dots,n-1$
is compatible with the positive orientation on $C$, i.e., for $1 \le i
\le n$ the arc orarc $(p_{i-1},p_i$) (the indices are taken modulo
$n$, so that $p_n = p_0$) does not contain any point $p_j \, (1
\le j \le n)$ in its relative interior.
For $1 \le i \le n$ choose a point $v_i$ in the $(x,y)$-plane such that
\begin{itemize}
\item[i)] $\|v_i - p_{i-1}\| = \|v_i - p_i\| = 1$, and
\item[(ii)] $\|v_i - x\| > 1$ for $x \in$ relint(orarc$(p_{i-1}, p_i)$).
\end{itemize}
Conditions (i) and (ii) uniquely determine $v_i$, and they imply
that $v_i$ and orarc$(p_{i-1}, p_i$) lie on opposite sides of the
chord $[p_{i-1},p_i]$ of $C$; see Figures 3 and 4.

Thus the unit ball $B(v_i)$ misses relint(orarc$(p_{i-1},p_i)$) and includes its complementary arc
on $C$, namely orarc$(p_i, p_{i-1}$).

Choose any subset $W \subseteq \{v_0,v_1, \dots, v_{n-1}\}$ and
define $V=_{\rm def} \{p,q\} \cup W$. Then cr$(V)<1$, $V$ is tight, and
\begin{equation}\label{5}
\begin{array}{l}
F_p \cap F_q = C_{pq} \cap {\cal B}(V)=
\left(\bigcup\{\mbox{orarc} (p_{i-1},p_i):\right.\\
~\\
\left.v_i \notin W, \, 1 \le i \le n\}\right) \cup \{p_i: v_i,v_{i+1} \in W, \,
1 \le i \le n\}\,.
\end{array}
\end{equation}

Hence, by choosing $W$ at will, we have full control of which arc
orarc$(p_{i-1},p_i)$ of $C$ belongs to $F_p \cap F_q = C \cap
{\cal B}(V)$, and of which not. E.g., if $W = \{v_0,v_1, \dots,
v_{n-1}\}$, then $F_p \cap F_q = \{p_0,p_1,\dots,p_{n-1}\}$
because $B(v_i)$ misses relint(orarc$(p_{i-1},p_i)$) for $1 \le i
\le n$.

Note that if $v_i \in W \subset V$, then the facet $F_{v_i}$ of
${\cal B}(V)$ has exactly two ``edges'', both of which are short arcs with
endpoints $p_{i-1},p_i$ of the circles $C_{pv_i} = S(p) \cap
S(v_i)$ and $C_{qv_i} = S(q) \cap S( v_i)$, respectively. Thus
$F_{v_i}$ is a ``digonal'' facet of ${\cal B}(V)$.

The notion of
``edge'' will be defined in the next paragraph, and Example 5.1 will
be generalized in Example 6.2 below.

\section{The face structure of Ball polytopes -- continuation}

\textbf{Vertices:}

Next we turn our attention to the \emph{vertices} of ${\cal
B}(V)$. As in the previous paragraph we assume (as usual) that $V
\subset {\mathbb R}^3$ is finite, non-empty, tight and satisfies cr$(V) < 1$.

\textbf{Definition 6.1 (a vertex of \boldmath${\cal
B}(V)$\unboldmath:} A boundary point $z$ of ${\cal B}(V)$ is a
\emph{vertex} of ${\cal B}(V)$ if either $z$ belongs to three or
more distinct facets of ${\cal B}(V)$, in which case $z$ is a
\emph{principal vertex}, or $z \in V \cap {\cal B}(V)$ and $z$
belongs to exactly two facets of ${\cal B}(V)$, in which case $z$
is a \emph{dangling vertex}. Denote by vert${\cal B}(V)$ the set
of vertices of ${\cal B}(V)$. In other words, $z \in {\rm
vert}{\cal B}(V)$ if and only if $z \in {\cal B}(V)$ and $\| z - p
\| = 1$ holds for at least three points $p \in V$, or if $z \in V
\cap {\cal B} (V)$ and $\|z-p\| = 1$ holds for exactly two points
$p \in V$\footnote{As will become apparent below, a dangling vertex $z\, (\in V)$ always
belongs to the relative interior of the intersection $F_p \cap F_q$ of the two facets $F_p$ and $F_q$,
\emph{unlike} the case of $3$-polytopes, where every vertex is a
``corner'' of the body. Nevertheless, this apparently strange notion of ``dangling vertex'' proves
to the fruitful. E.g., in the duality theory developed in \S~8 below for every extremal $V$, $\cB(V)$ is
self-dual without any artificial restrictions as given in \cite{Bez-Nasz}, p. 260, lines 19-22 (a restriction
that $\cB(V)$ has no digonal facets). A similar notion of dangling vertex can also relieve the Bezdek-Nasz\'odi duality theorem
in \cite{Bez-Nasz}, p. 256, Theorem 0.1 (2), from such an artificial restriction (no digonal facets). See footnote 22 below.}.

The number of vertices is finite: $V$ is finite, and three
distinct unit spheres in ${\mathbb R}^3$ have at most two points
in common. Note that if $\# V = 2$, then ${\cal B}(V)$ is
``lentil'' shaped with two small caps as facets and no
vertices at all. Put vert$F_p =_{\rm def}$ (vert$\cB(V)) \cap F_p\, (p \in V)$ and call a point of vert$\,F_p$
\emph{vertex} of $F_p$.

\textbf{Proposition 6.1:} \emph{If $\# V \ge 3$, then the following properties hold true.}
\begin{itemize}
\item[(i)] \emph{A facet
$F_p$} ($p \in V$) \emph{of ${\cal B}(V)$ has at least two vertices which are principal vertices of $\cB(V)$}.
\item[(ii)] ${\rm vert}F_p\subset {\rm relbd}F_p$.
\end{itemize}

\textbf{Proof:} (i) $F_p$ is strictly spherically convex with
a non-empty interior relative to $S_p$ (= the unit sphere with center
$p$). Hence its relative boundary
on $S_p$ is a union of at least two closed circular arcs (cf. Theorem
5.4, Proposition 5.1~(ii), and Example 5.1 above), each arc having two endpoints. Each such endpoint is a
principal vertex, by Proposition 5.1~(iii).

(ii) Since relint$F_p = \{x \in S_p : \|x-q\| < 1$ for $q \in V \setminus \{p\}\}$, clearly
relint$F_p \cap {\rm vert}{\cal B}(v) = \emptyset$, i.e., $F_p \cap {\rm vert}{\cal B}(V)
\subset {\rm relbd}F_p$. \hfill \rule{2mm}{2mm}

A facet may contain exactly two principal vertices. In this case
and if it contains no dangling vertices, it is a \emph{digonal
facet}. Its two ``edges'' are circular arcs (not necessarily
of the same radii). In fact, ${\cal B}(V)$ may have only two vertices
$p,q$ and many digonal facets.

\textbf{Example 6.1:} Let $p = p (0,0,h), \, q = q (0,0,-h), \, 0
< h < 1$, and let $C = C_{pq} = S(p) \cap S(q)$ be a circle in the
$(x,y)$-plane centered at the origin of radius $r_h$ with $0 <
r_h = \sqrt{1-h^2} < 1$. Let $V = \{p_1,p_2,\dots,p_n\}, n \ge 3$, be
the set of vertices of a regular $n$-gon inscribed to $C$. Then
$V$ is tight and $p,q$ are principal vertices of ${\cal B}(V)$
(there are no other principal vertices), and all $n$ facets
$F_{p_1}, F_{p_2}, \dots, F_{p_n}$ of ${\cal B}(V)$ are
digonal (one can say that ${\cal B}(V)$ is ``rugby-ball''
shaped). When $n \ge 3$ is odd, $h$ with $0 < h < 1$ can be ``tuned''
so that all the principal diagonals of the regular polygon
$[p_1,p_2,\dots,p_n]$ are of length 1, and then all points
$p_1,p_2,\dots,p_n$ become dangling vertices of ${\cal B}(V)$.
(Simple calculation shows that this happens for $h =
\sqrt{1-\frac{1}{4 \cos^2 \frac{180^o}{2n}}}$.) In this case
the intersection of ${\cal B}(V)$ with the $(x,y)$-plane is a Reuleaux polygon of order $n$,
whose vertices are the dangling vertices of $\cB(V)$.

\textbf{Edges:}

Next we turn our attention to the \emph{edges} of ${\cal B}(V)$.
We would like to define an edge as the intersection $F_p \cap F_q$ of two facets
$F_p$ and $F_q$, provided the intersection is
one-dimensional and connected. But it may happen that the
intersection $F_p \cap F_q$ is not connected, and that some
connectivity components of $F_p \cap F_q$ are singletons (i.e.,
$0$-dimensional); see Example 5.1 above
for such a possibility. Hence it is temptable to define an edge
of ${\cal B}(V)$ as a one-dimensional connectivity component (in
the topological sense) of the intersection $F_p \cap F_q$ of two
facets. Again this may fail, since a
connectivity component of $F_p \cap F_q$ may contain dangling
vertices in its relative interior, which is un-acceptable for an
edge, of course. This leads us to the following refined definition of an edge.

\textbf{Definition 6.2 (edge of \boldmath${\cal B}(V)$\unboldmath):} An \emph{edge} of ${\cal B}(V)$ is the closure of a connectivity
component of $(F_p \cap F_q) \setminus ({\rm vert} {\cal B}(V))$, where $\{p,q\}$ ranges over all pairs of distinct points of
$V$.

Note that this definition obliterates the possibility of an edge being a singleton.

The intersection $F_p \cap F_q$ satisfies
\begin{equation}\label{sternsternsternstern}
F_p \cap F_q = C_{pq} \cap \bigcap \limits_{v \in V \setminus
\{p,q\}} B(v)\,,
\end{equation}
where $C_{pq}$ is the circle of radius $r_{pq}$, $0 < r_{pq} =
\sqrt{1 - \frac{1}{4}\|p-q\|^2} < 1$, centered in $\frac{1}{2}
(p+q)$ lying on the plane that bisects the segments
$[p,q]$ perpendicularly. For $\# V = 2$ (i.e., $V = \{p,q\})$,
$C_{pq}$ is the only edge of ${\cal B}(V)$ (and vert${\cal B}(V) =
\emptyset$). For $\# V \ge 3$, an edge is a relatively closed (short or long) arc of
$C_{pq}$ with two endpoints, both being vertices of $\cB(V)$; call such a point \emph{vertex} of the edge.
(By
Proposition 5.1 (i) an edge cannot be the whole of $C_{pq}$, and by our definition it cannot be a
singleton; so an edge has two vertices.) From this we infer:
\begin{enumerate}
\item[(A)] $F_p \cap F_q$ is the union of at most $\#V-2$ edges.
\item[(B)] Every relative boundary point of $F_p \cap F_q$
(relative to $C_{pq}$) is a principal vertex of ${\cal B}(V)$. In
particular, an isolated point of $F_p \cap F_q$ (relative to
$C_{pq}$) is a principal vertex. (Such a point is possible only for
$\# V \ge 4$, by Proposition 5.1 (iv).)
\item[(C)] A dangling vertex
of ${\cal B}(V)$ is contained in the relative interior of a non-singleton
connectivity component of $F_p \cap F_q$ for some $p,q \in V$.
\end{enumerate}

\textbf{Proof of (C):} A dangling vertex $v$ is a point of $V$
satisfying $\|v-p\|=\|v-q\|=1$ for some $p,q \in V \setminus
\{v\}, p \not= q$, and $\|v-w\| < 1$ for $w \in V \setminus
\{v,p,q\}$. Thus $v \in F_p \cap F_q = C_{pq} \cap {\cal B}(V)$,
and a slight movement of $v$ on $C_{pq}$ ``back and forth'' leaves
it in $F_p \cap F_q$, i.e., there is a neighborhood of $v$
relative to $C_{pq}$ contained in $F_p \cap F_q$, i.e., $v \in
\mbox{ relint } (F_p \cap F_q)$. \hfill \rule{2mm}{2mm}
\begin{enumerate}
\item[(D)] A dangling vertex $v$ of ${\cal B}(V)$ is incident with
exactly two edges of ${\cal B}(V)$, both contained in the connectivity
component of $v$ in $C_{pq} \cap {\cal B}(V)$.
\item[(E)] For $\#
V \ge 3$ the boundary of a facet $F_p$ is a circuit of (at least two)
edges separated by vertices. These vertices include at least two principal vertices
and possibly dangling
vertices appearing in the relative interiors of
connectivity components of $F_p \cap F_q$, where $q$ ranges over $V
\setminus \{p\}$.
\end{enumerate}

\textbf{Definition 6.3 (\boldmath$1$\unboldmath-skeleton):} The
geometric graph whose set of vertices is vert$\cB(V)$ and whose edges are
the edges of ${\cal B}(V)$ is the \emph{$1$-skeleton} of ${\cal
B}(V)$, denoted by skel$_1 \cB(V)$.

As we already saw (cf. Examples 5.1 and 6.1 above), skel$_1 \cB(V)$
may have digons, as well as edges
of arbitrary (finite) multiplicity. See also Example 6.2 below.
For $p \in V$ we have
\begin{equation}\label{3a}
{\rm skel}_1 \cB(V) = ({\rm skel}_1 \cB(V \setminus \{p\}) \cap B(p)) \cup ({\rm relbd} F_p)\,,
\end{equation}
and the union is disjoint except for the set of principal vertices of $F_p$.

In fact, the set of (at least two) principal vertices of $F_p$ is (skel$_1 \cB(V \setminus \{p\})) \cap ({\rm relbd}F_p)$ and,
by Corollary 6.1 below, for $x \in ({\rm vert}\cB(V \setminus \{p\})) \cap F_p$ the valence of $x$ in skel$_1 \cB(V)$ is greater by $1$
than its valence in $\cB(V \setminus \{p\})$. Similarly, if $x \in {\rm vert} F_p \subset {\rm vert} \cB(V)$ is of valence $\ge 4$ in
skel$_1 \cB(V)$, then $x \in {\rm vert} \cB(V \setminus \{p\})$ and $x$ is a principal vertex of $\cB(V \setminus \{p\})$.

If $x \in {\rm vert} F_p$ and the valence of $x$ in skel$_1 \cB(V)$ is $3$, then either $x$ lies in the relative interior relint$e$ of some edge
$e$ of skel$_1 \cB(V \setminus \{p\})$, or $x$ is a dangling vertex of $\cB(V \setminus \{p\})$.

\textbf{Definition 6.4 (The face complex \boldmath${\cal S}{\cal F} ({\cal
B}(V))$\unboldmath):}

The set of faces of ${\cal B}(V)$, including facets, edges,
vertices, and the improper faces ${\cal B} (V)$ and $\emptyset$,
is the \emph{spherical face complex} of ${\cal B}(V)$ denoted by
${\cal S}{\cal F} ({\cal B}(V))$.

The complex ${\cal S}{\cal
F} ({\cal B}(V))$ is a poset (under inclusion) which may differ considerably from
the face lattice of ordinary convex $3$-polytopes. The intersection of two
facets may have more than one component (cf. Example 5.1 above), there may be digonal
facets and there may be only two vertices (cf. Example 6.1 above).
Consequently, the set of faces ${\cal S}{\cal F}({\cal
B}(V))$, partially ordered by inclusion, need not be a lattice (cf. footnote 5 above),
and its $1$-skeleton skel$_1{\cal B}(V)$ is a planar graph which
is not necessarily $3$-connected. Nevertheless, we have

\textbf{Theorem 6.1:} \emph{For a finite set $V \subset {\mathbb
R}^3\,, \, \# V \ge 3$, which is tight with} cr$(V) < 1$,
\emph{the $1$-skeleton ${\rm skel}_1 ({\cal B}(V))$ of ${\cal
S}{\cal F}({\cal B}(V))$ is $2$-connected $($in the sense of graph
theory$)$}\footnote{In \cite{Bezdek}, p. 220, Claim 9.4, it is proved that the $1$-skeleton of a so-called ``standard''
ball polytope is $3$-connected, where ``standard'' is a ball polytope every supporting sphere of which
(i.e., a sphere which bounds the polytope but does not meet its interior)
 intersects it in a connected set, i.e., a topological
disk (of some dimension, cf. Definition 6.4 there). This assumption is too restrictive since, as Example 5.1 above shows,
the $1$-skeleton of a ball polytope is
\emph{not} $3$-connected generically. So
the word ``standard'' there should not be taken literally.}.

\textbf{Proof:} We use induction on $n=_{\rm def} \# V$. For $n=3$, ${\rm skel}_1 \cB(V)$ has two
(principal) vertices (and an edge of multiplicity $3$ between them, i.e., $3$ edges), which is a $2$-connected graph
(in the void sense).

\textbf{Induction step \boldmath$n \to n+1 \, (n \ge 3)$\unboldmath:} Let $\{x,y,z\} \subset {\rm vert} \cB(V), \, x \not= y \not= z \not= x$.
It is necessary (and sufficient) to show that there is a path $\gamma_{xy} \subset {\rm skel}_1 \cB(V)$ with endpoints $x,y$ which avoids $z$.

The basic idea is easily illustrated for the case that $\{x,y,z\} \subset {\rm vert} \cB(V \setminus \{p\})$ for some $p \in V$. Then there is a path
$\tilde{\gamma}_{xy} = {\rm skel}_1 \cB(V \setminus \{p\})$ with endpoints $x,y$ which avoids $z$ (induction hypothesis). If $\tilde{\gamma}_{xy} \subset
{\rm skel}_1 \cB(V)$, define $\gamma_{xy} =_{\rm def} \tilde{\gamma}_{xy}$. Otherwise, we shall modify $\tilde{\gamma}_{xy}$ using (10)
(and the discussion thereof). Assuming that $\tilde{\gamma}_{xy}$ is parameterized by $t, \, a \le t \le b$ (``time''), with $\tilde{\gamma}_{xy} (a) = x$ and
$\tilde{\gamma}_{xy} (b) = y$, define
\[
t_{\min} = \min \{a \le t \le b: \tilde{\gamma} (t) \in {\rm relbd} F_p\}, \, x' =_{\rm def} \tilde{\gamma}_{xy} (t_{\min})
\]
and
\[
t_{\max} = {\max} \{a \le t \le b: \tilde{\gamma}_{xy} (t) \in {\rm relbd} F_p\}, \, y' =_{\rm def} \tilde{\gamma}_{xy} (t_{\max})\,.
\]
We say that $x'\,[y']$ is the \emph{first} [\emph{last}] point of $\tilde{\gamma}_{xy}$ which lies on relbd$F_p$ (equivalently: on $S_p$), or that
$\tilde{\gamma}_{xy}$ leaves [returns to] $B(p)$ at $x'\, [y']$ for the first [last] time (possibly $x' = y'$). Since $x',y' \subset ({\rm skel}_1 \cB
(V \setminus \{p\})) \cap ({\rm relbd} F_p)$,  $x'$ and $y'$ are (principal) vertices of $\cB(V)$, and they divide relbd$F_p$ ($\subset {\rm skel}_1 \cB(V)$,
cf. (10) above) into two paths $A^1_{x'y'}, A^2_{x'y'}$ both with endpoints $x',y'$ $(A^1_{x'y'} \cap A^2_{x'y'} =
\{x',y'\})$, at least one of which, say $A^1_{x'y'}$, avoids $z$ (if $z \notin {\rm relbd} F_p$, then both $A^i_{x'y'} \, (i = 1,2)$ avoid $z$).
Denote by $\tilde{\gamma}_{xx'}\, [\tilde{\gamma}_{yy'}]$ the part of $\tilde{\gamma}_{xy}$ between $x$ and $x'$, $a \le t \le t_{\min}$ $[y$ and $y',
t_{\max} \le t \le b$] and define $\gamma_{xy} =_{\rm def} \tilde{\gamma}_{xx'} \circ A^1_{x'y'} \circ \tilde{\gamma}_{y'y}$ (``$\circ$'' means
concatenation). Note that we used the facts that relbd$F_p$ is a circuit in skel$_1 \cB(V)$ and that the intersection of a path in skel$_1 \cB(V \setminus
\{p\})$ with $S_p$ is a vertex of $F_p$ (contained in relbd$F_p$).

This idea will reappear now in variants adapted to the various cases to be distinguished in the following.

Consider first the case that $\{x,y\} \subset {\rm relbd} F_p \,(\subset {\rm skel}_1 \cB(V))$ for some facet $F_p$ $(p \in V)$ of $\cB(V)$. Then
$x,y$ divide relbd$F_p$ into two paths $A^1_{xy}, A^2_{xy}$ both with endpoints $x,y\, (A^1_{xy} \cap A^2_{xy} = \{x,y\})$, at least one of which,
say $A^1_{xy}$, avoids $z$ (if $z \notin {\rm relbd} F_p$, i.e.,
$z \notin S_p$, then both $A^i_{xy}, i=1,2$, avoid $z$), and define $\gamma_{xy} =_{\rm def} A^1_{xy}$.

From now on we assume that $\{x,y\} \not\subset F_p$ for all $p \in V$; in particular, $x$ and $y$ are not both vertices of the same edge of skel$_1 \cB(V)$ nor
do they belong to the same edge of skel$_1 \cB(V \setminus \{p\})$. Fix $p \in V$. Our case analysis begins with

\textbf{Case 1:} $\{x,y\} \subset {\rm vert} \cB(V \setminus \{p\})$.

\textbf{Subcase 1.1:} $z \in {\rm vert} \cB(V \setminus \{p\})$. This case was fully discussed above.

\textbf{Subcase 1.2:} $z \notin {\rm vert} \cB(V \setminus \{p\})$. Then $z \in {\rm relbd} F_p$, and by (10) (and the discussion
thereof) there is an edge $e_{uv} \in {\rm skel}_1 \cB(V \setminus \{p\})$ with endpoints $u,v \in {\rm vert} \cB(V \setminus \{p\})$ such that
$z \in {\rm relint} e_{uv}$. Since $\{x,y\} \not= \{u,v\}$ (otherwise $x,y$ lie on the same edge $e_{uv}$ of skel$_1 \cB(V \setminus \{p\})$, contrary
to our assumption), $u \notin \{x,y\}$ or $v \notin \{x,y\}$, say $v \notin\{x,y\}$. We have now $v\not= x \not= y \not= v, \{x,y,v\} \subset
{\rm vert} \cB(V \setminus \{p\})$, and by the induction hypothesis there is a path $\tilde{\gamma}_{xy} \subset {\rm skel}_1 \cB(V \setminus \{p\})$ with
endpoints $x,y$ such that $v \notin \tilde{\gamma}_{xy}$. Clearly $z \notin \tilde{\gamma}_{xy}$. If $\tilde{\gamma}_{xy} \cap {\rm relbd} F_p = \emptyset$,
then $\tilde{\gamma}_{xy} \subset {\rm skel}_1 \cB(V)$, and we define $\gamma_{xy} =_{\rm def} \tilde{\gamma}_{xy}$.

Assume $\tilde{\gamma}_{xy} \cap {\rm relbd} F_p \not= \emptyset$, let $x'\, [y']$ be the first [last] point of $\tilde{\gamma}_{xy}$ which lies
on relbd$F_p$ (possibly $x' = y'$; surely $z \notin \{x',y'\}$ (clear)), and denote by $\tilde{\gamma}_{xx'}\, [\tilde{\gamma}_{yy'}]$ the part
of $\tilde{\gamma}_{xy}$ between $x$ and $x'\, [y$ and $y'$]. Points $x',y'$ divide relbd$F_p$ into two paths $A^1_{x'y'}$ and $A^2_{x'y'}$, both
with endpoints $x',y'$ $(A^1_{xy} \cap A^2_{x'y'} = \{x',y'\})$, at least one of which, say $A^1_{x'y'}$, avoids $z$. Define $\gamma_{xy} = \tilde{\gamma}_{xx'}
\circ A^1_{x'y'} \circ \tilde{\gamma}_{y'y}$.

\textbf{Case 2:} $\{x,y\} \not\subset {\rm vert} \cB(V \setminus \{p\})$, say $x \notin \cB(V \setminus \{p\})$ (note that the roles
of $x,y$ are interchangeable). By (10) (and the discussion thereof) $x$ is a (principal) $3$-valent vertex of skel$_1 \cB(V)$ lying on relbd$F_p \subset
{\rm skel}_1 \cB(V)$ and $x \in {\rm relint} e_{uv}$, where $e_{uv} \in {\rm skel}_1 (\cB \setminus \{p\})$ is an edge not lying on relbd$F_p$, with
vertices $u,v \in {\rm vert} \cB(V \setminus \{p\})$. Denote by $e_{xu}\, [e_{xv}]$ the part (circular arc) of $e_{uv}$ between $x$ and $u$ $[x$
and $v$]. Since $e_{uv} \not\subset F_p$ is a circular arc, it may have at most one point of intersection with relbd$F_p$ in addition to $x$. Call this
point $w$ -- if it exists -- and assume, w.l.o.g., that $w \in e_{xv}$, i.e., $w \notin e_{xu}$ (if $w$ exists; possibly $w = v$). It
follows from this convention that $u \notin {\rm relbd} F_p$, hence $u \notin S_p$, i.e., either $u \in {\rm int} B(p)$ or $u \in \bR^3 \setminus B(p)$
and $e_{xu}$ is an edge of $\cB(V)$. We also make the following
convention and notation: since $x,w \in {\rm relbd} F_p$ (if $w$ exists), the points $x,w$ divide relbd$F_p$ into two paths $A^1_{xw}, A^2_{xw}$ both with
endpoints $x,w\, (A^1_{xw} \cap A^2_{xw} = \{x,w\})$ at least one of which, say $A^1_{xw}$, does not contain $z$ in its relative interior (possibly $w = z$;
if $z \notin {\rm relbd} F_p$, then both $A^i_{xw} (i = 1,2)$ do not contain $z$). Note that $A^1_{xy} \subset {\rm skel}_1 \cB(V)$.

Since $y \notin {\rm relbd} F_p$ (by our basic assumption $x,y$ do not belong to the same facet of $\cB(V)$), $y \in {\rm int} \cB(p)$ and $y \in {\rm vert}
\cB(V \setminus \{p\})$. If $u=y$, then $e_{xu} = e_{xy}$ is an edge of skel$_1 \cB(V)$, a contradiction to our assumption; hence $u \not= y$. If $v=y$,
then $e_{xv} = e_{xy}$ is an edge of skel$_1 \cB(V \setminus \{p\})$, contradicting our assumption that $x,y$ do not lie on the same edge of skel$_1 \cB(V
\setminus \{p\})$. To sum up,
\begin{equation}\label{aaa}
y \notin \{u,v\}\,.
\end{equation}
We distinguish between two subcases.

\textbf{Subcase 2.1:} $z \in {\rm vert} \cB(V \setminus \{p\})$.

Again we distinguish two subcases.

\textbf{Sub-subcase 2.1.1:} $z \not= u$.

We have $u \not= z \not= y \not= u$ (see (\ref{aaa})) and $\{u,y,z\} \subset {\rm skel}_1 \cB(V \setminus \{p\})$. By the induction hypothesis there is a path
$\tilde{\gamma}_{uy} \subset {\rm skel}_1 \cB(V \setminus \{p\})$ with endpoints $u,y$ such that $z \notin \tilde{\gamma}_{uy}$. If $\tilde{\gamma}_{uy}
\cap {\rm relbd} F_p = \emptyset$ (i.e., $\tilde{\gamma}_{uy} \cap S_p = \emptyset$), then $u \in {\rm int} B(p)$ (recall that $y \in {\rm int} B(p)$), hence
$e_{xu}$ is an edge of $\cB(V)$ and we define $\gamma_{xy} =_{\rm def} e_{xu} \circ \tilde{\gamma}_{uy}$. Assume now that $\tilde{\gamma}_{uy} \cap S_p \not=
\emptyset$; call $y'$ the last point of $\tilde{\gamma}_{uy}$ which lies on relbd$F_p$ (possibly $y'=x$). Since $z\not= x$ and $z \not= y'$, the points $x,y'$
divide the circuit relbd$F_p$ into two paths $A^1_{xy'}, A^2_{xy'}\, (A^1_{xy'} \cap A^2_{xy'} = \{x,y'\})$ at least one of which, say $A^1_{x,y'}$,
avoids $z$ (if $z \notin {\rm relbd} F_p$, i.e., $z \notin S_p$), then both paths $A^i_{xy'}\, (i = 1,2)$ avoid $z$). Denote by $\tilde{\gamma}_{yy'}$ the
part of $\tilde{\gamma}_{uy}$ between $y'$ and $y$ and define $\gamma_{xy} =_{\rm def} A^1_{xy'} \circ \tilde{\gamma}_{y'y}$.

\textbf{Sub-subcase 2.1.2:} $z = u$.

Then $z \in ({\rm relint} B(p)) \cap ({\rm vert} \cB(V \setminus \{p\}))$. If $v = y$, then $v \in B (p)$ and $w$ is defined;
denote by $e_{wv}$ the part of $e_{uv}$ between $w$ and $v$ (in case $w = v$, $e_{wv}$ is just a point); note that $e_{wv}$ is an edge of skel$_1 \cB(V)$ (or
a vertex in case $w = v$). Define $\gamma_{xy} = A^1_{xw} \circ e_{wv}$. Assume now that $v \not= y$. Then $v \not= y \not= z \not= v$ (recall that
$z = u \not= v$) and $\{v,y,z\} \subset {\rm vert} \cB(V \setminus \{p\})$. By the induction hypothesis there is a path $\tilde{\gamma}_{vy} \subset {\rm skel}_1
\cB(V \setminus \{p\})$ with endpoints $v,y$ which avoids $z$. If $\tilde{\gamma}_{vy} \cap ({\rm relbd} F_p) = \emptyset$, then $v \in B(p)$, $w$ is defined and we define
$\gamma_{xy} =_{\rm def} A^1_{xw} \circ e_{wv} \circ \tilde{\gamma}_{vy}$.

Assume now that $\tilde{\gamma}_{vy} \cap ({\rm relbd} F_p) \not= \emptyset$; call $y'$ the last point of $\tilde{\gamma}_{vy}$ which lies on
relbd$F_p$. Points $x,y'$ divide relbd$F_p$ into two paths $A^1_{xy'},A^2_{xy'}$, both with endpoints $x,y'$~
$(A^1_{xy'} \cap A^2_{xy'} = \{x,y'\})$ both of which avoid $z$ (since $z = u \notin {\rm relbd} F_p$).
Denote by $\tilde{\gamma}_{y'y}$ the part of $\tilde{\gamma}_{vy}$ between $y'$ and $y$ and define $\gamma_{xy} =_{\rm def} A^1_{xy'} \circ \tilde{\gamma}_{y'y}$.

\textbf{Subcase 2.2:} $z \notin {\rm vert} \cB(V \setminus \{p\})$.

Then $z \in {\rm relbd} F_p, z$ is a (principal) vertex of valence $3$ in skel$_1 \cB (V)$ and there is an edge $e_{st} \in {\rm skel}_1 \cB(V \setminus
\{p\})$ with endpoints $s, t \in {\rm vert} \cB(V \setminus \{p\})$ such that $z \in {\rm relint} e_{st}\, (z \notin \{s,t\})$. Now we distinguish again two cases.

\textbf{Sub-subcase 2.2.1:} $u \notin \{s,t\}$.

Since $y \not= s$ or $y \not= t$ (or both), we assume, w.l.o.g., that $y \not= s$. Then $s \not= u \not= y \not= s$ (see (\ref{aaa})), $\{u,s,y\} \subset
{\rm vert} \cB(V \setminus \{p\})$, and by the induction hypothesis there is a path $\tilde{\gamma}_{uy} \subset {\rm skel}_1 \cB(V \setminus \{p\})$ with endpoints
$u$ and $y$ which avoids $s$. Since $e_{st} \notin \tilde{\gamma}_{uy}$ and $\tilde{\gamma}_{uy}$ does not use any edge of relbd$F_p$, $z \notin \tilde{\gamma}_{uy}$.
If $\tilde{\gamma}_{uy} \cap {\rm relbd} F_p = \emptyset$, define $\gamma_{xy} =_{\rm def} e_{xu} \circ \tilde{\gamma}_{uy} \subset {\rm skel}_1 \cB(V)$.

Assume now that $\tilde{\gamma}_{uy} \cap {\rm relbd} F_p \not= \emptyset$ and call $y'$ the last point of $\tilde{\gamma}_{uy}$ which lies on relbd$F_p$
(equivalently: on $S_p$). Points $x,y'$ divide relbd$F_p$ into two paths $A^1_{xy'}, A^2_{xy'}$, both with endpoints $x,y'$ $(A^1_{xy'} \cap A^2_{xy'}
= \{x,y'\})$, (exactly) one of which, say $A^1_{xy'}$, avoids $z$. Denote by $\tilde{\gamma}_{y'y}$ the part of $\tilde{\gamma}_{xy'}$ between $y'$ and $y$
and define $\gamma_{xy} =_{\rm def} A^1_{xy'} \circ \tilde{\gamma}_{y'y}$.

\textbf{Sub-subcase 2.2.2:} $u \in \{s,t\}$, say $u = s$.

Since $u \not= y$ (cf. (\ref{aaa})), we have $u = s \not= y$. Again we distinguish two subcases.

\textbf{Sub-sub-subcase 2.2.2.1:} $y \not= t$.

Then $(s=)\, u \not= y \not= t \not= u$ and $\{u,y,t\} \subset {\rm skel}_1 \cB(V \setminus \{p\})$. By the induction hypothesis there is a path $\tilde{\gamma}_{uy}
\subset {\rm skel}_1 \cB(V \setminus \{p\})$ with endpoints $u$ and $y$ which avoids $t$. The rest of the argument is exactly as in Sub-subcase 2.2.1 above
from the wording ``Since $\tilde{\gamma}_{wy}$ does not use $e_{st}$...'' to the end ``Define $\gamma_{xy} =_{\rm def} A^1_{xy'} \circ \tilde{\gamma}_{y'y}$''.

\textbf{Sub-sub-subcase 2.2.2.2:} $y = t$.

Then $v \not= t$ (recall that $v \not= y$; cf. (\ref{aaa})), and we have $(s=)\, u \not= v \not= y\, (=t) \not= u\, (=s)$ and $\{u,v,y\} \subset {\rm vert} \cB(V \setminus
\{p\})$. By the induction hypothesis there is a path $\tg_{vy} \subset {\rm skel}_1 \cB(V \setminus \{p\})$ with endpoints $v$ and $y$ $(=t)$ which avoids
$u\, (=s)$. Since $e_{st} \notin \tg_{vy}$ and $\tg_{vy}$ does not use any edge of relbd$F_p, z \notin \tg_{vy}$. If $\tg_{vy} \cap {\rm relbd} F_p = \emptyset$,
then $v \in {\rm int} B(p)$ (recall that $y \in {\rm int} B(p))$, $w$ exists (see above) and, by our conventions definition of $w$ above, $z \notin A^1_{xw}$. Call
$e_{wv}$ the part of $e_{uv}$ between $w$ and $v$ ($e_{wv} \in {\rm skel}_1 \cB(V)$) and define $\gamma_{xy} =_{\rm def} A^1_{xw} \circ e_{wv} \circ
\tg_{vy}$.

Assume now that $\tg_{vy} \cap {\rm relbd} F_p \not= \emptyset$. Call $y'$ the last point of $\tg_{vy}$ which lies on relbd$F_p$, and denote by $\tg_{y'y}$
the part of $\tg_{vy}$ between $y'$ and $y$. Points $x,y'$ divide relbd$F_p$ into two paths $A^1_{xy'}, A^2_{xy'}$, both with endpoints $x$ and $y'$
$(A^1_{xy'} \cap A^2_{xy'} = \{x,y'\})$, (exactly) one of which, $A^1_{xy'}$ say, avoids $z$. Define $\gamma_{xy} =_{\rm def} A^1_{xy'} \circ \tg_{y'y}$.
\hfill \rule{2mm}{2mm}

\textbf{Remarks 6.1:}
\begin{enumerate}
\item An alternative approach to Theorem 6.1, via the dual graph of $\cB(V)$, can be given along the following lines. First prove as a lemma that the dual
graph of $\cB(V)$ is connected. Let $F,G$ be two facets of $\cB(V)$, let $x \in {\rm relint} F, \,y \in {\rm relint} F$, and let $H$ be a plane passing
through $x$ and $y$ such that $H \cap {\rm vert} \cB(V) = \emptyset$. Then $H$ intersects int$\cB(V)$ and $\gamma =_{\rm def} H \cap ({\rm bd} \cB(V))$ is
a (closed) convex curve avoiding vert$\cB(V)$, and $x,y$ divide $\gamma$ into two paths $\gamma^1_{xy}, \gamma^2_{xy}$ both with endpoints $x,y \, (\gamma^1_{xy}
\cap \gamma^2_{xy} = \{x,y\})$. The facets of $\cB(V)$ met by $\gamma^1_{xy}$ form a succession of facets from $F$ to $G$, each two successive facets
share an edge. This is a path between $F$ and $G$ in the dual graph of $\cB(V)$. Let $\{x,y,z\} \subset {\rm vert} \cB(V) \, (x \not= y \not= z \not= x)$,
let $F,G$ be facets of $\cB(V)$ containing $x$ and $y$, respectively, and let $(F_0 = F,F_1,\dots,F_n = G)$ be a path in the dual graph of $\cB(V)$ from
$F$ to $G$ (assured by the lemma described above). The union of the relative boundaries of $F_i, \bigcup \limits^n_{i=0} {\rm relbd} F_i$, is clearly a connected subgraph
of skel$_1\cB(V)$, hence there is a path $\tg_{xy}$ between $x$ and $y$ in skel$_1 \cB(V)$.

If $z \notin \tg_{xy}$, then $\tg_{xy}$ is a path in skel$_1 \cB(v)$ which avoids $z$. If $z \in \tg_{xy}$, then $z \in {\rm relbd} F_i$ for some $0 \le i \le n$,
and (relbd$F_i$)$\cap \tg_{xy}$ is a path $A^1_{ab}$ on relbd$F_i$ with endpoints $a,b$, say $(z \in A^1_{ab}$). Replace $A^1_{ab}$ by its complementary path
$A^2_{ab} =_{\rm def}{\rm cl}$(relbd$F_i \setminus A^1_{ab}$). Then $z \notin A^2_{ab}$, and the obtained path $\gamma_{xy}$ connects $x$ and $y$ in
skel$_1 \cB(V)$ and avoids $z$.
\item The dual graph of $\cB(V)$ is $2$-connected as well (as skel$_1 \cB(V)$). For a proof one can use an approach similar to the one given in the
previous remark. A shorter proof, using already established results, runs as follows. The Bezdek-Nasz\'{o}di theorem, described in footnote 22 below,
assures that the dual graph of $\cB(V)$ is isomorphic to the $1$-skeleton of the ball polytope $\cB$(vert$\cB(V)$) (in fact, $\cB(V)$
and $\cB$(vert$\cB(V)$) have dual face
structures), and by Theorem 6.1 above skel$_1 \cB$(vert$\cB(V)$) is $2$-connected.
\end{enumerate}

\textbf{Question 6.1:} \emph{Is every $2$-connected planar graph
realizable as the $1$-skeleton of some Ball polytope ${\cal B}(V)$}\footnote{See \cite{Bezdek},
Problem 9.5 (p. 221), for a similar question; note that the condition ``with no loops'' there is superflous
and does not appear in our formulation.}?

Since skel$_1 \cB(V)$ is a finite connected (even $2$-connected, by Theorem 6.1) planar graph embedded on the boundary bd$\cB(V)$ of the convex
body $\cB(V)$, we have:

\textbf{Proposition 6.2:} \emph{If $V \subset {\mathbb R}^3$ is
finite with $\# V \ge 3$ and} cr$(V) < 1$, \emph{then the face
numbers $v$ (vertices), $e$ (edges) and $f$ $($facets$)$ of ${\cal
S}{\cal F} ({\cal B}(V))$ satisfy the Eulerian relation}
$v-e+f=2$\footnote{This appears as ``the Euler-Poincare theorem'' in \cite{Bez-Nasz}, p. 256, with an alleged proof
given in p. 260. Instead of saying that this relation follows merely from the planarity and connectivity
of the graph skel$_1 \cB(V)$, an argument is given based on the sophisticated notion of $CW$-decomposition from
homology theory, not necessary in this simple case. Another drawback of this proof is that the whole face structure of
$\cB(V)$ is compressed there into two sentences.}.

\textbf{Example 6.2:} The following Ball polytope $\cB(V)$, generalizing Example 5.1
above, exhibits the various notions of faces (vertices, edges, facets) of $\cS\cF(\cB(V))$ defined above in a non-trivial
way. We first briefly introduce the basic notion of $2$-tuples in equilateral position, which is fully defined
in Definition 8.2 below, in connection with the notion of dual arcs.

Let $a,b,x,y \in \bR^3$ be four points satisfying $\|x-a\|=\|x-b\|=\|y-a\|=\|y-b\|=1$.
The couple of $2$-tuples $(a,b;x,y)$ are in \emph{equilateral position}. Denote by $C_{xy}$, resp. $C_{ab}$, the
circle centered at $\frac{1}{2}(x+y)$, resp. $\frac{1}{2}(a+b)$, that passes through $a,b$, resp. $x,y$.

In Lemma 8.1 below it is proved in detail that if $c \in C_{xy}$, resp. $z \in C_{ab}$, and if both $c$ and $z$ lie on
the relative interiors of the short [resp. both lie in the long] arcs of $C_{xy}$, resp. $C_{ab}$, determined by
$a,b$, resp. $x,y$, then $\|c-z\|>1$. Otherwise $\|c-z\|\le 1$, and this is satisfied as an equality iff $c \in \{a,b\}$
or $z \in \{x,y\}$. We use this to generalize Example 5.1 above. With the notation $p = p(0,0,h), q = q(0,0,-h), 0 < h < 1,
C_{pq} = S(p) \cap S(q), p_0,p_1,\dots,p_{n-1} \subset C_{pq} \, (n \ge 2)$ and ${\rm orarc} (p_{i-1},p_i)$, as introduced
there, define $n$ relatively open arcs $A_i, \, 1 \le i \le n$ (replacing the $n$ points $v_i, \, 1 \le i \le n$, of
Example 5.1) in the following way: For $1 \le i \le n$, if ${\rm orarc} (p_{i-1},p_i)$ is a short [resp. long] arc
of $C_{pq}$, then $A_i$ is the relatively open short [resp. long] arc $A^i_{pq}$ [resp. $C_{p_{i-1}p_i} \setminus A^i_{pq}$]
of the circle $C_{p_{i-1}p_i}=_{\rm def} S(p_{i-1}) \cap S(p_i)$ whose endpoints are $p,q$. Note that the $v_i$'s of Example 5.1
satisfy $v_i \in A_i  \, (1 \le i \le n)$ and that since ${\rm orarc}(p_{i-1},p_i)$ is a long arc for at most one
$i, \, 1 \le i \le n$, at most one $A_i \, (1 \le i \le n)$ is a long arc. It follows from Lemma 8.1 (as explained above)
that for $v \in A_i$ the unit ball $B(v,1)$ misses the relative interior of ${\rm orarc} (p_{i-1},p_i)$ and includes
(in its interior) the relative interior of the complementary arc ${\rm orarc} (p_i,p_{i-1})$. Let $W \subset
\bigcup \limits^n_{i=1} A_i$ be any finite set $(0 \le \# W < \infty)$ and put $V=_{\rm def} \{p,q\} \cup W$. Then ${\rm cr}
(V) < 1,\, V$ is tight and, similar to (8) in \S~5 above,
\begin{equation}\label{star1}
\begin{array}{l}
F_p \cap F_q = C_{pq} \cap \cB(V)=\\
~\\
\bigcup \left\{{\rm orarc} (p_{i-1},p_i): W \cap A_i = \emptyset, 1 \le i \le n\}
\cup \{p_i: W \cap A_i \not= \emptyset\,,\right.\\
~\\
\mbox{ and }\left.W \cap A_{i+1} \not= \emptyset\,, \, 1 \le i \le n \right\}\,.
\end{array}
\end{equation}
The set of principal vertices prin-vert$\cB(V)$ of $\cB(V)$ is contained in the point set $\{p_1,p_2,\dots,p_n=p_0\}$, and for $1 \le i \le n$
the point $p_i$ is a principal vertex of $\cB(V)$
$\Leftrightarrow W \cap A_i \not= \emptyset \, \mbox{ or } \, W \cap A_{i+1} \not= \emptyset$ (or both).

Before describing the dangling vertices of $\cB(V)$, we describe its $1$-skeleton and its facets. The closure of a
connectivity component of skel$_1 (\cB(V)) \setminus \mbox{prin-vert}\cB(V)$ is a circular arc with two endpoints in
$\mbox{prin-vert}\cB(V)$, and it may contain dangling vertices in its relative interior as well. We call such a
connectivity component a \emph{hedge} (to differentiate it from ``edge'').

A hedge is the union of consecutive edges separated by dangling vertices. Hence a hedge is an edge iff it does not contain a dangling
vertex (in its relative interior), and clearly the body of skel$_1 \cB(V)$ is the union of its hedges. Every
connectivity component of $F_p \cap F_q$ is either a hedge or a principal vertex, and these connectivity
components are cyclically ordered on $C_{pq}$.

Two consecutive points $p_{i-1},p_i$ of $\{p_1,\dots,p_n=p_0\}$ belong to two consecutive components of $F_p \cap F_q$
iff $A_i \cap W \not= \emptyset$ (in which case relint(${\rm orarc} (p_{i-1},p_i)$) is truncated), and if
$A_i \cap W \not= \emptyset$, then $p_{i-1}, p_i \in \mbox{prin-vert}\cB(V)$, and $p_{i-1},p_i$ are connected in skel$_1 \cB(V)$
by $\# (A_i \cap W)+1$ hedges all of which are disjoint from $C_{pq}$ except for their endpoints $p_{i-1},p_i$. These
$\#(A_i \cap W)+1$ hedges are the boundaries of $\#(A_i \cap W)$ facets $F_v$ for $v \in A_i \cap W$. Every such facet
is ``di-hedged'', since its boundary consists of two such hedges (with common endpoints $p_{i-1},p_i$); and besides
$F_p$ and $F_q$ these are all the facets of $\cB(V)$. Thus the hedges of $\cB(V)$ are either included in $C_{pq}$ or
have endpoints $p_{i-1},p_i$ for some $1 \le i \le n$ such that $A_i \cap W \not= \emptyset$.

A point $v \in W$ is a dangling vertex of $\cB(V)$ iff either (i) $v \in C_{pq} \setminus \mbox{prin-vert}\cB(V)$ (in
which case $[v,p]$ and $[v,q]$ are edges of $D(V))$, or (ii) $v$ lies in the relative interior of a hedge
with endpoints $p_{i-1},p_i$ such that $A_i \cap W \not= \emptyset$ for some $1 \le i \le n$ (i.e., ${\rm orarc} (p_{i-1},p_i)$
is truncated). Dangling vertices satisfying (i) are easy to construct, but the realization of dangling vertices satisfying
(ii) turns out to be a delicate question, which we will not address. The following criterion says when $p$ and/or $q$ are
dangling vertices.

\textbf{Criterion:} $p$, resp. $q$, is a dangling vertex of $\cB(V)$ iff $\|p-q\|=1, \, \#(W \cap C_{pq}) = 1$,
and the (unique) arc $A_i \, (1 \le i \le n)$ that contains the point $W \cap C_{pq}$ satisfies: $W \cap A_i$ is
contained in the closed half space bounded by aff$C_{pq}$ which contains $p$, resp. $q$. (Thus, if both
$p,q$ are dangling vertices, then $\#(W \cap A_i) = \#(W \cap C_{pq}) = 1$.)

This criterion finishes our Example 6.2.

The following is fundamental for all later development.

\textbf{Proposition 6.3:} \emph{Let $V \subset {\mathbb R}^3$ be
finite with $\# V \ge 3$ and tight with ${\rm cr}(V)<1$.
Let $F_p$ be the facet that corresponds to $p \in V$. Then}
\begin{itemize}
\item[(i)] cone$_p F_p$ \emph{is apexed at} $p$.
\item[(ii)] \emph{For an edge
$e$ of $F_p$ with vertices $r$ and $s$}, \emph{the set} cone$_p(r,s)$ \emph{is a facet of}
cone$_p$ (vert$F_p$), and relint$\,e \cap {\rm cone}_p({\rm vert} F_p) = \emptyset$.
\item[(iii)] \emph{The rays} cone$_p r, r \in {\rm
vert} F_p$, \emph{are the extreme rays of} cone$_p$(vert$F_p$), \emph{and their
cyclic order on} bd(cone$_p$(vert$F_p$)) \emph{is compatible with the cyclic
order of the vertices of $F_p$ on} relbd$F_p$.
\end{itemize}

\textbf{Proof:} (i) By Corollary 5.2, $F_p$ is contained in a small cap of $S_p$.\\
(ii) There is a facet $F_q \,(q \in V \setminus \{p\})$ such that $e
\subset F_p \cap F_q$ and $F_p$ is contained in the small cap of
$S_p = \{x \in {\mathbb R}^3 : \| x-p\| = 1\}$ bounded by the circle
$C_{pq} ={\rm def}\{x \in {\mathbb R}^3$: $\|x-p\|=\|x-q\|=1\}$
(of radius $0 < r_{pq} =_{\rm def} \sqrt{1-\frac{1}{4} \|p-q\|^2} < 1$). Denote
this small cap of $S_p$ by cap$(p)$. $C_{pq} = {\rm relbd}\,{\rm cap}(p)$,
the cap cap$(p) \supset F_p$ lies beyond the plane $H_{pq} =_{\rm
def} {\rm aff} C_{pq}$ relative to $p$, and the plane aff$(p,r,s)
\supset {\rm cone}_p (r,s)$ intersects $S_p$ in a great circle,
the short arc between $r$ and $s$ of which is contained in cap$(p)$.
Denote by arc$(r,s)$ this short arc on cap$(p)$. Note that arc$(r,s)$
and $e$ share $r$ and $s$ as common vertices.

Since ${\rm arc}(r,s)$ is a geodetic arc on $S_p$ both of whose endpoints belong to $F_p$,
${\rm arc}(r,s) \subset F_p$ (by the spherical convexity of $F_p$), and since $e \subset F_p$,
sph-conv$(e \cup {\rm arc} (r,s)) \subset F_p$ as well (again by the spherical convexity of $F_p$). Since $F_p$
is strictly spherically convex, ${\rm arc}(r,s) \setminus \{r,s\} \subset {\rm relint} F_p$, and
sph-conv$(e \cup {\rm arc}(r,s)) \setminus e \subset {\rm relint} F_p$. Hence
\begin{equation}\label{triangle1}
\mbox{sph-conv} (e \cup {\rm arc}(r,s)) \cap {\rm vert} F_p = \{r,s\}
\end{equation}
(cf. Proposition 6.1~(ii); note that ${\rm relint}\,e \cap {\rm vert} F_p = \emptyset$, since $e$ is an
edge of $F_p$).
The plane $H=_{{\rm def}} {\rm aff} (p,r,s)$ divides $\bR^3$ into two closed half-spaces, one of which contains $e$;
denote this half-space by $H^+$. By $F_p \subset {\rm cap} (p)$ we have ${\rm vert} F_p \subset {\rm cap} (p)$, hence
\begin{equation}\label{triangle2}
H^+ \cap {\rm vert} F_p \subset H^+ \cap {\rm cap} (p)\,.
\end{equation}
\textbf{Claim:}
\begin{equation}\label{triangle3}
H^+ \cap {\rm cap} (p) = \mbox{sph-conv} (e \cup {\rm arc} (r,s))\,.
\end{equation}
\textbf{Proof:} $H^+ \cap {\rm cap} (p) = (H^+ \cap S_p) \cap {\rm cap} (p)$, which is the intersection of two spherically
convex sets, hence $H^+ \cap {\rm cap} (p)$ is spherically convex and ``$\supset$'' in (\ref{triangle3}) follows. Let
$m \in H^+ \cap {\rm cap} (p)$ and denote by $n$ the center of $\capp (p)$ on $S_p$. The great circle of $S_p$ passing
through $n$ and $m$ cuts $H^+ \cap \capp (p)$ in a short arc $\arc (x,y)$, one of whose endpoints $x \in \arc(r,s)$ and the other
endpoint $y \in e$. Hence $m \in \arc (x,y) \subset \mbox{sph-conv} (e \cup \arc (r,s))$.
This proves ``$\subset$'' in (\ref{triangle3}), proving the claim.

It follows from (\ref{triangle1}), (\ref{triangle2}), and (\ref{triangle3}) that ${\rm vert} F_p \cap H^+ = \{r,s\}$, hence
$H$ supports ${\rm vert} F_p$ and $\cone_p (r,s)$ is a facet of $\cone_p ({\rm vert} F_p)$

(iii) It follows from (ii) that for each edge $e$ of $F_p$ with
vertices $r,s$, say, the rays cone$_p r$ and cone$_p s$ are
extreme rays of cone$_p$(vert$F_p$), and since relbd$F_p$ is a
circuit of edges separated by vertices, the cyclic order of the
rays cone$_p r, r \in {\rm vert}F_p$ on bd~cone$_p$(vert$F_p$)
is compatible with the cyclic order of the vertices of $F_p$ on relbd$F_p$.
\hfill \rule{2mm}{2mm}

The following proposition gives a simple criterion for an edge of $\cB(V)$ being short or
long (as an arc of its loading circle); it is needed, e.g., in the proof of Theorem 8.1 below.

\textbf{Proposition 6.4:} \emph{With the usual assumptions on $V$ (see} Propositions
6.2 \emph{and} 6.3 \emph{above), let $e$ be an edge of ${\cal B}(V)$ with endpoints
$a,b$, and assume that $e \subset C_{xy}$, where $x,y \in V$. Assume
that $\#({\rm vert}F_x) \ge 3$ and put $H=_{\rm def} {\rm
aff}(x,a,b)$. The following four conditions are equivalent}:
\begin{itemize}
\item[(i)] \emph{$e$ is a short arc of $C_{xy}$}.
\item[(ii)] \emph{The edge $e$ and the point $y$ lie in different sides of the plane $H$}.
\item[(iii)] \emph{The edge $e$ and the set ${\rm vert}F_x \setminus \{a,b\}$ lie in
different sides of $H$}.
\item[(iv)] \emph{The point $y$ and the set ${\rm vert}F_x \setminus \{a,b\}$ lie in
the same side of $H$}.
\end{itemize}

\textbf{Proof:} The facet $F_x$ is contained in the small cap
cap$(x) =_{\rm def} S(x,1) \cap B(y,1)$ on the unit sphere
$S(x,1)$, relbd cap$(x) = C_{xy}$, and both $y$ and cap$(x)$ lie
beyond aff$C_{xy}$ relatively to $x$. Put $y^* =_{\rm def} {\rm
cone}_x (y) \cap {\rm aff}C_{xy}$ and $y^{**} =_{\rm def} {\rm
cone}_x(y) \cap {\rm cap}(x)$. $y^*$ is the center of the disc$D_{xy}
=_{\rm def} {\rm conv}C_{xy}$, and $y,y^*,y^{**}$ lie in the same open
half space bounded by $H$. The center $y^*$ of $D_{xy}$ lies in
the plane aff$C_{xy}$ beyond the line aff$(a,b)$ relatively to the
short arc of $C_{xy}$ whose endpoints are $a$ and $b$. And since
aff$(a,b)=H \cap {\rm aff} C_{xy}$ (in fact, $[a,b] = H \cap
D_{xy})$, $y^*$ lies (in the $3$-space) beyond $H$ relatively to
this short arc. This shows that (i) $\Leftrightarrow$ (ii). Since
the facet $F_x \subset {\rm cap}(x)$ is strictly spherically convex
on $S(x,1)$, the geodetic line $\gamma = H \cap {\rm cap}(x)$ on
cap$(x)$ is contained in relint$F_x$ except for its
endpoints $a$ and $b$, and $\gamma$ separates relbd$F_x \subset
{\rm cap}(x)$ into two parts; one is $e$, and the other part,
containing vert$F_x \setminus \{a,b\}$, lies on cap$(x)$ in the
same side of $\gamma$ as $y^{**}$. Thus vert$F_x \setminus
\{a,b\}$ and $y^{**}$ lie in the same half space bounded by $H$,
which contains $y$ and $y^*$ as well. This shows that (ii)
$\Leftrightarrow$ (iii) $\Leftrightarrow$ (iv). \hfill
\rule{2mm}{2mm}

Lemmata 6.1 and 6.4 below establish a connection between the local
face structure of ${\cal B}(V)$ near a vertex $x$ and the points
of $V$ at distance $1$ from $x$.

\textbf{Lemma 6.1:} \emph{Suppose $c$ is a vertex of ${\cal
B}(V)$, where $V \subset {\mathbb R}^3$, $3 \le \# V < \infty$, $V$ is
tight and} cr$(V)<1$. \emph{Let} $\{v \in V :
\|v-c\| = 1\} = \{v_1, v_2, \dots, v_m\}$ ($m \ge 3$ \emph{for $c$
principal, and $m=2$ for $c$ dangling, $v_i \not= v_j$ for $
1 \le i < j \le m$$)$. Then the polyhedral cone} cone$_c (v_1, \dots, v_m)$
\emph{is apexed $($at $c$$)$, and the rays} cone$_c(v_i)\,, \,
i=1,2,\dots,m$, \emph{are the extreme rays of it}\footnote{Note the similarity between
Lemma 6.1 and Proposition 6.3 (iii) above. The full understanding of this similarity becomes apparent in the
rigorous proof of the Bezdek-Nasz\'{o}di duality theorem; see footnote 22 below.}.

\textbf{Proof:} Put $V(c) =_{\rm def} \{v_1,\dots,v_n\}$. First we
show that the cone cone$(V(c))$ is apexed (at $c$) (cf. Definition 4.2).
Put $z=_{\rm def} cc(V)$, the circumcenter of $V$,
and let $r$ be its circumradius ($r < 1$, by assumption). For $v \in V
(c)$ we have $\|v-c\| = 1$ and $\|v-z\| \le r < 1$. Hence $v$ lies
in the open halfspace bounded by the plane $H$ that bisects the
segment $[c,z]$ in its midpoint perpendicularly which does not
contain $c$; see Figure 5. Hence $H$ strictly separates $c$ and $V(c)$.
\vspace{1.5cm}

\begin{picture}(280,130)
\put(50,0){\includegraphics[scale=0.5]{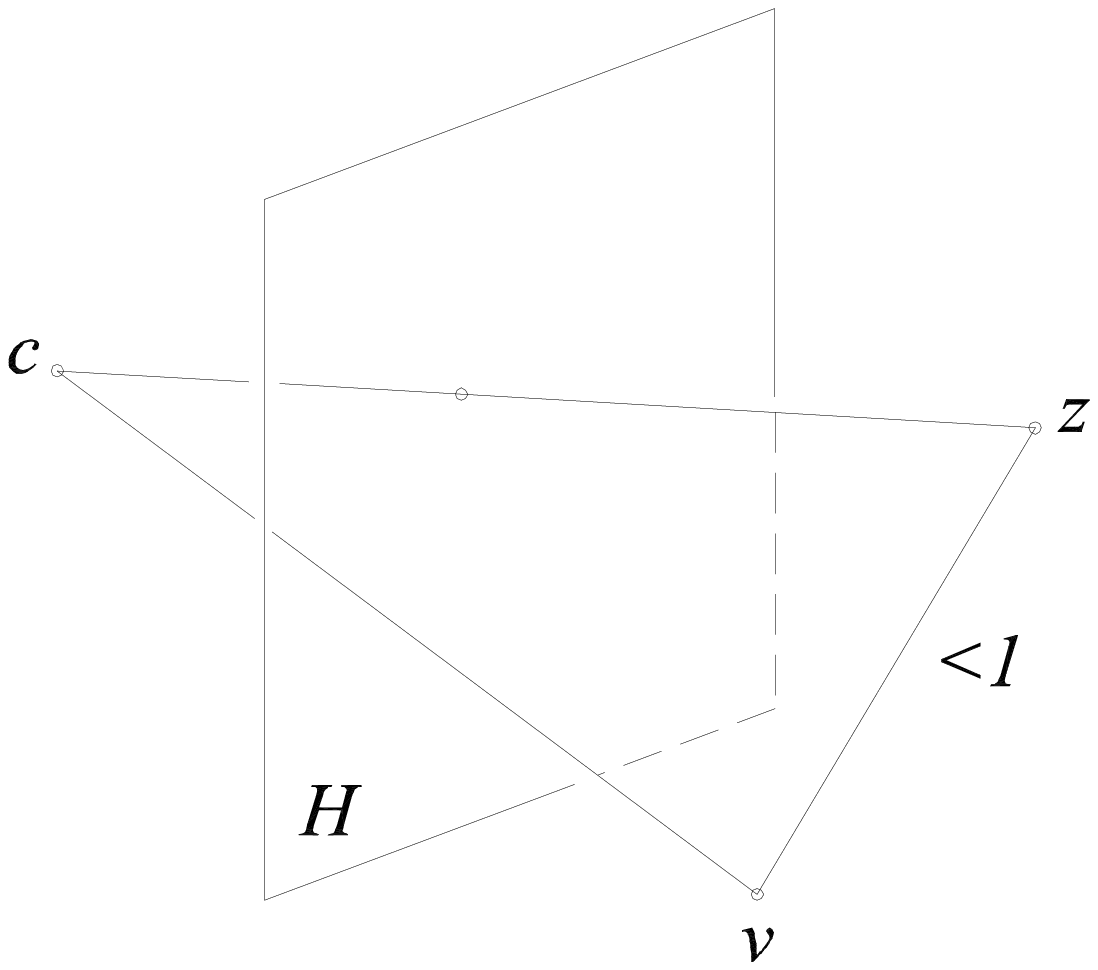}} 
\end{picture}

\begin{center}
-- Figure 5 --
\end{center}

Thus cone$_c (V (c))$ is apexed (at $c$). Next we show that
for $v \in V (c)$ the ray cone$_c (v)$ is extreme in cone$_c
(V(c))$. Assume, by r.a.a., that $\exists v_0 \in V(c)$ such that
cone$_c(v_0)$ is not an extreme ray of cone$_c (V(c))$. (Clearly,
every extreme ray of cone$_c (V(c))$ is of the form cone$_c (v)$
for some $v \in V (c)$.) Then cone$_c (V(c)) = {\rm cone}_c (V(c)
\setminus \{v_0\})$ and
\begin{equation}\label{1a}
v_0 \in S(c) \cap {\rm cone}_c (V(c) \setminus \{v_0\})\,,
\end{equation}
($S(c)$ is the unit sphere centered at $c$). The
set $S(c) \cap {\rm cone}_c (V(c) \setminus \{v_0\})$ is spherically
convex in $S(c)$ by Proposition 4.1 (ii). The r.a.a. is accomplished by the following

\textbf{Claim:} $v_o$ \emph{is inessential in $V$, i.e.},
\begin{equation}\label{1b}
{\cal B}(V \setminus \{v_0\}) = {\cal B}(V)\,.
\end{equation}

\textbf{Proof:} We prove ``$\subset$'' in (\ref{1b}) (the
containment ``$\supset$'' is clear, as the operation $S \to
{\cal B}(S)$ reverses inclusion). Assume $x \in {\cal B}(V
\setminus \{v_o\})$. In order to prove that $x \in {\cal B}(V)$, we
use the equivalence $x \in {\cal B} (T) \Leftrightarrow T \subset
B (x)$, from which it follows that $x \in {\cal B}(V \setminus
\{v_0\}) \Leftrightarrow V \setminus \{v_0\} \subset B (x)$, and
$x \in {\cal B}(V) \Leftrightarrow V \subset B(x)$, where $B (x)$
is the unit ball centered at $x$. Hence
\begin{equation}\label{1c}
V \setminus \{v_0\} \subset B (x) \mbox{ and, in particular, } V (c) \setminus \{v_0\} \subset B (x)\,.
\end{equation}
By the first containment in (\ref{1c}) it
remains to show that $v_0 \in B(x)$. The set $B (x)$ is
spindle-convex. By Proposition 4.2, $S(c) \cap
B(x)$ is spherically convex in $S(c)$, hence by the second containment in (\ref{1c})

\begin{equation}\label{1d}
\mbox{sph-conv}(V(c) \setminus (v_0\}) \subset S (c) \cap B(x)\,,
\end{equation}
(the left-hand side is the spherical convex hull of $V(c) \setminus \{v_0\}$ in
$S(c)$).

By Proposition 4.1 (iii), sph-conv$(V(c) \setminus \{v_0\})
= S(c) \cap {\rm cone}_c (V(c) \setminus \{v_0\})$, hence by
(\ref{1d})
\begin{equation}\label{1e}
S(c) \cap {\rm cone}_c (V(c) \setminus \{v_0\}) \subset S(c) \cap
B(x)\,.
\end{equation}
It follows from (\ref{1a}) and (\ref{1e}) that
$v_0 \in S(c) \cap B(x) \subset B(x)$,
hence $v_0 \in B(x)$; this proves the claim.

As already said, the claim brings us to a
contradiction $(V$ is tight). \hfill \rule{2mm}{2mm}

\textbf{Lemma 6.2 (Spatial arm-lemma):} \emph{Let $C$ be a circle
in ${\mathbb R}^3$ centered at $c$, let $L$ be the line passing through $c$ perpendicular to {\rm aff}$\,C$,
and let $v \in {\mathbb R}^3 \setminus L$.
Define a function $f: C \to {\mathbb R}$ by
$f(p) = \|p-v\|$ for $p \in C$. Then}
\begin{enumerate}
\item[a)] $f$ \emph{has (exactly) two extremum points, one minimum
$p_{\rm min}$ and one maximum $p_{\rm max}$ obtained in the
following way: let $v'$ be the orthogonal projection of $v$ on the
plane {\rm aff}$C$ and denote by ${\rm pr}_c (\cdot): {\rm aff}C\to C$ the central
projection on $C$ with center $c$ in the plane} aff$C$. \emph{Then
$p_{\rm min} = {\rm pr}_c(v')$ and $p_{\rm max} = c + (c-p_{\rm min})=2c-p_{\rm \min}$
$($= the diametral point of $p_{\rm min}$ on the circle}
$C$),
\item[b)] $f(p)$ \emph{strictly increases as the point $p$
moves from $p_{\rm min}$ to $p_{\rm max}$ on a
half-circle of $C$ with endpoints $p_{\rm min}$ and} $p_{\rm
max}$. (\emph{There are two such halves}.)
\end{enumerate}

\textbf{Proof:} By $(f(p))^2 = \|p-v\|^2 = \|p-v'\|^2 + \|v'-v\|^2$ it is sufficient
to prove the lemma for the case where $v=v'$ (i.e., $v \in {\rm aff} C$).
Denote by $r$ the radius of $C$ and by $\alpha$ the small angle
between the vectors $v'-c \,\, (= v-c)$ and $p-c \,\, (0 \le \alpha
\le 180^o)$. Then
\begin{equation}\label{1f}
~\!\!\!\!\!\!\!\!\!\!\!\!\!\!\!(f(p))^2 = \|p-v'\|^2 =
\|p-c\|^2 + \| v'-c\|^2 - 2 \|p-c\| \cdot
\| v' - c \| \cos \alpha
\end{equation}
(cosinus theorem).

\vspace{1cm}

\begin{picture}(280,130)
\put(55,0){\includegraphics[scale=0.5]{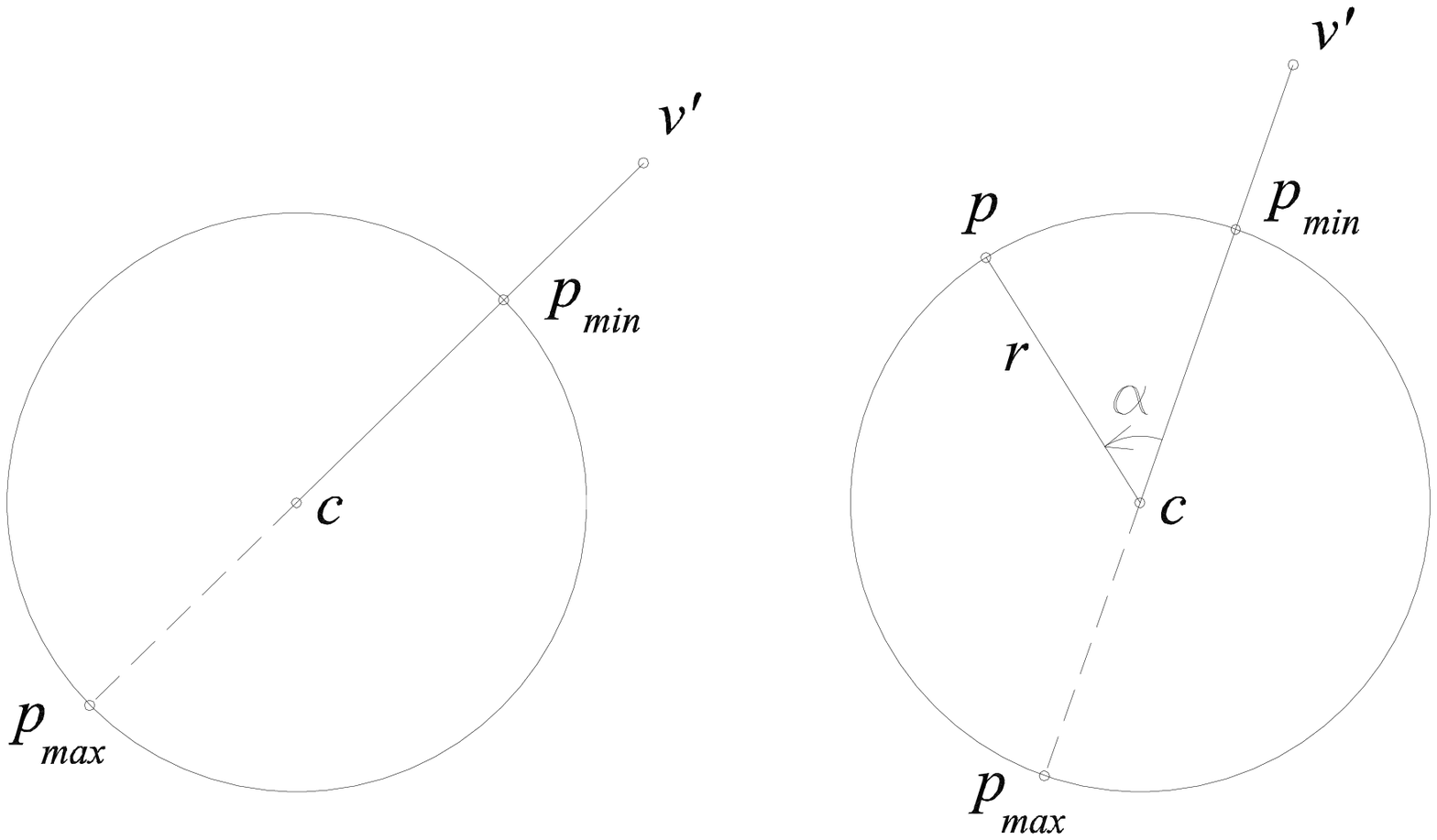}} 
\end{picture}

\begin{center}
-- Figures 6 and 7 --
\end{center}

As the point $p$ moves on one half of  the
circle $C$ from $p_{\rm min} \, (\alpha = 0^o)$ to $p_{\rm max}\,
(\alpha = 180^o)$, $\cos \alpha$ stricly decreases, and the
right-hand side of (\ref{1f}) increases. This proves a)
and b). \hfill \rule{2mm}{2mm}

\textbf{Lemma 6.3:} \emph{Let $C$ be a circle in ${\mathbb R}^3$
centered at $c$, and fix $p_0 \in C$. Let $H$ be the plane in
${\mathbb R}^3$ that passes through $p_0$ and $c$ perpendicularly
to} aff$C$. \emph{Let $H^+$ be an open halfspace
bounded by $H$ and put $C^+ =_{\rm def} C \cap H^+$. Then}
\begin{equation}\label{1g}
\begin{array}{l}
(\forall v \in H^+) (\exists \varepsilon > 0) (\forall p \in
C^+)\\
~\\
(\|p-p_0\| < \varepsilon \Rightarrow \|p-v\| < \|p_0 - v\|).
\end{array}
\end{equation}

\vspace{1.5cm}

\begin{picture}(280,130)
\put(50,0){\includegraphics[scale=0.6]{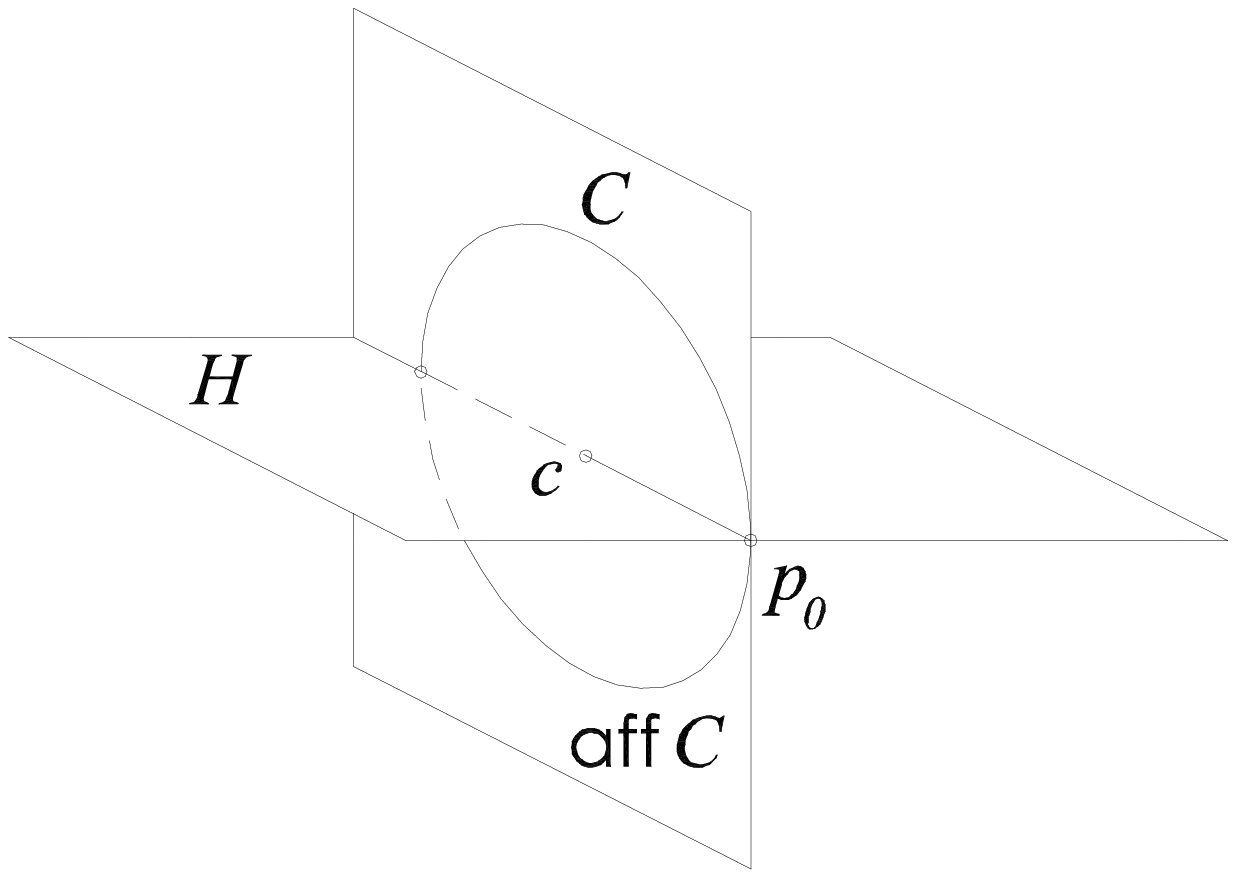}} 
\end{picture}
\vspace{-1cm}

\begin{center}
-- Figure 8 --
\end{center}

\textbf{Remark:} In Figure 8 the plane $H$ is ``horizontal'' and $H^+$ lies
``above'' it, while in Figure 9 below the plane $H$ is perpendicular to the plane of the figure.

\textbf{Proof:} Let $v'$ be the orthogonal projection of $v$ on
aff$C$. Since $v \in H^+$ and $H$ is orthogonal to aff$C$,
$v' \in {\rm aff} C \cap H^+$. By the foregoing Lemma 6.2 the function $f(p)
=_{\rm def} \|p-v\|, p \in C$, attains its minimum at $p_{\rm min}
= {\rm pr}_c(v')$ where pr$_c (\cdot): {\rm aff}C \to C$ is the central projection
on $C$ with center $c$ ($c \in {\rm aff}C$), and it attains its maximum at $p_{\rm max}
= c+(c-p_{\rm min})$.

\vspace*{1.5cm}

\begin{picture}(380,130)
\put(50,0){\includegraphics[scale=0.5]{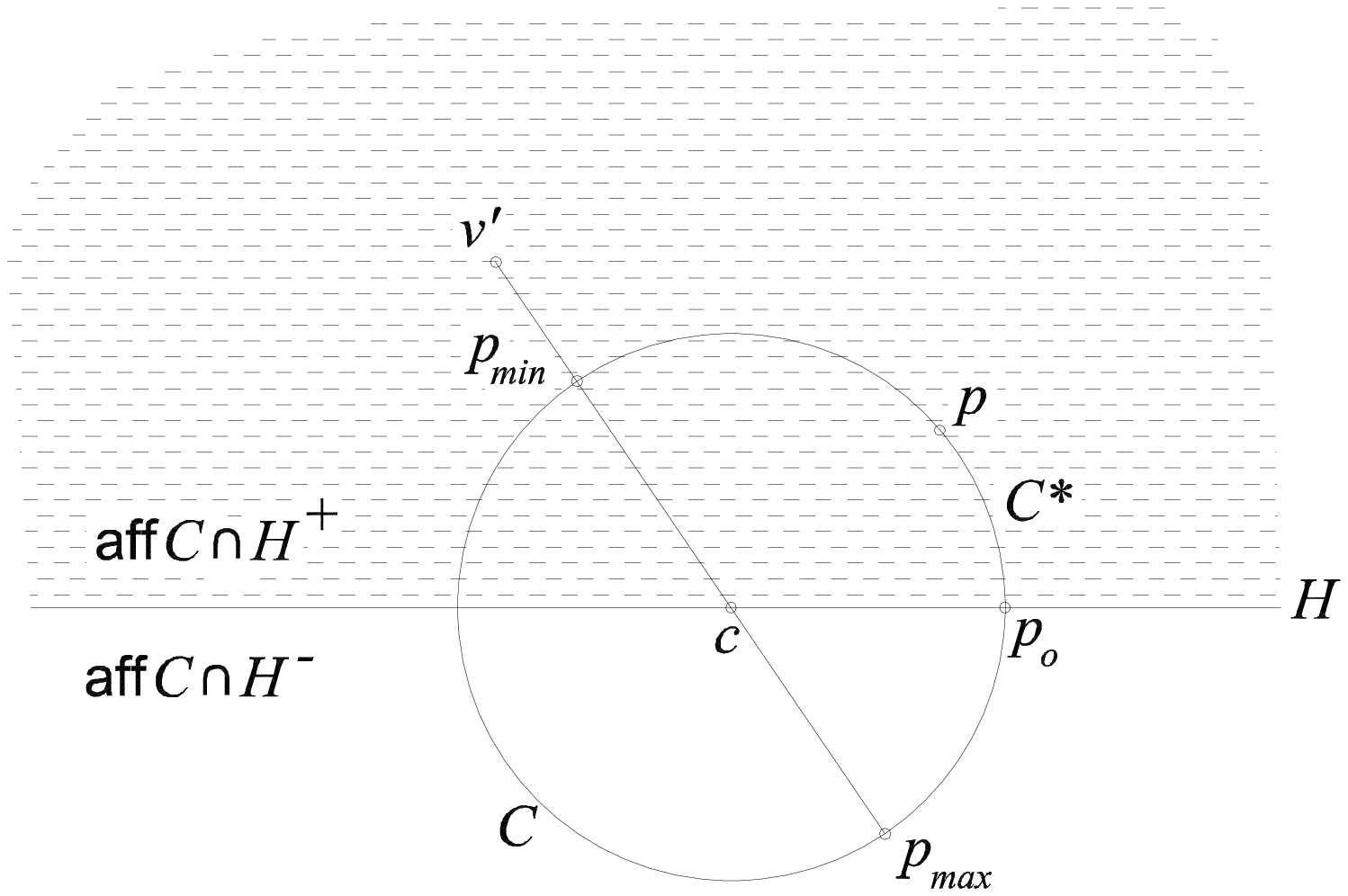}} 
\end{picture}

\begin{center}
-- Figure 9 --
\end{center}

Since $v' \in {\rm aff}C \cap H^+$,
$p_{\rm min} \in {\rm aff} C \cap H^+$ and $p_{\rm max} \in {\rm
aff}C \cap H^-$, where $H^-$ is the open half-space bounded by $H$
opposite to $H^+$. Denote by $C^*$ the half-circle of $C$ with
endpoints $p_{\rm min}$ and $p_{\rm max}$ which contains $p_0$
(note that $p_0 \notin \{p_{\rm min}, p_{\rm max}\}$ since $v \in
H^+$; hence $v \notin H$ and $v' \notin {\rm aff} \{p_0,c\}$). As
the point $p$ moves on $C^*$ starting at $p_{\rm min}$, the
distance $f(p) = \|p-v\|$ strictly increases (by Lemma 6.2). Hence
$f$, restricted to $C^+ \cap C^*$ (= the short arc on $C$ between
$p_{\rm min}$ and $p_0$), attains its maximum at $p_0$. Thus
$(\forall p \in C^+ \cap C^*) \,\, (\|p-v\| < \|p_0-v\|)$.
This proves (\ref{1g}) with $\varepsilon =_{\rm def}
\|p_{\rm min} - p_0\|$. ~\phantom{0} \hfill \rule{2mm}{2mm}

\textbf{Lemma 6.4:} \emph{With the assumptions and notation of}
Lemma 6.1 \emph{above assume that the rays} cone$_c (v_i) \, (i =
1,2, \dots, m,1)$ \emph{appear in this circular order on the
boundary of} cone$_c (v_1, \dots,v_m)$. \emph{I.e., for $m \ge 3$
the facets of this cone are} cone$_c (v_i, v_{i+1}) \,
(i=1,\dots,m)$ (\emph{the indices taken} modulo $m$, i.e., $v_{m+1} =
v_1$). (\emph{For $m=2$ the cone} cone$_c (v_1,v_2)$ \emph{is
$2$-dimensional.$)$ Then}:
\begin{itemize}
\item[a)] \emph{The facets of ${\cal B}(V)$ that meet at $c$ are
$F_{v_i},\, i=1,\dots,m$.}
\item[b)] \emph{These facets are
separated by edges $e_i$ with $e_i \subset F_{v_i} \cap
F_{v_{i+1}}\,, \, i=1,\dots,m$.}
\item[c)] \emph{For $m \ge 3$ and $1 \le i \le m$ the
edge $e_i$ initiates from $c$ on the circle $C_{v_iv_{i+1}} =_{\rm def}
S (v_i) \cap S(v_{i+1})$ to
the half space bounded by $H=_{\rm def} {\rm aff}(c,v_i,v_{i+1})$
which contains the set $\{v_1,v_2,\dots,v_{i-1},\hat{v}_i,
\hat{v}_{i+1}, v_{i+2},\dots,v_m\}$ $($and the tangent to $e_i$ at
$c$ is perpendicular to $H$$)$}.
\item[d)] \emph{For $m=2$
$($i.e., when $c$ is a dangling vertex$)$, there are exactly two
edges $e_1,e_2$ incident with $c$, both lying in $F_{v_1} \cap
F_{v_2} \subset C_{v_1v_2}$, and each edge $e_i \, (i=1,2)$
initiates from $c$ to a different side of the plane $H={\rm
aff}(c,v_1,v_2)$, $e_1$ and $e_2$ have a common tangent at $c$ which is the
tangent to $C_{v_1v_2}$ at $c$, and
this tangent is perpendicular to $H$.}
\end{itemize}

\textbf{Proof:} Put $V(c) =_{\rm def} (v_1, \dots, v_m)$ and
define
$\delta =_{\rm def} 1-\max \{\|v-c\| : v \in V \setminus V (c)\}$.
Then $0 < \delta < 1$ and
\begin{equation}\label{1h}
\forall x \in {\mathbb R}^3: \|x-c\| < \delta \Rightarrow x \in
{\rm int} \, {\cal B}(V \setminus V(c))\,.
\end{equation}
Fix $i \in \{1,\dots,m\}$, and put $C_i =_{\rm def} C_{v_i v_{i+1}}=
S(v_i) \cap S(v_{i+1})$ ($S (v_i)$ and $S(v_{i+1})$ are the
unit spheres centered at $v_i$ and $v_{i+1}$, respectively).
$C_i$ is a circle of radius $\sqrt{1-\frac{1}{4}
\|v_{i+1}-v_i\|^2}$ centered at $c_i =_{\rm def} \frac{1}{2} (v_i
+ v_{i+1})$, passing through $c$. Consider two cases:

\textbf{Case 1 (\boldmath$m=2$\unboldmath):} Then $\{v_i,
v_{i+1}\} = \{v_1,v_2\} = V(c)$, and cone$_c V(c)$ is
$2$-dimensional. Since $c$ is a vertex of ${\cal B}(V)$ that
belongs to exactly two facets of ${\cal B}(V)$, $c$
is a dangling vertex of ${\cal B}(V)$, hence $c \in V$.

By (\ref{1h}) there is an arc $A$ of length $> 2 \delta$ centered
at $c$ on the circle $C_i = C_1$ wholly contained in
${\cal B}(V \setminus V(c))$. Hence $F_{v_1} \cap F_{v_2} = C_1
\cap {\cal B}(V)$ has a non-degenerate connectivity component $e$
which is a closed arc of $C_1$ containing $c$ in its relative
interior. Since $c$ is a
dangling vertex of ${\cal B}(V)$, $e$ is the union of the two
edges of ${\cal B}(V)$ that meet at $c$, and $F_{v_1}$ and
$F_{v_2}$ are the only facets of ${\cal B}(V)$ that meet at $c$.

\textbf{Case 2 (\boldmath$m \ge 3$\unboldmath):} By Lemma 6.1 the
cone cone$_c (V(c))$ is full dimensional, and the rays cone$_c
(v_i), i = 1,\dots,m$, are extreme rays of it. Fix $i \in \{1,...,m\}$.
The plane aff$C_i$ bisects the segment $[v_i, v_{i+1}]$
perpendicularly at its midpoint denoted by $c_i$, hence it is
perpendicular to the plane $H_i =_{\rm def} {\rm aff} (v_i,
v_{i+1},c)$. The intersection of these two planes is the line
aff$(c_i,c)$ (this line is a perpendicular bisector of the
segment$[v_i,v_{i+1}]$; see Figure 10).

\vspace*{1.5cm}

\begin{picture}(280,130)
\put(50,0){\includegraphics[scale=0.6]{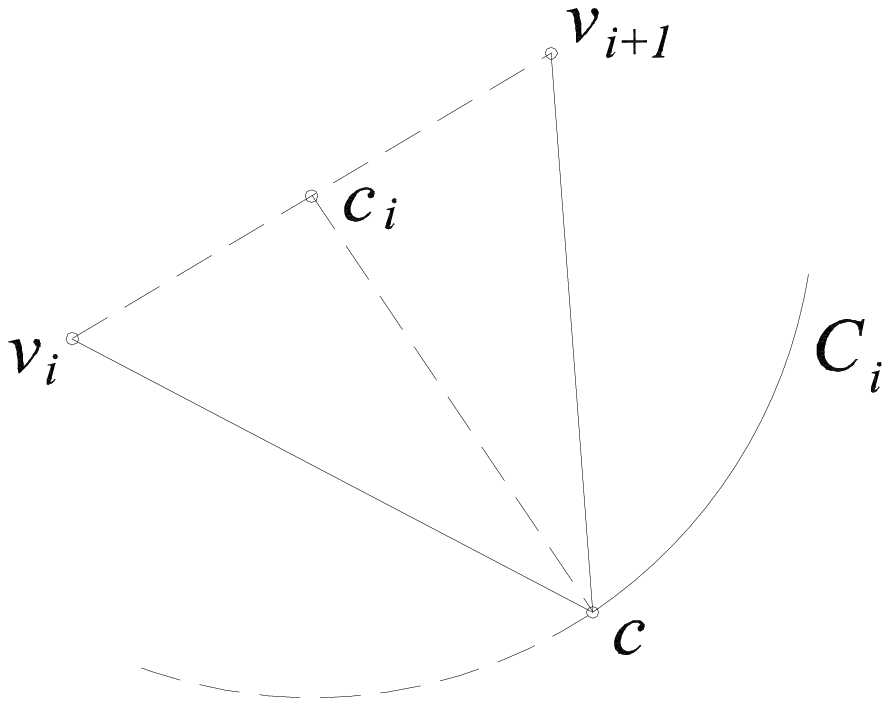}} 
\end{picture}

\vspace*{-0.3cm}
\begin{center}
-- Figure 10 --
\end{center}

Since the rays cone$_c(v_i)$ and cone$_c(v_{i+1})$ are extreme
rays of cone$_c V(c)$, cone$_c(v_i,v_{i+1})$ is a
($2$-dimensional) facet of cone$_c(V(c))$, hence the plane $H_i$
of this facet strictly supports $V(c) \setminus \{v_i,v_{i+1}\}$.

Denote by $H^+_i$, resp. $H^-_i$, the open halfspace bounded by
$H_i$ which contains, resp. does not contain, $V(c) \setminus
\{v_i,v_{i+1}\}$ and put $C^+_i = C_i \cap H^+$. By Lemma 6.3
(with $C = C_i, \, p_0 = c$ and $v = v_j$)
\begin{equation}\label{1i}
\begin{array}{l}
(\forall v_j \in V(c) \setminus \{v_i,v_{i+1}\}) (\exists
\varepsilon_j > 0) (\forall p \in C^+_i)\\[0.2cm]
((\|p-c\| < \varepsilon_j) \Rightarrow (\|p-v_j\| < 1))\,.
\end{array}
\end{equation}
Put $\varepsilon' =_{\rm def} \min \{\varepsilon_j: j \in
\{1,\dots,m\} \setminus \{i,i+1\}\}$ and $\varepsilon'' = _{\rm
def} \min \{\varepsilon',\delta\}$. Then (\ref{1h}) and
(\ref{1i}) imply
\begin{equation}\label{1j}
(\forall v \in V \setminus
\{v_i,v_{i+1}\}) (\forall p \in C^+_i)
((\|p-c\| < \varepsilon'') \Rightarrow (\|p-v\|< 1))\,.
\end{equation}
It follows that there is an arc $A$ of length $> \varepsilon''$ on
$C_i$ with one endpoint at $c$ such that
\begin{equation}\label{1k}
\mbox{relint} A \subset H^+_i \cap \mbox{ int } {\cal B}(V
\setminus \{v_i,v_{i+1}\})\,.
\end{equation}
Thus $C_i \cap {\rm bd} {\cal B}(V)$ has a component which
contains $A$, i.e., there is an edge $e_i$ of ${\cal B}(V)$ on
$C_i$ which is incident with $c$, and $e_i$ is common to $F_{v_i}$
and $F_{v_{i+1}}$. This proves a), and $e_i$ satisfies all the requirements
of b) and c), as it is easy to check. \hfill \rule{2mm}{2mm}

From Lemma 6.4 we have

\textbf{Corollary 6.1:} \emph{If $V \subset {\mathbb R}^2$ with $\# V \ge 3$ is
finite, tight, and ${\rm cr}(V) < 1$, then the
valence of a vertex $x$ of ${\cal B}(V)$ in the $1$-skeleton of
${\cal S}{\cal F}({\cal B}(V))$ is equal to the number of points
of $V$ at distance $1$ from $x$}\footnote{This corollary was observed (without a proof) by Gr\"unbaum, Heppes, and Straszewicz
independently  of each other for the case that $V$ is extremal (for the V\'{a}zsonyi problem). The corollary generalizes
their observation to any feasible set $V$ (i.e., $3 \le \# V < \infty$, $V$ is tight and satisfies cr$(V) < 1$).}.

\section{On \boldmath$\cB(V)$\unboldmath ~for an extremal set \boldmath$V$\unboldmath,
and an extended GHS-Theorem}

In this paragraph we use the foregoing results to
establish an extended form of the GHS-theorem. Recall that an
\emph{extremal configuration} (for the V\'{a}zsonyi problem) is a
finite set $V \subset {\mathbb R}^3$ with $\# V = n \ge 4$ such that
$e(V) = e(3,n)$.
We \textbf{ignore} now the GHS-Theorem ($e(3,n) = 2n-2$) in order
to achieve a detailed proof of an extended form of it. As said in the introduction, this is
essentially an elaboration of the proof given by Gr\"{u}nbaum, Heppes, and Straszewicz.
The following result was derived in Proposition and Definition 2.1 above as
a consequence of the GHS-Theorem. Now we derive it independently
of it.

\textbf{Proposition 7.1:} \emph{Let $V \subset {\mathbb R}^3\,, \,
\# V = n \ge 4$, be an extremal configuration for the V\'{a}zsonyi
problem. Then the valence of a point $v \in V$ in the diameter
graph $D(V)$ is $\ge 2$, i.e.}, \emph{if} diam$V=1$, \emph{then}
dist$(v,V \setminus \{v\})=1$, \emph{and this distance is attained
at least twice}.

\textbf{Proof:} Assume, without loss of generality, that diam$V
= 1$. If there is a point $v \in V$ such that dist$(v,V \setminus \{v\})
< 1$, then choose a line $l$ through $v$ and move $v$ on $l$ until
it reaches for the first time a point $v'$ such that dist$(v',V
\setminus \{v\}) = 1$. Put $V' =_{\rm def} (V \setminus \{v\})
\cup \{v'\}$. Then diam$V' = {\rm diam}V = 1$, $\# V' = \#
V = n$, and $e(V') > e (V)$, contrary to our
hypothesis (that $V$ is extremal).

A similar contradiction is obtained when dist$(v, V \setminus \{v\})
= 1$ and this distance is attained only once, i.e., if there is a
point $p \in V$ such that $\| v-p \| = 1$ and dist$(v, V \setminus
\{p,v\}) < 1$. Let $q$ be the symmetric point of $p$ relative to
$v$ on the line aff$(p,v)$, and let $C$ be a semi-circle of radius $1$
with center $v$ and endpoints $p,q$. Since $v$ is the center of
the segment $[p,q]$ whose length is $2$ and since dist$(p, V \setminus \{p,v\}) \le 1$,
clearly dist$(q, V \setminus \{p,v\})
> 1$.

Move $p$ on $C$ until it reaches for the first time a point $p'$
such that dist$(p', V \setminus \{p,v\}) = 1$. This will happen
before $p$ reaches $q$ since dist$(q, V \setminus \{p,v\})> 1$.
Put $V' =_{\rm def} (V \setminus \{p\}) \cup \{p'\}$. Then diam$V'
= {\rm diam}V = 1$, $\# V' = \# V$, and $e(V') > e(V)$, a
contradiction. Thus dist$(v, V \setminus \{v\}) = 1$, and this
distance is achieved at least twice. \hfill \rule{2mm}{2mm}

\textbf{Proposition 7.2:} \emph{Let $V \subset
{\mathbb R}^3$ with $\infty > \# V = n \ge 4$ be an extremal configuration
for the V\'{a}zsonyi problem such that} diam$V = 1$. \emph{Then}:
\begin{enumerate}
\item[a)] ${\cal B}(V)$ \emph{is full-dimensional and $V$ is
tight}.
\item[b)] \emph{ ${\cal B}(V)$ has $n$ facets}.
\item[c)] $V \subseteq {\rm vert} {\cal B}(V)$.
\end{enumerate}

\textbf{Proof:}
\begin{enumerate}
\item[(a)] The first part of a) follows from Theorem 5.2 a). The
second part of a) follows from Theorem 5.2~b) using Proposition 7.1
above.
\item[(b)] This follows from the fact that $V$ is tight.
\item[(c)] The valence of $v \in V$ in $D(V)$ is
$\ge 2$ for $v \in V$ (Proposition 7.1). If $v$ has valence
$2$ in $D(V)$, then it belongs to exactly two facets of ${\cal
B}(V)$, i.e., it is a dangling vertex of ${\cal B}(V)$. Otherwise
$v$ is of distance $1$ from at least three distinct points $p,q,r
\in V$, i.e., $v \in F_p \cap F_q \cap F_r$ and $v$ is a principal
vertex of ${\cal B}(V)$. \hfill \rule{2mm}{2mm}
\end{enumerate}

\textbf{Proposition 7.3:} \emph{Let $V \subset {\mathbb
R}^3$} \emph{be tight with} $\infty > \# V = n \ge 1$, cr$(V) < 1$, \emph{and assume that} $V
\subseteq {\rm vert} {\cal B}(V)$. \emph{Then} diam $V=1$,
\begin{equation}
e(V) \le 2n-2\,,
\end{equation} \emph{and this inequality is satisfied as an equality
iff} $V = {\rm vert} {\cal B}(V)$.

\textbf{Proof:} Since cr$(V) < 1$, ${\cal B}(V)$ is
full-dimensional. Since $V \subseteq {\rm vert} {\cal B}(V)$,
${\cal B}(V)$ has a non-empty set of vertices, hence $\# V = n \ge
3$.
The face numbers of ${\cal S}{\cal F}({\cal B}(V))$ satisfy the
Eulerian relation $v - e + f = 2$
(cf. Proposition 6.2). Here $v$ is the number of (principal or dangling) vertices
of $\cB(V)$, $e$ is the number of edges, and $f$ is the
number of facets which is a-priori $\le n$, but since $V$ is tight, $f=n$.
Since $V \subseteq {\rm vert}{\cal B}(V) \subset
{\cal B}(V)$, diam$V \le 1$. Each point $v \in V$, being a
vertex of ${\cal B}(V)$, is at distance $1$ from at least two
points of $V$, hence diam$V=1$.
Put $\widetilde{V} =_{\rm def} {\rm vert} {\cal B}(V)$. Thus $V
\subseteq \widetilde{V}$ and $\# \widetilde{V} = v = n + m$ for some $m
\ge 0$. The number $e$ of edges of ${\cal B}(V)$ is one-half the
sum of valences of the vertices of ${\cal B}(V)$ in the $1$-skeleton of
${\cal S}{\cal F}({\cal B}(V))$. The valence of
each vertex $x$ of ${\cal B}(V)$, henceforth denoted by val$(x)$, is
equal to the number of points of $V$ at distance 1 from $x$ (cf.
Corollary 6.1 in the end of \S~6).
Thus
\[
\begin{array}{lll}
2e & = & \sum \{{\rm val} (x) : x \in \widetilde{V}\}=\\
& = & \sum \{{\rm val} (x) : x \in V\} + \sum \{{\rm val} (x) : x
\in \widetilde{V} \setminus V\}\\
& \ge & \sum \{{\rm val} (x) : x \in V\} + 3 \cdot \#
(\widetilde{V} \setminus V)\\
& = & 2 \cdot \# \{\{u,w\} \subset V : \| u-w\| = 1\} + 3m\\
& = & 2 \cdot e(V) + 3 m\,.
\end{array}
\]
It follows, using the Eulerian relation $v-e+f=2$, that
\[
\begin{array}{lll}
 2 \cdot e (V) & \le & 2e-3m\\
& = &  2 (v+f-2)-3m \quad \mbox{(by} \,  v-e+f=2)\\
& = & 2 (n+m+f-2) - 3m\\
& = & 2 (2n-2)-m \quad \mbox{(since} \, f=n)\,.
\end{array}
\]
Thus $e(V) \le 2n-2$,
and this is satisfied as an equality iff $m = 0$, i.e.,
iff $V = {\rm vert} {\cal B}(V)$. \hfill
\rule{2mm}{2mm}

This leads to

\textbf{Theorem 7.1 (Extended GHS-Theorem):} \emph{Let $V \subset
{\mathbb R}^3$ be finite with $\# V = n \ge 4$ and} diam$V = 1$.
\emph{The following three statements are equivalent}:
\begin{enumerate}
\item[(i)] $V$ \emph{is extremal for the V\'{a}zsonyi
problem, i.e.}, $e(V) = e(3,n)$;
\item[(ii)] $e(V) = 2n-2$;
\item[(iii)] $V$ \emph{is tight and} $V = {\rm vert} {\cal B}(V)$\footnote{
(a) The equivalence of statements (i) and (ii) is what we call the
GHS-Theorem (see Corollary 7.1 below), and their equivalence to (iii) is what ``Extended'' stands for in the name of the
theorem.\\
(b) Sallee in \cite{Sa}, p. 318, calls a Ball polytope satisfying (iii) ``a frame of a Reuleaux polytope'',
and he points on their usefulness to construct concrete bodies of constant width in $\bR^3$.
This is analogous to the planar case, where Reuleaux polygons are paradigmatic for (planar) bodies of constant width.
(He does not point on
the connection of his ``frames of Reuleaux polytopes''
with the V\'azsonyi problem, however.) We may come to this topic in the mentioned subsequent paper.}.

\end{enumerate}

\textbf{Proof:} diam$V =1 \Rightarrow {\rm cr}(V) < 1$ (Theorem 5.2, a)), hence $\cB(V)$ is
full-dimensional.

\textbf{Proof of (i) $\Rightarrow$ (ii) and (i) $\Rightarrow$
(iii):} Assume (i). By Proposition 7.2 a), c) the set $V$ is tight, and $V \subseteq {\rm vert}
{\cal B}(V)$. Hence, by Proposition 7.3, $e(V) \le 2n-2$. On the
other hand, $e(V) = e(3,n) \ge 2n-2$, by the examples given in
\S~1. Thus $e(V) = 2n-2$ and, again by the suffix to Proposition
7.3, $V = {\rm vert} {\cal B}(V)$.

\textbf{Proof of (ii) $\Rightarrow$ (i):} By (i) $\Rightarrow$ (ii) proved above $e(3,n) = 2n-2$.
Hence $e(V) =e(3,n)$, i.e., $V$ is extremal.

\textbf{Proof of (iii) $\Rightarrow$ (ii):} Assume that $V$ is tight and $V = {\rm vert}
{\cal B}(V)$. Thus $V \subseteq {\rm vert} {\cal B}(V)$
and the conditions of Proposition 7.3 are fulfilled. By the suffix to this proposition
$V = {\rm vert}\cB(V) \Rightarrow e(V) = 2n-2$.
\hfill \rule{2mm}{2mm}

From Theorem 7.1 follows

\textbf{Corollary 7.1 (GHS-Theorem):} $e(3,n) = 2n-2$.

\textbf{Corollary 7.2:} \emph{Assume that $V \subset {\mathbb R}^3, 4 \le \# V < \infty$}, diam$V=1$
\emph{and $V$ is extremal for the V\'{a}zsonyi problem. Then every edge of ${\cal B}(V)$ is a} short
\emph{arc of its loading circle}.

\textbf{Proof:} Let $e$ be an edge of ${\cal B}(V)$ with endpoints $a,b \in {\rm vert}{\cal B}(V) = V$.
There are points $x,y \in V$ such that $e \subset C_{xy} =_{\rm def}S_x \cap S_y$ ($S_x,S_y$
are the spheres of radius $1$ centered at $x$ and $y$, respectively). $C_{xy}$ is a circle of
radius $r_{xy} = \sqrt{1-\frac{1}{4}\|x-y\|^2}$ centered at $p=\frac{1}{2}(x+y)$, $a,b \in C_{xy}$, and
$e$ is either the short arc between $a$ and $b$ on $C_{xy}$, or the respective long arc. We will
show that the second case is impossible.

\vspace{1cm}

\begin{picture}(280,130)
\put(53,0){\includegraphics[scale=0.6]{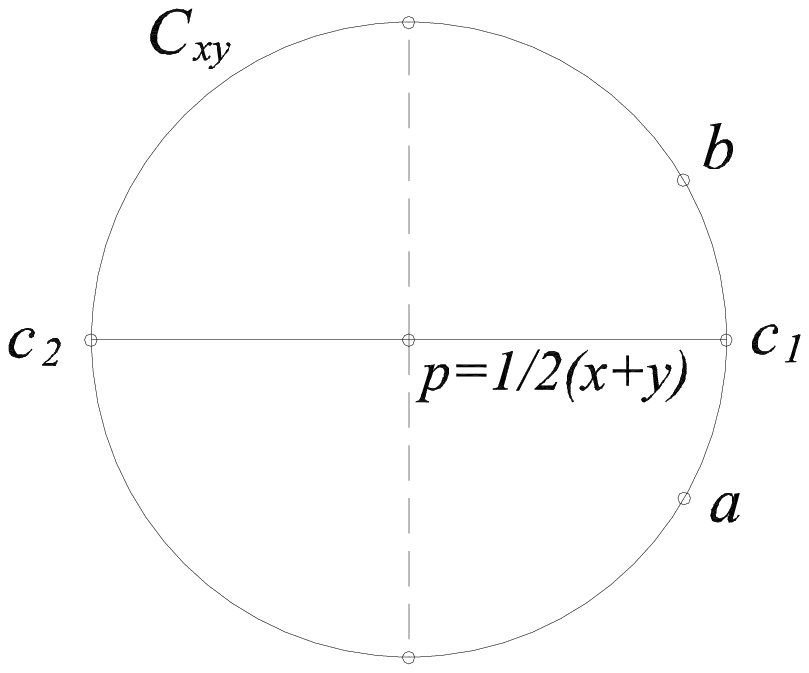}} 
\end{picture}

\vspace{-1cm}
\begin{center}
-- Figure 11 --
\end{center}

Denote by $c_1$ [resp. $c_2$] the midpoint of the short [resp. long]
arc between $a$ and $b$ on $C_{xy}$ (cf. Figure 11). $c_1$ and $c_2$ are diametrical points of $C_{xy}$. Since $\|x-y\| \le 1$,
$r_{xy} = \sqrt{1-\frac{1}{4}\|x-y\|^2} \ge \frac{\sqrt{3}}{2} > \frac{1.7}{2} = 0.85$; hence
diam$C_{xy} = 2 r_{xy} \ge 1.7 > 1$. Since $c_2$ is the midpoint of the long arc between $a$ and
$b$ on $C_{xy}$, clearly $\wk c_2pa > 90^o$.
It follows that
\[
\begin{array}{lll}
\|c_2-a\|^2 & = & \|c_2-p\|^2+\|a-p\|^2-2\|c_2-p\| \cdot \|a-p\| \cdot \cos \wk c_2 pa\\
& > & \|c_2-p\|^2+\|a-p\|^2=2(r_{xy})^2 \ge 2 \cdot \left(\frac{\sqrt{3}}{2}\right)^2 = \frac{3}{2} > 1\,.
\end{array}
\]
Hence $\|c_2-a\|> 1$ and $c_2 \notin B(a,1)$.
Since $a \in {\rm vert}\cB(V)=V$ (by Theorem 7.1 (iii)), this implies $c_2 \notin {\cal B}(V)$, hence $c_2 \notin e$, i.e., $e$
is \emph{not} the long arc between $a$ and $b$ on $C_{xy}$. Thus $e$ is the short arc between
$a$ and $b$ on $C_{xy}$. \hfill \rule{2mm}{2mm}

\section{The canonical self-duality \boldmath$\varphi$\unboldmath~ of~ \boldmath${\cal S}
{\cal F}({\cal B}(V))$\unboldmath\\ for an extremal set \boldmath$V$\unboldmath}

In this paragraph we assume that $V$ is a set of $n$ points in ${\mathbb R}^3\,, \, 4 \le n
< \infty$, diam$V = 1$, and $V$ is extremal for the
V\'{a}zsonyi problem. Under these assumptions $V$ coincides with
the set of vertices of its Ball polytope ${\cal B}(V)$ (Extended GHS-Theorem, Theorem 7.1 (iii) above). The poset
${\cal S}{\cal F}({\cal B}(V))$, partially ordered by inclusion,
consists of vertices, edges and facets of $\cB(V)$ (and, if one wishes, the
improper faces $\emptyset$ and ${\cal B}(V)$ as well).

\textbf{Definition 8.1:} An \emph{involutory self-duality} $($or \emph{antimorphism}$)$ of
${\cal S}{\cal F}({\cal B}(V))$ is an order reversing map (i.e., a
duality) $\varphi: {\cal S}{\cal F} ({\cal B}(V)) \longrightarrow {\cal
S}{\cal F}({\cal B}(V))$ of order two $(\varphi^2 =$ identity).

The correspondence $\{x\}
\stackrel{\varphi}{\leftrightarrow} F_x \,\, (x \in V)$ (and, if
necessary, $\emptyset \stackrel{\varphi}{\leftrightarrow} {\cal
B}(V))$ is, in fact, order reversing, i.e.,
\[
\begin{array}{lll}
\{x\} \subset F_y \Leftrightarrow \| x-y \| = 1 &
\Leftrightarrow & \{y\} \subset F_x\\
& \Leftrightarrow & \varphi (F_y) \subset \varphi (\{x\})\,.
\end{array}
\]
And it is, of course, of order two.
This correspondence can be extended to an involutory
self-duality of ${\cal S}{\cal F}({\cal B}(V))$
if and only if the following statement holds.

\textbf{Theorem 8.1 (dual edges):} \emph{Let $e$ be an edge of
${\cal B}(V)$, with endpoints $a$ and $b$, and assume that $e$ is included in
the facets $F_x$ and $F_y$ of ${\cal B}(V)$ $(x,y \in V, x \not= y)$. Then
there is another edge $e'$ of $\cB(V)$, with endpoints $x$ and $y$, included in the facets
$F_a$ and $F_b$. The edge $e'$ is, of course, uniquely
determined by $e$}.

This theorem is indeed true. Its proof takes some effort,
however, and will be postponed (until after Lemma 8.1 below). From this theorem we conclude

\textbf{Corollary -- Definition 8.1 (canonical
self-duality \boldmath$\varphi$\unboldmath~of \boldmath${\cal
S}{\cal F} ({\cal B}(V))\unboldmath$:}

\emph{If $V \subset {\mathbb R}^3$, with $\# V \ge 4$ and}
diam$(V) = 1$, \emph{is an extremal configuration for the V\'{a}zsonyi
problem, then the correspondence $x
\stackrel{\varphi}{\leftrightarrow} F_x \, (x \in V)$ has a unique
edge-extension $\widetilde{\varphi}, e
\stackrel{\widetilde{\varphi}}{\leftrightarrow} e'$, which is an
involutory self-duality of ${\cal S}{\cal F} ({\cal B}(V))$. This
self-duality $\widetilde{\varphi}$ is the} canonical
self-duality \emph{of ${\cal S}{\cal F}({\cal B}(V))$. Abusing notation
we drop henceforth the tilde $\,\widetilde{~}$ and denote $\widetilde{\varphi}$ by
$\varphi$}\footnote{This self-duality was partially
observed by Sallee in \cite{Sa}, p. 322, Remark 2). Our $\cB(V)$ is called there
``a frame of a Reuleaux polytope'' (cf. footnote 20 above), and something like our
Theorem 8.1 above is formulated in (6.2) in p. 320 there. His treatment is sketchy, however,
overlooking some delicate points.}~\footnote{Also, this corollary turns out to be a partially special case of the following
remarkable duality theorem of \cite{Bez-Nasz}, p. 256: If $\cB(V)$ is a ball polytope without digonal facets (in $\bR^3$), then
$\cB(V)$ and $\cB$(vert$\cB(V)$) have dual face structures. The proof given there (pp. 261-262) is un-convincing. To see this, 
it is sufficient to point out that even such most basic
properties as enunciated in Proposition 6.3 (iii) and Lemmata 6.1, 6.3, and 6.4 above are nowhere mentioned there, without which
a rigorous treatment of the face structure of $\cB(V)$ is impossible. Also the treatment of the edge-edge duality in
p. 262 there is un-convincing merely because of the lack of the important geometric inequality enunciated in Lemma 8.1
below (pointwise distance between dual arcs), which takes care of the delicate possibility that what is called the ``desired edge''
(p. 262, line 16) contains no point of $X$ (= our $V$) in its relative interior. Also dual edges of $\cB(V)$ and $\cB$(vert$\cB(V)$) are
both either short or long (as arcs of their corresponding loading circles) -- which follows from Proposition 6,4 and Lemma 6.4 above, a fact not
mentioned there. Let us also mention that by defining a variant of our dangling vertices (cf. Definition 6.1 above) one can
relieve this Bezdek-Nasz\'{o}di duality theorem
from the manifestly artificial assumption that $\cB(V)$ contains no digonal facets.}.

The proof of Theorem 8.1 makes use of a geometric fact for which
we need the following

\textbf{Definition 8.2 (dual arcs):} Assume that $a,b,x,y$ are
four distinct points in ${\mathbb R}^3$, such that $\| x-a \| = \|x-b \| = \| y-a \| = \| y-b \| = 1$.
The (spatial) quadrilateral ($a,x,b,y)$ is an \emph{equilateral quadrilateral}, and
the couple of $2$-tuples $(a,b;x,y)$ are in \emph{equilateral
position}. Denote by $C_{xy}$, resp. $C_{ab}$, the circle centered at $\frac{1}{2}
(x+y)$, resp. $\frac{1}{2}(a+b)$, that passes through $a$ and $b$, resp. $x$ and $y$.
Let $A_{ab}$ [resp. $A_{xy}$] be the short open circular arc of
$C_{xy}$ [resp. $C_{ab}$] with endpoints $a,b$ [resp. $x,y$]. The
circular arcs $A_{ab}$ and $A_{xy}$ are \emph{short dual arcs}. The complementary arcs
$C_{xy} \setminus ({\rm relint}A_{a,b})$ and $C_{ab} \setminus ({\rm relint}A_{xy})$ are
\emph{long dual arcs}.

Note that the radii $r_{ab}$ and $r_{xy}$ of the circles $C_{ab}$ and $C_{xy}$, respectively, are both $< 1$
and, of course, these radii may be different from each other. Note also that the
planes aff$C_{xy}$ and aff$C_{ab}$ are perpendicular to each other.

The following lemma partially appears as Lemma 4.2 in \cite{Sch-Pe-Ma-Ku}, with an analytic proof. We extend it
here and give a purely geometric proof.

\textbf{Lemma 8.1 (pointwise distance between dual arcs and their
complements):} \emph{Let $A_{ab}$ and $A_{xy}$ be two (circular) short dual arcs, lying on the
circles $C_{xy}$ and $C_{ab}$,
respectively, and let $c \in C_{xy}$ and $z \in C_{ab}$}.
\begin{itemize}
\item[a)]\emph{If $c \in {\rm relint}A_{ab}$ $[$resp. $c \in C_{xy} \setminus A_{ab}$$]$ and
$z \in {\rm relint}A_{xy}$, then $\|c-z\|>1$ $[$resp. $\|c-z\|<1$$]$}.
\item[b)]\emph{If $c \in C_{xy} \setminus A_{ab}$ and $z \in C_{ab} \setminus A_{xy}$, then
$\|c-z\|>1$}\footnote{This lemma is fundamental for ball polytopes. E.g., it has to be used in any rigorous proof of the fundamental
Bezdek-Nasz\'{o}di duality theorem (cf. the former footnote), especially for the edge-edge duality, alas it is lacked in \cite{Bez-Nasz}, p. 262,
which is one of the reasons we find this proof un-convincing.}.
\end{itemize}

In other words: If $c \in C_{xy}$ resp. $z \in C_{ab}$, and if either both $c$ resp. $z$ lie on the relative
interiors of the short [resp. long] arcs of $C_{xy}$ resp. $C_{ab}$, with endpoints $a,b$ resp. $x,y$, then
$\|c-z\|>1$. Otherwise $\|c-z\| \le 1$, and this is satisfied as an equality iff $c \in \{a,b\}$ or $z \in \{x,y\}$
(or both).

\vspace{1cm}

\begin{picture}(280,130)
\put(60,0){\includegraphics[scale=0.6]{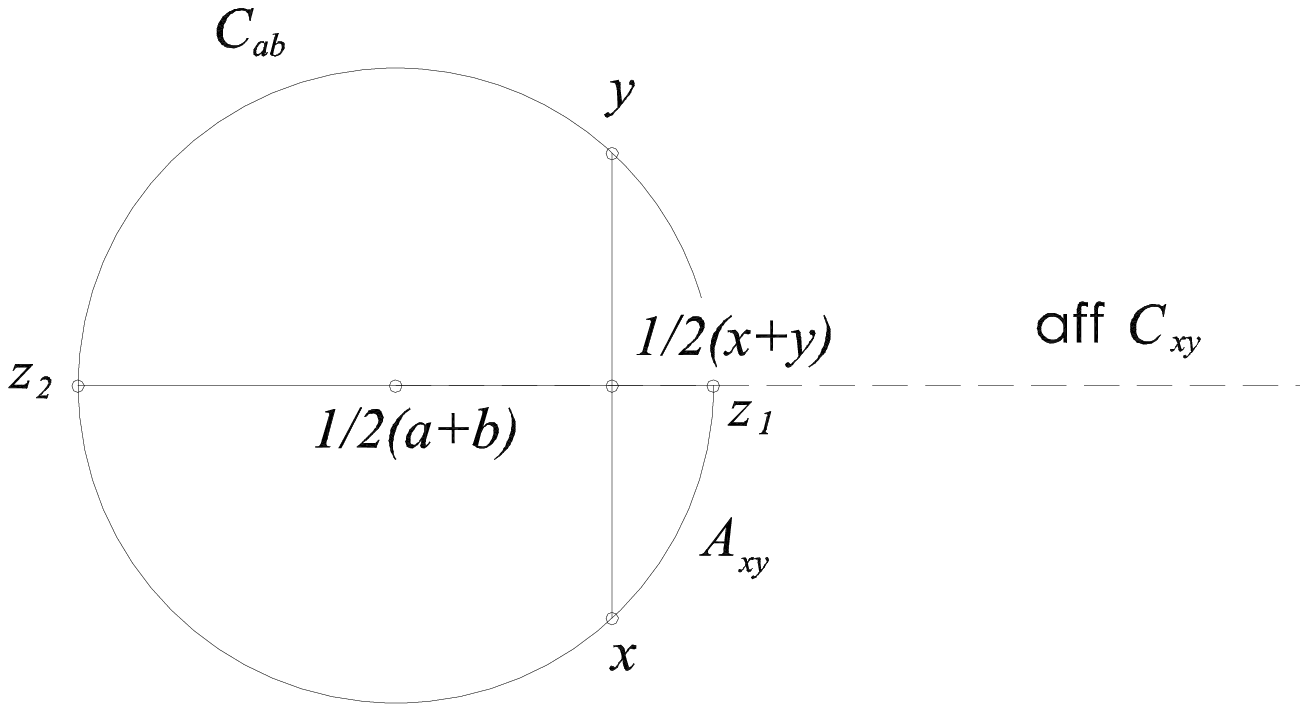}} 
\end{picture}

\vspace{-1cm}
\begin{center}
-- Figure 12 --
\end{center}

\textbf{Proof:} The line aff$(\frac{1}{2}(a+b), \frac{1}{2}(x+y))$
passing through the center $\frac{1}{2}(a+b)$ of the circle
$C_{ab}$ and through the midpoint $\frac{1}{2}(x+y)$ of its chord
$[x,y]$ intersects the circle $C_{ab}$ in two diametrical points
$z_1$ and $z_2$; we choose the notation so that $z_1 \in A_{xy}$ (in
fact, $z_1$ is the midpoint of the arc $A_{xy}$) and $z_2 \in
C_{ab} \setminus A_{xy}$ ($z_2$ is the midpoint of this
long arc of $C_{ab}$). The planes $H_{xy} =_{\rm def} {\rm aff}
C_{xy}$ and $H_{ab} =_{\rm def} {\rm aff} C_{ab}$ intersect
perpendicularly in the line aff$(z_1,z_2)$, and $\frac{1}{2}(a+b)$
is the midpoint of both segments $[a,b]$ and $[z_1,z_2]$ (cf. Figure 12).

\textbf{Claim:} aff$(a,b) \perp {\rm aff}(z_1,z_2)$

\textbf{Proof:} Since $z_1,z_2 \in {\rm aff} \{\frac{1}{2} (x+y), \frac{1}{2}(a+b)\}$, it
is sufficient to show that $\langle a-b,
\frac{1}{2}(x+y)-\frac{1}{2}(a+b)\rangle = 0$. This is equivalent to
\[
\begin{array}{llll}
&                 &   & \langle a-b, x-a+y-b\rangle = 0\\
& \Leftrightarrow &   & \langle 2a-2b, x-a+y-b\rangle = 0\\
& \Leftrightarrow &   & \langle a-x+a-y+x-b+y-b, x-a+y-b\rangle = 0\\
& \Leftrightarrow &   & \langle a-x,x-a\rangle + \langle y-b,y-b\rangle + \langle a-x,y-b\rangle + \langle y-b,x-a\rangle\\
&                 & + & \langle a-y+x-b, x-a+y-b\rangle = 0\\
& \Leftrightarrow & - & 1+1 +  \langle y-b,a-x+x-a\rangle\\
&                 & + & \langle a-y+x-b, x-b+y-a\rangle = 0\\
& \Leftrightarrow &   & \langle a-y+x-b,x-b+y-a\rangle = 0\\
& \Leftrightarrow &   & \langle a-y,y-a\rangle + \langle x-b,x-b\rangle\\
&                 & + & \langle a-y,x-b\rangle + \langle x-b, y-a\rangle = 0\\
& \Leftrightarrow & - & 1+1+\langle x-b,a-y+y-a\rangle = 0\\
& \Leftrightarrow &   & \langle x-b,0\rangle = 0 \Leftrightarrow 0 = 0\,,
\end{array}
\]
which proves the claim.

Let pr$(\cdot) : \bR^3 \to H_{ab}$ be the orthogonal projection on the plane $H_{ab}$. Since
$H_{ab} \perp H_{xy}$,
\begin{equation}\label{21}
{\rm pr}(C_{xy}) \subset H_{ab} \cap H_{xy} = {\rm aff}(z_1,z_2)\,,
\end{equation}
and by the foregoing claim
\begin{equation}\label{22}
{\rm pr}(a) = {\rm pr}(b) = {\small \frac{1}{2}} (a+b)\,.
\end{equation}
Since $A_{ab}$ is a short arc on $C_{xy}$, relint$A_{ab}$ lies beyond the line aff$(a,b)$ relatively to
the center $\frac{1}{2}(x+y)$ of $C_{xy}$ in the plane $H_{xy}$. Hence, by (\ref{22}),
\begin{equation}\label{22a}
{\rm pr}(A_{ab}) \subset-{\rm cone}_{\small \frac{1}{2}(a+b)} {\small \frac{1}{2}}(x+y)
\end{equation}
and
\begin{equation}\label{22b}
{\rm pr}(C_{xy} \setminus A_{ab}) \subset {\rm cone}_{\frac{1}{2}(a+b)} {\small \frac{1}{2}} (x+y)\,.
\end{equation}
Similarly, since $A_{xy}$ is a short arc on $C_{ab}$, relint$A_{xy}$ lies beyond the line aff$(x,y)$
relatively to the center $\frac{1}{2}(a+b)$ of $C_{ab}$ in the plane $H_{ab}$. Hence, in particular,
\begin{equation}\label{22c}
z_1 \in {\rm cone}_{\frac{1}{2}(a+b)} {\small \frac{1}{2}}(x+y)
\end{equation}
and
\begin{equation}\label{22d}
z_2 \in {\rm cone}_{\frac{1}{2}(a+b)} {\small \frac{1}{2}} (x+y)\,.
\end{equation}
It follows from (\ref{22a}), (\ref{22b}), (\ref{22c}), and (\ref{22d}) that
\begin{equation}\label{23a}
{\rm pr}(A_{ab}) \subset{\rm cone}_{\frac{1}{2}(a+b)} z_2
\end{equation}
and
\begin{equation}\label{23b}
{\rm pr}(C_{xy} \setminus A_{ab}) \subset {\rm cone}_{\frac{1}{2}(a+b)} z_1\,.
\end{equation}

Back to the proof of Lemma 8.1, let $c \in {\rm relint}A_{ab}$
[resp. $c \in C_{xy} \setminus A_{ab}]$. Then, by (\ref{23a}) and (\ref{23b}),
$c' =_{\rm def} {\rm pr}(c)$ lies on
cone$_{\frac{1}{2}(a+b)} z_2$ [resp.
cone$_{\frac{1}{2}(a+b)}z_1$]. Hence, by Lemma 6.2 (Spatial
arm-lemma), $f:C_{ab} \to \bR$, defined by $f(z)=_{\rm def} \|z-c\| \, (z \in C_{ab})$, attains its
minimum at $z_2$ [resp. $z_1$] and its maximum at $z_1$ [resp. $z_2$], and as
the point $z$ moves on half of
the circle $C_{ab}$ from $z_2$ to $z_1$, $f(z)$ increases [resp.
decreases]. Since for $z=x$ we have $f(x) = \|x-c\|=f(y) = \|y-c\|=1$
(recall that $c \in C_{xy}$ whether $c \in A_{ab}$ or $c \in
C_{xy} \setminus A_{ah}$) as $z$ moves from $z_2$ to $z_1$ on half
of the circle $C_{ab}$ that contains $x$ we have by increasing [resp. decreasing]
monotonicity of $f(z)$ that $f(z) = \|z-c\|<\|x-c\|=1$
[resp. $f(z) = \|z-c\|>\|x-c\|=1$] when $z$ lies on the short arc
arc$(z_2,x)$ of $C_{ab}$, and $f(z)=\|z-c\|>\|x-c\|=1$ [resp. $f(z)
= \|z-c\|<\|x-c\|=1$] when $z$ lies on the short arc arc$(x,z_1)$ of
$C_{ab}$. Similar inequalities hold as $z$ moves from $z_2$ to
$z_1$ on half of the circle $C_{ab}$ that contains $y$. This
proves Lemma 8.1 a).

To prove b), assume that the points $c,z$ are as stated in the
premises. Now we have $c'=_{\rm def} {\rm pr}(c) \in {\rm
cone}_{\frac{1}{2}(a+b)} z_1$, and by Lemma 6.2,
as $z$ moves from $z_1$ to $z_2$ on half of the circle $C_{ab}$
that contains $x$, $f(z) =_{\rm def}\|z-c\| \, (z \in C_{ab})$
increases, and since it reaches the value $1$ at $z=x \,
(\|x-c\|=1$, since $c \in C_{xy}$), $f(z) = \|z-c\|>\|x-c\|=1$ when
$z$ lies on the short arc arc$(x,z_2)$ of $C_{ab}$. A similar
inequality holds when $z$ moves from $z_1$ to $z_2$ on half of the
circle $C_{ab}$ that contains $y$. This shows that $\|z-c\|>1$ for
$z \in C_{ab} \setminus A_{xy}$ (and $c \in C_{xy} \setminus
A_{ab}$). \hfill \rule{2mm}{2mm}

\textbf{Proof of Theorem 8.1:} The edge $e$ (whose vertices are
$a,b$) is a short arc of the circle $C_{xy}$, by Corollary 7.2
above. Hence $e=A_{ab}$, where $A_{ab}$ and $A_{xy}$ are the short
circular (dual) arcs defined for the couple of $2$-tuples
$(a,b;x,y)$ in Definition 8.2 above. Put $H=_{\rm def} {\rm
aff}(x,a,b)$ and let $H^+$ be the closed half space bounded by $H$ which
contains $y$. Since $A_{xy}$ is a short arc of $C_{ab}$, clearly
$A_{xy} \subset H^+$.

\textbf{Claim 1:} (relint$A_{xy}$) $\cap\, V = \emptyset$.

\textbf{Proof:} If $A_{xy}$ contains a point $v \in V$ in its
relative interior, then by Lemma 8.1~a) above $B(v,1) \cap {\rm
relint}A_{ab} = \emptyset$, hence $e = A_{ab} \nsubseteq {\cal
B}(V)$, which is a contradiction ($e$ is an edge of ${\cal
B}(V)$). This proves Claim 1.

\textbf{Claim 2:} \emph{There is an edge of ${\cal B}(V)$ incident with
$x$ whose initial arc from $x$ lies on $A_{xy}$}.

\textbf{Proof:} Put $m=_{\rm def} \#({\rm vert}F_x$). If $m=2$,
i.e., vert$F_x = \{a,b\}$, then by Lemma 6.4~d) there are two
edges of ${\cal B}(V)$ incident with $x$ which have initial arcs
from $x$ on $C_{ab}$, lying in different sides of $H$; hence one
of them has an initial arc from $x$ lying on $A_{xy}(\subset
C_{ab})$ contained in $H^+$. If $m \ge 3$ then, by Proposition 6.4 (iv), the point
$y$ and (vert$F_x) \setminus \{a,b\}$ lie on the same side of the plane $H$,
and since $y \in H^+$, (vert$F_x) \setminus \{a,b\} \subset H^+$.
By Lemma 6.4~c) there is an edge $e'$ of ${\cal B}(V)$ that has in
initial arc from $x$ lying on $C_{ab}$ contained in $H^+$, i.e.,
the initial arc of $e'$ from $x$ is contained in $A_{xy}$. This
proves Claim 2.

\textbf{Claim 3:} $A_{xy} = e'$.

\textbf{Proof:} $e'$ has an initial arc from $x$ lying on $A_{xy}$
(Claim 2), hence its other endpoint (in int$H^+$) is either $y$ or a
point $y' \in {\rm relint}A_{xy}$. The second case is impossible
by Claim 1. This proves Claim 3, concluding the proof of Theorem 8.1.
~\phantom{00} \hfill \rule{2mm}{2mm}

The following lemma is needed in the proof of Theorem 8.2 below.

\textbf{Lemma 8.2:} \emph{Let $V \subset \bR^3$ be an extremal set for the Vazs\'{o}nyi problem with $\infty > \#V = n \ge 4$, and
$e,f$ be two dual edges of $\cB(V)$ with endpoints $(x,y)$ and $(a,b)$, respectively, and let
$u \in {\rm relint}e$. Then}:
\begin{enumerate}
\item[a)] \emph{$V \cup \{u\}$ is an extremal set}.
\item[b)] \emph{All edges of $\cB(V)$ are edges of $\cB(V \cup \{u\})$ and vice versa}, except
\emph{for $e$ and $f$, for which}
\begin{enumerate}
\item[(i)] \emph{the edge $e$ of $\cB(V)$ is splitted in $\cB(V \cup \{u\})$ into two edges $e_1,e_2$
with endpoints $(x,u)$ and $(u,y)$, respectively, and}
\item[(ii)] \emph{$\cB(V \cup \{u\}) \cap {\rm relint} f = \emptyset$, and instead of the edge
$f$ of $\cB(V), \cB(V \cup \{u\})$ has two new edges $f_1,f_2$ which are the duals of the two edges
$e_1,e_2$ to which $e$ was splitted, respectively}.
\end{enumerate}
\end{enumerate}

\textbf{Proof:}
\begin{enumerate}
\item[a)] $e$ is an arc of the circle $C_{ab} = S_a \cap S_b$, hence $\|a-u\|=\|b-u\|=1$, and since
$e$ is an edge of $\cB(V)$ and $u \in {\rm relint}\,e$, $u$ is not a vertex of $\cB(V)$, hence $\|u-z\|<1$ for
$z \in V \setminus \{a,b\}$. Thus ${\rm diam} (V \cup \{u\}) = 1$ and
$e(V \cup \{u\}) = e(V) +2 = (2n-2) + 2 = 2(n+1) - 2 = e(n+1,3)$ (see Corollary 7.1 above),
i.e., $V \cup \{u\}$ is extremal.
\item[b)] By Lemma 8.1~a) we have $\|u-z\|>1$ for $z \in {\rm relint} f$, i.e., $B (u) \cap {\rm relint} f = \emptyset$,
hence $\cB(V \cup \{u\}) \cap {\rm relint} f = \emptyset$ and $f$ is not an edge of $\cB(V \cup \{u\})$. Now $u$ is a dangling vertex of $
\cB(V \cup \{u\})$, hence the edge $e$ of $\cB(V)$ is splitted by $u$ into two (new) edges $e_1,e_2$ of $\cB(V \cup \{u\})$.
Clearly, their dual edges in $\cB(V) \cup \{u\}$, namely
$f_1$ and $f_2$, respectively, in $\cB(V \cup \{u\})$ are not edges of $\cB(V)$, i.e., they are also new edges
relatively to $\cB(V)$. It remains to show  that every edge $g$ of $\cB(V)$ other  than $e$ and $f$ is also an edge
of $\cB(V \cup \{u\})$. It suffices to show that ${\rm relint} g \subset {\rm int} B(u)$, i.e., that the following
claim holds.
\end{enumerate}

\textbf{Claim:} $\|u-d\|<1$ \emph{for} $d \in {\rm relint}\,g$.

\textbf{Proof:} Let $(s,t)$ be the two endpoints of $g$ and consider two cases:

\textbf{Case I:} $\{s,t\} = \{a,b\}$, i.e., $g$ and $f$ share the same endpoints. Let $h$ be the dual edge of $g$ in $\cB(V)$.
Clearly, $h$ is a short arc of the circle $C_{ab}$, same as $e$ (which is the dual arc of $f$), so both $e$ and $h$
lie on $C_{ab}$, and clearly ${\rm relint}\,e \cap {\rm relint}\,h = \emptyset$. Let $x',y'$ be the endpoints of
$h\,(x',y' \in C_{ab})$.

It follows that $u$, as a point from the relative interior of $e \subset C_{ab}$, lies on the \emph{long} arc between $x'$ and
$y'$ of the circle $C_{ab}$. Since $(a,h,x',y')$ are in equilateral position (Definition 8.2), we can apply Lemma 8.1~a) to
the dual arcs $g=A_{ab} \subset C_{x'y'}$ and $h=A_{x'y'} \subset C_{ab}$, with $c = u \in C_{ab} \setminus A_{x'y'} =
C_{ab} \setminus h$ and $z = d \in {\rm relint} A_{ab} = {\rm relint}\,g$, to conclude that $\|z-d\| < 1$ ($\|c-z\|<1$ in the
notations of Lemma 8.1~a)).

\textbf{Case II:} $\{s,t\} \not= \{a,b\}$. Then $\|u-s\|<1$ or $\|u-t\|<1$ (or both), say $\|u-s\|<1$. Assume, by r.a.a.,
that there exists a $d \in {\rm relint} g$ such that $\|u-d\|>1$. Then there is a point $z \in {\rm relint} g$ on the part of $g$
between $s$ and $d$ such that $\|u-z\|=1$. Then $\|z-a\|=\|z-b\|=\|z-u\|=1$, and $z \in {\rm bd} \cB(V \cup \{u\})$. Thus
$z \in {\rm vert} \cB(V \cup \{u\})$. On the other hand, $z \in {\rm relint} g$, hence $z \notin {\rm vert} \cB(V) = V$ (the
equality holds by Theorem 7.1~(iii) -- the extended GHS-Theorem). It follows that $z \notin V \cup \{u\}$, hence $V \cup \{u\} \not=
{\rm vert} \cB(V \cup \{u\})$. Since $V \cup \{u\}$ is extremal (part a)) this is a contradiction to Theorem 7.1~(iii),
proving the Claim. \hfill \rule{2mm}{2mm}

The following theorem shows how to obtain an extremal set $V \cup W$ from a given extremal set $V$ by
adding a set $W$ of ``dangling vertices'' at will -- just observing the simple rule that no two new points belong
to dual edges of $\cB(V)$.

\textbf{Theorem 8.2:} \emph{Let $V \subset {\mathbb R}^3$ be an extremal
set for the Vazs\'{o}nyi problem with $\infty > \# V = n \ge 4$.
Let $W \subset {\rm skel}_1{\cal B}(V)$, where $W$ is finite such that
for any two dual edges $e,f$ of ${\cal B}(V)$ either $W \cap e =
\emptyset$ or $W \cap f = \emptyset$. Then $V \cup W$ is
extremal}.

\textbf{Proof:} Let $W = \{v_1,v_2,\dots,v_m\}, 1 \le m < \infty$, and let $e_i$ be the edge of $\cB(V)$ that contains
$v_i$, $1 \le i \le m$. Assume, w.l.o.g., that $v_i \in {\rm relint}\,e_i \, (1 \le i \le m)$.
By assumption there are no two dual edges among $e_1,e_2,\dots,e_m$. Put $V_0 =_{\rm def} V\,, \, V_1 =_{\rm def}
V_0 \cup \{v_1\}, V_2 =_{\rm def} V_1 \cup \{v_2\}, \dots, V_m = V_{m-1} \cup \{v_m\} = V \cup W$.

\textbf{Claim:} \emph{$V_i$ is an extremal set for $0 \le i \le m$,
and $e_{i+1}$ is an edge of $\cB(V_i)$ for $0 \le i \le m-1$}.

\textbf{Remark:} The part of the claim that $e_{i+1}$ is an edge of $\cB(V_i)$ is inserted only for technical reasons:
to carry on the induction step (see below).

For $i=0$ there is nothing to prove.

\textbf{Induction step \boldmath$i \longrightarrow i+1\,, \, 0 \le i < m$:\unboldmath}

Assume that $V_i$ is extremal and that $e_{i+1}$ is an edge of $\cB(V_i)$. Consider $V_{i+1} = V_i \cup \{v_{i+1}\}$. Since
$v_{i+1} \in {\rm relint}\,e_{i+1}$ and $e_{i+1}$ is an edge of $\cB(V_i),\,V_{i+1} = V_i \cup \{v_{i+1}\}$ is
extremal, by Lemma 8.2~a). By assumption the edge $e_{i+2}$ is not dual to $e_{i+1}$ (if $i \le m-2$; if $i=m-1$, this part
of the induction proof is finished), and we conclude from Lemma 8.2~b) that $e_{i+2}$ is an edge of $\cB(V_{i+1})$ as well.
This proves the claim.

It follows from the claim that, in particular, $V_m = V \cup W$ is extremal. This proves
Theorem 8.2. \hfill \rule{2mm}{2mm}

\section{Fixed-point free property of \boldmath$\varphi$\unboldmath~ (strong self-duality)}

In this paragraph we define the first barycentric subdivision of
${\cal S}{\cal F}({\cal B}(V))$ for an extremal set $V$ and show that the
canonical self-duality $\varphi$ of ${\cal S}{\cal F}({\cal
B}(V))$ is fixed-point free when it acts as an automorphism of
this barycentric subdivision.

Let $V \subset {\mathbb R}^3$ be an extremal set with $\#V=n \ge
4$. We are about to define the \emph{barycentric subdivision} of
the spherical complex ${\cal S}{\cal F}({\cal B}(V))$. This will
be done in two stages. First we define an abstract two-dimensional
simplicial complex ${\cal K}$. Then we define a (curvilinear)
realization of ${\cal K}$ as a subdivision of ${\cal S}{\cal
F}({\cal B}(V))$.

\textbf{I. Vertices of \boldmath${\cal K}$\unboldmath:}
Corresponding to each non-empty face $\Phi$ of ${\cal S}{\cal
F}({\cal B}(V)), {\cal K}$ has a vertex $z(\Phi)$. Thus we have
vertices $z(x)$ for vertices $x \in V (= {\rm vert}\cB(V)$, by Theorem 7.1 (iii)),
vertices $z(e)$ for edges
$e$ of ${\cal S}{\cal F}({\cal B}(V))$, and vertices $z(F)$ for
facets $F=F_x$ of ${\cal S}{\cal F}({\cal B}(V))$. The
correspondence $\Phi \to z(\Phi)$ is 1-1.

\textbf{II. Faces of \boldmath${\cal K}$\unboldmath:} For each
\emph{chain $($flag} in the current terminology) ${\cal C}$ of
faces of ${\cal S}{\cal F}({\cal B}(V))$ the set $\{z (\Phi):\Phi
\in {\cal C}\}$ is a \emph{face} of ${\cal K}$.

The empty chain yields the empty face $\emptyset$ of ${\cal K}$. Chains
of size 1 yield faces consisting of single vertices of ${\cal K}$.
The edges of ${\cal K}$ ($1$-faces) are of three kinds:
\begin{itemize}
\item $\{z(x),z(e)\}$, where $e$ is an edge of ${\cal S}{\cal
F}({\cal B}(V))$, and $x \in V$ is an endpoint of $e$;
\item $\{z(e),
z(F)\}$, where $F$ is a facet of ${\cal S}{\cal F}({\cal B}(V))$,
and $e$ is an edge of $F$;
\item $\{z(x),z(F)\}$, where $F$ is a
facet, and $x$ is a vertex of $F$.
\end{itemize}
The $2$-\emph{faces} of ${\cal K}$ are of the form $\{z(x), z(e),
z(F)\}$, where $F$ is a facet of ${\cal S} {\cal F}({\cal B}(V))$,
$e$ is an edge of $F$, and $x \in V$ is an endpoint of $e$.

Now we present the realization. We start with the \emph{vertices
of} ${\cal K}$. For each vertex $x \in V$ of ${\cal S}{\cal
F}({\cal B}(V))$, we realize $z(x)$ by $x$.

For each edge $e$ of ${\cal S}{\cal F}({\cal B}(V))$, we realize
$z(e)$ by the midpoint of the circular arc $e$. Each facet $F=F_y$
of ${\cal S}{\cal F}({\cal B}(V))$, $y \in V$, is a compact,
strictly spherically convex subset of the unit sphere $S(y,1)$ (cf. Corollary 5.2 above). We realize
$z(F)$ by a point that lies in the relative interior of the
spherical convex hull of the vertices of $F$ on $S(y,1)$. (When
$F$ is a digon, then the spherical convex hull of its vertices is
just the geodetic arc on $S(y,1)$ that connects the two vertices.)

Now we come to the \emph{edges of} ${\cal K}$.
The edge $\{z(x),z(e)\}$ (where $e$ is an edge of ${\cal S}{\cal F}({\cal B}(V))$, and
$x$ is an endpoint of $e$) is realized by the half of $e$ with endpoints $x,z(e)$.
The edge $\{z(x),z(F_y)\}$ is realized by the geodetic arc on
$S(y,1)$ with endpoints $z(x)$ and $z(F_y)$. Similarly, the edge
$\{z(e),z(F_y)\}$ is realized by the geodetic arc on $S(y,1)$ with
endpoints $z(e)$ and $z(F_y)$.

We continue with the \emph{facets of} ${\cal K}$.
The facet $\{z(x),z(e),z(F_y)\}$ is realized by the union of all
geodetic arcs on $S(y,1)$ that connect $z(F_y)$ with points on the
half of $e$ with endpoints $x$ and $z(e)$. This is the closed
convex triangular region on $S(y,1)$ bounded by the realization of
the edges of $\{z(x),z(e),z(F_y)\}$.

To make the correspondence between (a rectilinear realization of)
${\cal K}$ and the curvilinear realization described above more
explicitly, we introduce the following notation: If $C$ is a
circular arc with endpoints $a$ and $b$, on a circle with center
$c$, we denote by ~``$1-\lambda$'' $ a+$``$\lambda$''$b \,\, (0
\le \lambda \le 1)$ the point $x$ on $C$ that satisfies $\wk acx =
\lambda \wk acb$. If $p = \alpha z(F_y) + \beta z(e) + \gamma
z(x)$ is a point in the triangle $[z(F_y),z(e), z(x)]$ of ${\cal
K}$ $(x \in e, e \in F_x, \alpha,\beta,\gamma \ge 0, \alpha +
\beta + \gamma = 1)$, then the corresponding point $p'$ in the
curvilinear realization is defined in two steps as follows:
\[
q' = {^"}\frac{\beta\,}{1-\alpha}{^"}\,\,z (e) +
{^"}\frac{\gamma}{1-\alpha}{^"}z(x)\,,
\]
where $q'$ lies on $e$, between $z(x)$ and $z(e)$,
\[
p' = "\alpha\," z(F_y) + "1-\alpha\," q'\,,
\]
where $p'$ lies on the great circle through $z(F_y)$ and $q'$ on
$S(y,1)$.

The canonical self-duality $\varphi$ of ${\cal S}{\cal F}({\cal
B}(V))$ induces an automorphism (to be denoted again by $\varphi$)
of order $2$ on the barycentric subdivision ${\cal K}$. To show
that this automorphism is fixed-point free, it is sufficient (and
also necessary) to confirm that no non-empty cell of ${\cal K}$ is
mapped to itself.

Indeed, if a cell $C$ of ${\cal K}$ contains $z(x), x \in V$, then
$\varphi(C)$ contains $\varphi(z(x)) = z(F_x)$ and, hence, does
not contain $z(x)$, since $x$ is \emph{not} a vertex of $F_x$, and
therefore $\varphi(C) \not= C$. Similarly, if $z(F_x) \in C$, then
$z(x) \in \varphi(C)$, and therefore $z(F_x) \notin \varphi (C)$,
hence $\varphi(C) \not= C$.
There remains the case where $C$ consists of a single vertex
$z(e)$ for some edge of ${\cal S}{\cal F}({\cal B}(V))$. But then
$\varphi(C) = \{z(e')\}$, where $e'$ is the edge dual to $e$.
Clearly, $e' \not= e$, hence $\varphi(C) \not= C$ also in this
case.

One can easily check that an even stronger statement holds: If $C$
is a cell of ${\cal K}$, then $C$ and $\varphi(C)$ have no vertex
in common.
We sum the foregoing discussion as

\textbf{Theorem -- Definition 9.1 (strong self duality):} \emph{The canonical $($involutory$)$ self-duality
$\varphi:\cS\cF(\cB(V)) \to \cS\cF(\cB(V))$ for $V$ extremal is fixed point free. An involutory self-duality
of a $2$-dimensional polyhedral $($maybe curvilinearly realized as $\cS\cF(\cB(V))$ above$)$ complex which is fixed
point free is a} strong self-duality. \emph{Thus $\cS\cF(\cB(V))$ is strongly self-dual}.

The following example, illustrating the notion of self-dual complex introduced above where the complex
is that of the boundary of a $3$-polytope,
is a generalization of an example given by Grünbaum and Shephard \cite{G-S}.

\textbf{Example 9.1 (An apexed \boldmath$3$\unboldmath-prism as a
\boldmath$3$\unboldmath-polytope with many involutory self-dualities and only one strong self
duality):}

Let $P$ be a $3$-polytope which is the union of a $3$-prism with (parallel) bases $Q,Q'$
and a $3$-pyramid with basis $Q'$ apexed at $q$, such that aff$Q'$ separates $Q$ and
$q$ (equivalently: relint$[q,v] \cap {\rm relint}Q' \not= \emptyset$ for $v \in {\rm vert}Q$). Note that $Q$
is a facet of $P$ and $Q'$ is not a facet. The polytope $P$ is an \emph{apexed prism} (of dimension $3$).

Put $n =_{\rm def} \#{\rm vert}Q \, (=\# {\rm vert}Q'), n \ge 3$, and denote by $a_1, \dots,
a_n$, resp. $b_1, \dots, b_n$ the vertices of $Q$, resp. $Q'$,
ordered cyclically on relbd$Q$, resp. relbd$Q'$. Assume, w.l.o.g.,
that $[a_i,b_i]$ is an edge of $P$ for $1 \le i \le n$ (otherwise
shift (cyclically) the indices of the $b_i$'s and, if necessary, invert their
order).

\vspace*{2cm}

\begin{picture}(280,130)
\put(40,0){\includegraphics[scale=0.6]{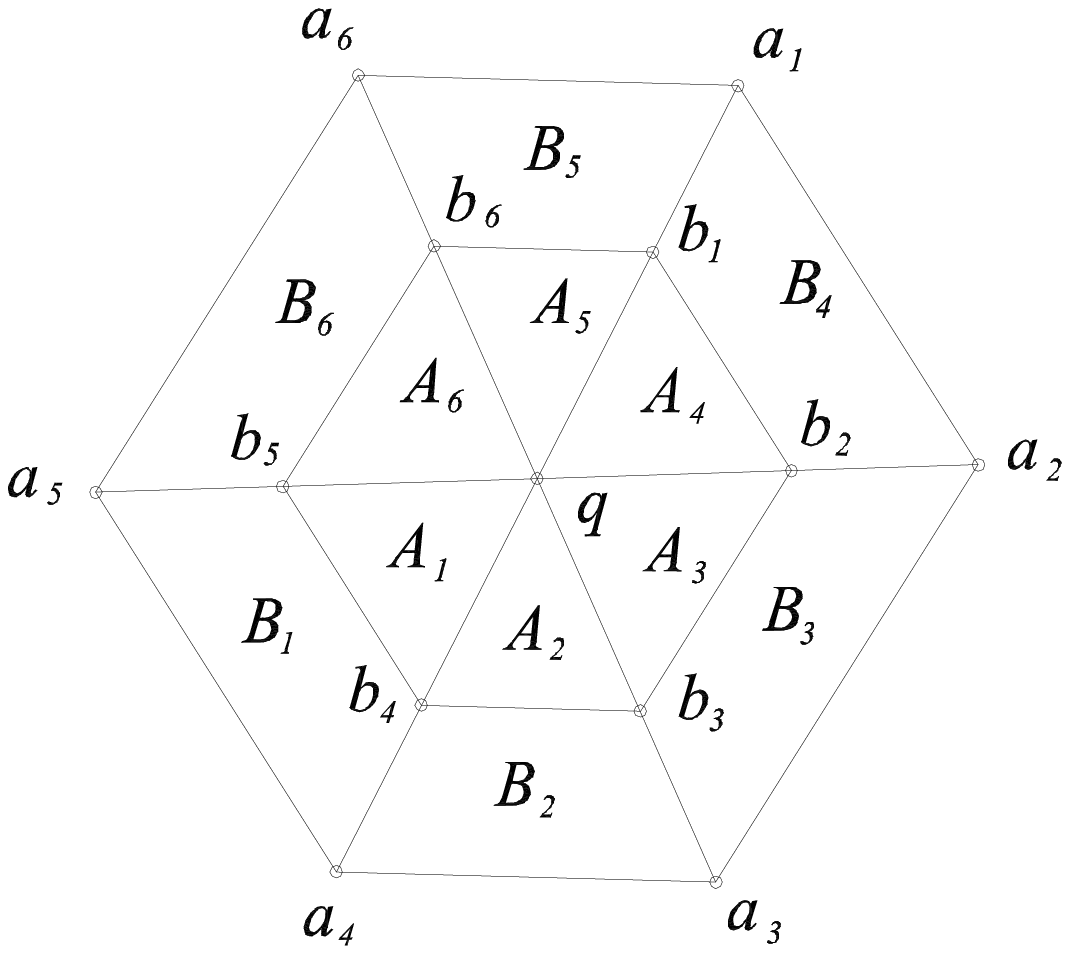}} 
\end{picture}

\vspace{-0.8cm}
\begin{center}
-- Figure 13 --
\end{center}

The $2$-dimensional face structure $\cF(P)$ of $P$ can be drawn
as a Schlegel diagram on the facet $Q$ (see Figure 13 for $n=6$;
ignore the $A_i$'s and $B_i$'s meanwhile). Assume, w.l.o.g., that
$q$ lies ``above'' the plane of the figure and that, looking on the diagram
from the side of $q$ (i.e., from ``above''), the ordering of
$a_1,\dots,a_n$, resp. $b_1,\dots,b_n$, is clockwise. Note that
$a_i$ is $3$-valent and $b_i$ is $4$-valent in skel$_1(P)$.

We will describe now all the involutory self-dualities of $\cF(P)$.
Assume that $\psi:\cF(P) \to \cF(P)$ is an involutory self-duality
of $\cF(P)$. Put $A_i = \psi (a_i)$ and $B_i = \psi (b_i), 1
\le i \le n$.

\textbf{Claim 1:} \emph{$q$ is a vertex of $A_i$ for $1 \le i \le n$}.

\textbf{Proof:} $a_i$ is $3$-valent in skel$_1 P$, hence $\psi
(a_i)$ is a triangular facet.
If $n \ge 4$, then the facets around $q$ are the only triangular
facets of $P$, hence $\psi (a_i) = A_i$ must have $q$ as a vertex.
For $n=3$ an extra consideration is needed, because then $Q$ being
a triangular facet is also a candidate for $A_i = \psi (a_i)$ for
some $i$. Assume, by r.a.a., that $\psi (a_i) = Q$ for some $1 \le
i \le 3$, say, w.l.o.g., $\psi (a_1) = Q$. Then $\psi (a_2)$ (and
$\psi (a_3)$) must have $q$ as a vertex and the edge $[a_1,a_2]$
of $P$ is mapped by $\psi$ to the edge $\psi (a_1) \cap \psi (a_2)
= Q \cap \psi (a_2)$. This is a contradiction, since $Q$ and a
facet having $q$ as a vertex do not share an edge (in fact, they
have an empty intersection).

\textbf{Claim 2:} ~~$q \stackrel{\psi}{\longleftrightarrow} Q$.

\textbf{Proof:} Since $a_i,\, 1 \le i \le n$, are the vertices of
$Q$, their images $\psi (a_i) = A_i$ must all share $\psi (Q)$ as a
vertex. But $A_i,  1 \le i \le n$, share only $q$ by Claim 1,
hence $\psi (Q) = q$.

\textbf{Claim 3:} $B_i \, (1 \le i \le n)$ \emph{and $Q$ share an edge}.

This follows from Claims 1 and 2 above.

\textbf{Claim 4:} $q$ \emph{is a vertex of the edge} $[a_i,a_{i+1}]
\stackrel{\psi}{\longleftrightarrow} A_i \cap A_{i+1} \, (1 \le i
\le n)$.

\textbf{Proof:} $[a_i,a_{i+1}]$ is an edge, hence its image by
$\psi$, $A_i \cap A_{i+1}$, is an edge, and by Claim 1 all
$A_i$'s share $q$ as a vertex. Thus also Claim 4 is verified.

It follows from Claim 4 that the facets $A_1, \dots, A_n$ are
cyclically ordered around $q$, either clockwise (same as the
$a_i$'s) or counterclockwise. Since $[a_1,b_1]$ is an edge, its
image under $\psi$ is the edge $\psi (a_1) \cap \psi (b_1) = A_1
\cap B_1$. Since $A_i \cap B_i$ has to be of the form
$[a_j,a_{j+1})$ for some $1 \le j \le n$ the edge $A_1 \cap B_1$
is of the form $[b_\nu, b_{\nu + 1}]$ for some $1 \le \nu \le n$
($n=6$ and $\nu = 4$ in Figure 16 above).
It follows that either
\begin{enumerate}
\item[(I)] $A_i \cap B_i = [b_{\nu-i+1}, b_{\nu-i+2}], \, 1 \le i
\le n$, or
\item[(II)] $A_i \cap B_i = [b_{\nu+i-1}, b_{\nu+i}], \, 1 \le i \le n$.
\end{enumerate}
We discuss these two cases separately.

\textbf{Case (I):} In this case the $A_i$'s and $B_i$'s $(1 \le i \le n)$ are cyclically ordered counterclockwise
around $q$ (\emph{opposite} to the cyclic order of the $a_i$'s and $b_i$'s).

Clearly, in this case
\begin{equation}\label{stern4}
\begin{array}{l}
{\rm (1)} \quad B_i \cap B_{i+1} = [a_{\nu -i+1}, b_{\nu-i+1}]\,,\\
{\rm (2)} \quad A_i \cap A_{i+1} = [q,b_{\nu-i+1}]\,,\\
{\rm (3)} \quad Q \cap B_i = [a_{\nu-i+1},a_{\nu-i+2}]\,.
\end{array}
\end{equation}

Since $i = \nu-i+1 \Leftrightarrow \nu$ is odd and $i = \frac{\nu+1}{2}$, and
$i = \nu-i+2 \Leftrightarrow \nu$ is even and $i = \frac{\nu+2}{2}$, it follows from (I) that
\begin{equation}\label{stern5}
b_i \in B_i \Leftrightarrow \left\{\begin{array}{c}
i = \frac{\nu + 1}{2} \, \mbox{ if } \, \nu \, \mbox{ is odd}\,,\\
i = \frac{\nu + 2}{2} \, \mbox{ if } \, \nu \, \mbox{ is even}\,.
\end{array}\right.
\end{equation}

Thus in case (I), (\ref{stern4}) + (I) give four necessary conditions for $\psi$ to be an involutory
self-duality of $\cF(P)$.

\textbf{Claim 5:} \emph{The conditions} (\ref{stern4}) + (I) \emph{are $($necessary and$)$ sufficient for $\psi$ to be an involutory
self duality of} $\cF(P)$.

\textbf{Proof:} With the notation of $q,Q,a_i,b_i, A_i, B_i \, (1 \le i \le n)$ and $\nu \,\, (1 \le \nu
\le n)$ satisfying (I) + (\ref{stern4}) introduced above, define an involutory self-duality $\psi_\nu: \cF(P) \to
\cF(P)$ in two stages. First define $\psi_\nu$ on the vertices and facets of $\cF(P)$ by
$q \stackrel{\psi_\nu}{\longleftrightarrow} Q\,, \,\,
a_i \stackrel{\psi_\nu}{\longleftrightarrow} A_i\,,\,\,
b_i \stackrel{\psi_\nu}{\longleftrightarrow} B_i \,\, (1 \le i \le n)$.

In order to define $\psi_\nu$ on the edges of $\cF(P)$, observe first that there
are four types of edges and that each edge has two forms, as a segment with two endpoints and dually as an intersection of
two neighboring facets (in the dual graph of $\cF(P)$). Here is the list of these four types, each
edge is given by its two dual forms (based on (I) + (\ref{stern4})).
\begin{equation}\label{stern6}
\begin{array}{l}
{\rm 1)} \quad [a_i,a_{i+1}] = Q \cap B_{\nu-i+1}\,,\\
{\rm 2)} \quad [a_i,b_i] = B_{\nu-i+1} \cap B_{\nu-i+2}\,,\\
{\rm 3)} \quad [b_i,b_{i+1}] = A_{\nu-i+1} \cap B_{\nu-i+1}\,,\\
{\rm 4)} \quad [q, b_i] = A_{\nu-i+1} \cap A_{\nu-i+2}\,.
\end{array}
\end{equation}
For each edge we use its left hand side form in (\ref{stern6}) to define its image by $\psi_\nu$.
\begin{equation}\label{stern7}
\begin{array}{lll}
{\rm 1)} \quad [a_i,a_{i+1}] &\stackrel{\psi_\nu}{\longleftrightarrow}& A_i \cap A_{i+1}\,,\\
{\rm 2)} \quad [a_i,b_i] & \stackrel{\psi_\nu}{\longleftrightarrow}& A_i \cap B_i\,,\\
{\rm 3)} \quad [b_i,b_{i+1}] & \stackrel{\psi_\nu}{\longleftrightarrow}& B_i \cap B_{i+1}\,,\\
{\rm 4)} \quad [q,b_i] & \stackrel{\psi_\nu}{\longleftrightarrow}& Q \cap B_i\,.
\end{array}
\end{equation}
In order to show that $\psi_\nu$ is well defined on the edges, we must check that the following
four diagrams are commutative.
\[
\begin{array}{l}
{\rm 1)} \quad [a_i,a_{i+1}] ~~ = ~ Q \cap B_{\nu-i+1}\\
\qquad \quad \,\,\,\, \psi_\nu \!\! \searrow \!\!\!\!\!\!\nwarrow ~~~~ \swarrow\!\!\!\!\!\!\nearrow\!\! \psi_\nu\\
\qquad \quad \quad \,\,\,\,\,\,A_i \cap A_{i+1}
\end{array}
\]
\[
\begin{array}{l}
{\rm 2)} \quad [a_i,b_i] ~ = ~ B_{\nu-i+1} \cap B_{\nu-i+2}\\
\qquad \quad \,\, \psi_\nu \!\! \searrow \!\!\!\!\!\!\nwarrow ~~ \swarrow\!\!\!\!\!\!\nearrow\!\! \psi_\nu\\
\qquad \quad \quad \,\,A_i \cap B_i
\end{array}
\]
\[
\begin{array}{l}
{\rm 3)} \quad [b_i,b_{i+1}] ~~ = ~ A_{\nu-i+1} \cap B_{\nu-i+1}\\
\qquad \quad \,\,\,\, \psi_\nu \!\! \searrow \!\!\!\!\!\!\nwarrow ~~~ \swarrow\!\!\!\!\!\!\nearrow\!\! \psi_\nu\\
\qquad \quad \quad \,\,\,B_i \cap B_{i+1}
\end{array}
\]
\[
\begin{array}{l}
{\rm 4)} \quad [q,b_i] ~ = ~ A_{\nu-i+1} \cap A_{\nu-i+2}\\
\qquad \quad  \psi_\nu \!\! \searrow \!\!\!\!\!\!\nwarrow ~~ \swarrow\!\!\!\!\!\!\nearrow\!\! \psi_\nu\\
\qquad \quad \quad \,\, Q \cap B_i
\end{array}
\]

Let us check this for diagram 1). The left double arrow in 1) follows from (\ref{stern7}), 1). To show the right
double arrow, write
\[
\begin{array}{lcl}
Q \cap B_{\nu-i+1} & \stackrel{\psi_\nu}{\longleftrightarrow} & [q,b_{\nu-i+1}]\,\,\, \mbox{(by (\ref{stern7}), 4))}\\
                   & = & A_{\nu-(\nu-i+1)+1} \cap A_{\nu-(\nu-i+1)+2}\,\,\, \mbox{(by (\ref{stern6}), 4))}\\
                   & = & A_i \cap A_{i+1}\,, \,\, {\rm i.e.,}\\
Q \cap B_{\nu-i+1} & \stackrel{\psi_\nu}{\longleftrightarrow} & A_i \cap A_{i+1}\,.
\end{array}
\]
Similarly one can check the commutativity of the remaining diagrams 2), 3), and 4).
This shows that $\psi_\nu \, (1 \le \nu \le n)$ is an involutory self-duality, proving
Claim 5. \hfill \rule{2mm}{2mm}

By (\ref{stern5}), $\psi_\nu$ has a fixed point (when acting as an automorphism of the first
barycentric subdivision of $\cF(P)$): the edge corresponding to the flag $(b_i, B_i)$, where
$b_i \in B_i \,\, (i = \frac{\nu+1}{2}$, resp. $i = \frac{\nu+2}{2}$, for $\nu$ odd, resp. even), is
mapped by $\psi_\nu$ to itself (since $b_i \stackrel{\psi_\nu}{\longleftrightarrow} B_i$). Thus all the
involutory self-dualities $\psi_\nu, 1 \le \nu \le n$, are \emph{not} strong self-dualities.

\textbf{Case (II):} In this case the $A_i$'s and $B_i$'s are cyclically ordered clockwise around
$q$ (\emph{same} as the cyclic order of the $a_i$'s and $b_i$'s).

\vspace*{2cm}

\begin{picture}(280,130)
\put(40,0){\includegraphics[scale=0.6]{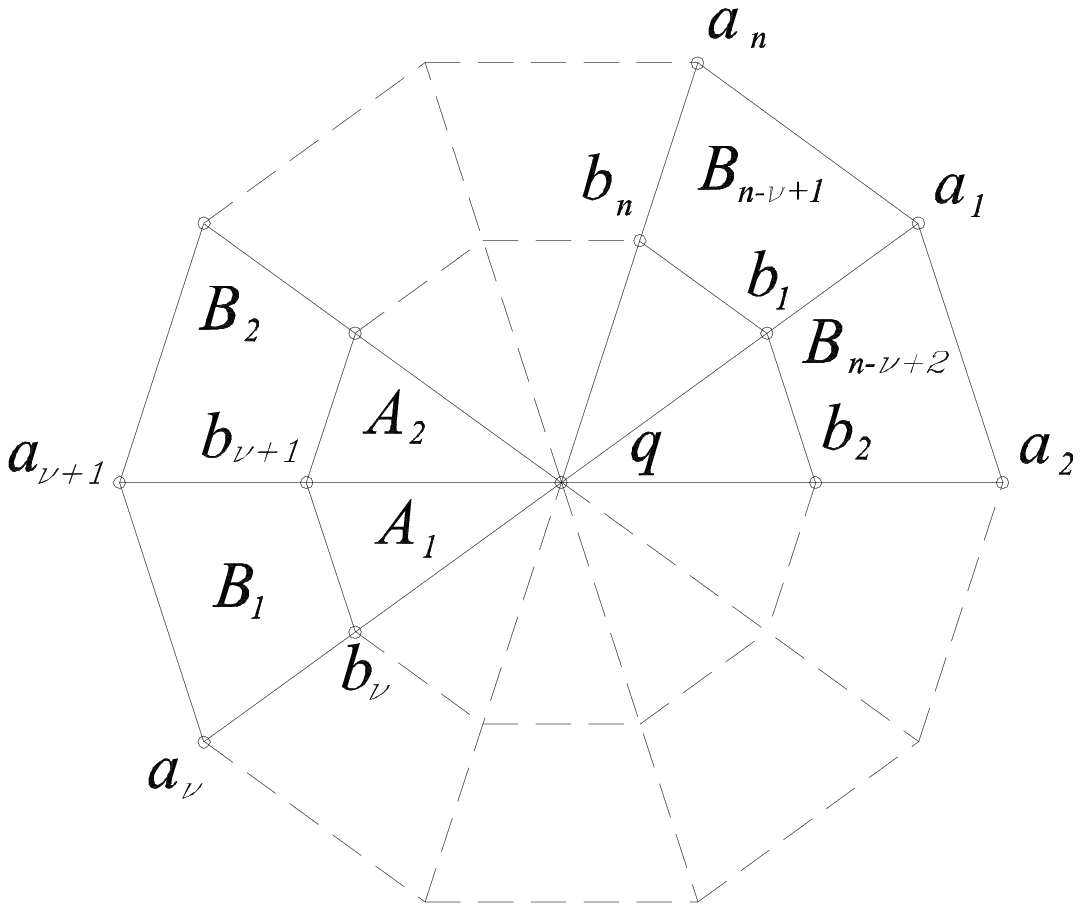}} 
\end{picture}
\vspace{-1cm}

\begin{center}
-- Figure 14 --
\end{center}

Figure 14 sketches the Schlegel diagram of $P$ in this case. We have now
\[
[a_1,b_1] = B_{n-\nu+1} \cap B_{n-\nu+2} \mbox{ and }
[a_1,b_1]  \stackrel{\psi}{\longleftrightarrow}  A_1 \cap B_1 = [b_\nu,b_{\nu+1}]\,,
\]
hence,
\[
[b_\nu,b_{\nu+1}]  \stackrel{\psi}{\longleftrightarrow}  B_{n-\nu+1} \cap B_{n-\nu+2}
                   \stackrel{\psi}{\longleftrightarrow}  [b_{n-\nu+1},b_{n-\nu+2}]\,.
\]
It follows that $[b_\nu,b_{\nu+1}] = [b_{n-\nu+1},b_{n-\nu+2}]$, hence $b_\nu = b_{n-\nu+1}$, i.e.,
$\nu=n-\nu+1$, i.e.,
\begin{equation}\label{stern8}
\nu = \frac{n+1}{2}\,.
\end{equation}
Thus case (II) is possible only if $n$ is odd and $\nu = \frac{n+1}{2}$ (unlike case (I) in which
all $\nu$'s, $1 \le \nu \le n$, are possible). In this case we have, in addition to (II), the following
three conditions (compare to conditions (I) + (\ref{stern4}) above):
\begin{equation}\label{stern9}
\begin{array}{lll}
{\rm 1)} \quad B_i \cap B_{i+1} & = & [a_{\nu+i},b_{\nu+i}]\,,\\
{\rm 2)} \quad A_i \cap A_{i+1} & = & [q,b_{\nu+i}]\,,\\
{\rm 3)} \quad Q \cap B_i & = & [a_{\nu+i-1},a_{\nu+i}]\,.
\end{array}
\end{equation}
Thus the conditions
\begin{equation}\label{stern10}
n \ge 3 \, \mbox{ is odd}, \, \nu = \frac{n+1}{2}, \,\mbox{ (II), and (\ref{stern9})}
\end{equation}
are necessary for $\psi$ to be an involutory self-duality of $\cF(P)$.

\textbf{Claim -- Definition 6:} \emph{The conditions} (\ref{stern10}) \emph{are (necessary and) sufficient for $\psi$ to be an involutory
self-duality of $\cF(P)$. Moreover, these conditions uniquely determine $\psi$. We denote $\psi$, uniquely determined
thus by $\psi^*$}.

\vspace*{2cm}

\begin{picture}(280,130)
\put(40,0){\includegraphics[scale=0.6]{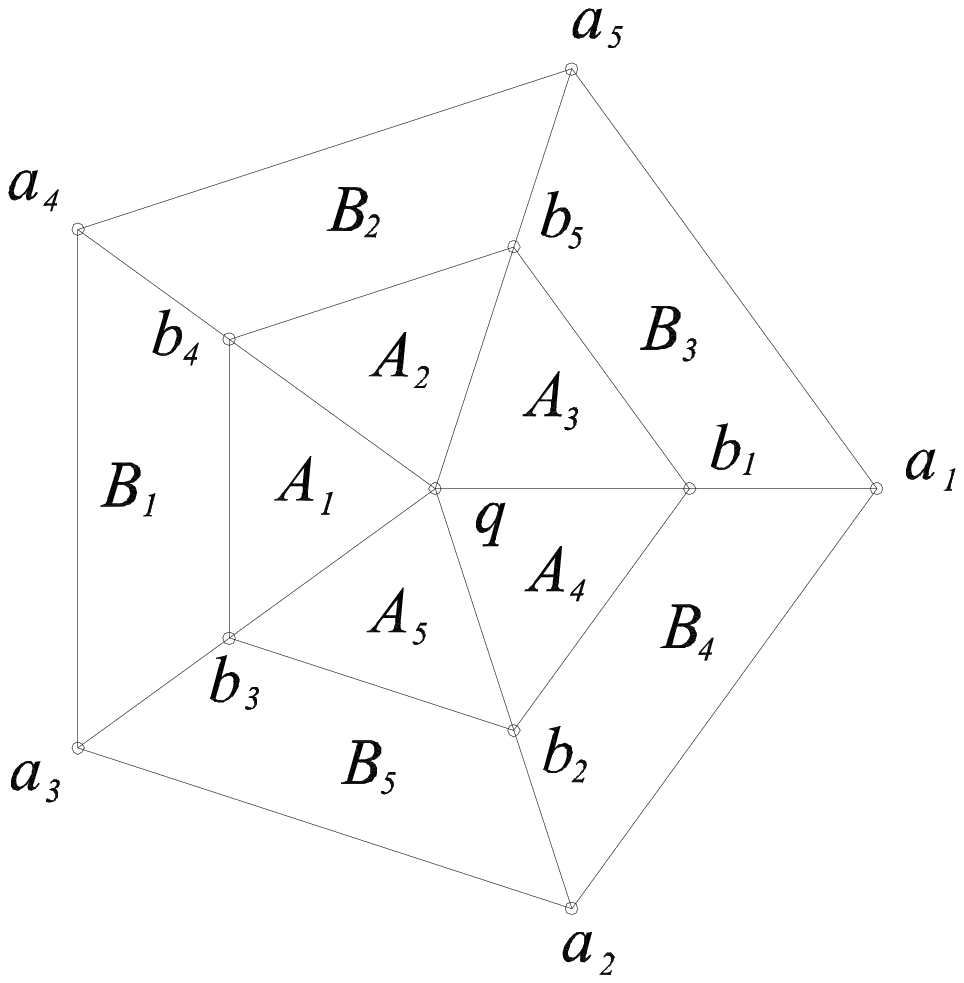}} 
\end{picture}

\vspace{-1cm}
\begin{center}
-- Figure 15 --
\end{center}

\textbf{Proof:} The proof is similar to that of Claim 5 above. One has to write four equalities like (\ref{stern6}), define
$\psi^*: \cF(P) \to \cF(P)$ similarly to the definition of $\psi_\nu$ in (\ref{stern7}),
and then to check the commutativity of four diagrams as above.
We kindly ask the reader to make this check. Figure 15 illustrates the case $n=5$
($\nu = \frac{5+1}{2}=3$).

$\psi^* : \cF(P) \to \cF(P)$ is defined by
$q \stackrel{\psi^*}{\longleftrightarrow} Q,
a_i \stackrel{\psi^*}{\longleftrightarrow} A_i$, and
$b_i \stackrel{\psi^*}{\longleftrightarrow} B_i \,\, (1 \le i \le 5)$,
and on the edges $\psi^*$ is defined as $\psi_\nu$ in (\ref{stern7}).

\textbf{Claim 7:} $\psi^*$ (\emph{assured by} Claim and Definition 6 above) \emph{is fixed point free (as
an induced automorphism of the first barycentric subdivision of} $\cF(P)$).

\textbf{Proof:} Again we leave this for the reader, using the equalities and definition of $\psi^*$ as in
his proof of Claim and Definition 6 above. Alternatively, the reader may check this directly for the case $n=5 \, (\nu = 3)$
depicted in Figure 15 above. \hfill \rule{2mm}{2mm}

Here is a summary of the foregoing discussion:

\textbf{Proposition 9.1:} \emph{Let $P$ be an apexed $3$-prism with $($polygonal$)$ basis of order $n \ge 3$}.
\begin{enumerate}
\item[a)] \emph{The face complex $\cF(P)$ of $P$ has exactly $n$ involutory self-dualities $\psi_\nu,\,\nu=1,2,\dots,n$,
each having a fixed point (in the sense of the barycentric subdivision given above)}.
\item[b)] \emph{For $n \ge 3$ \emph{odd}, $\cF(P)$ has a unique fixed point free involutory self-duality, and for $n > 3$} even
\emph{it has none}.
\item[c)] \emph{Hence $\cF(P)$ is strongly self-dual iff $n \ge 3$ is} odd.
\end{enumerate}
From Proposition 9.1~c) we get

\textbf{Corollary 9.1:} \emph{There is no extremal set $V \subset \bR^3\,, \, 4 \le \# V < \infty$, such that $\cS \cF (\cB(V))$
is isomorphic to the face complex $\cF(P)$ of an apexed prism $P$ whose basis is of an} even \emph{order} $n \ge 4$.

In view of Corollary 9.1, knowing that $\cS \cF(\cB(V))$ is strongly self dual for an extremal set $V$, a natural question
arises. Is there an extremal set $V \subset \bR^3$ such that $\cS \cF (\cB(V))$ is combinatorially equivalent
to the face complex of an
apexed prism whose basis is of an \emph{odd} order $n \ge 3$? The answer is ``yes'', as the following example shows.

\textbf{Example 9.2 (an extremal configuration the face complex of whose ball polytope is combinatorially equivalent to that of
an apexed \boldmath$3$\unboldmath-prism):}
Let $V = \{c\} \cup {\rm vert} P_{2k-1} \, (k \ge 2)$ be the set of vertices of the suspended
$(2k-1)$-gon described in Example 1.2 above. $\cS \cF(\cB(V))$ is combinatorially equivalent to the face complex $\cF(\Pi)$ of
a $3$-pyramid $\Pi$ whose basis is of order $2k-1$, as one can easily check. Put $p_c =_{\rm def} c-\frac{c}{\|c\|}$; $p_c$ is
the central point of the facet $F_c$ of $\cB(V) \, (\|p_c - c\|=1)$, and let $p'_c$ be a point such that
$p'_c =_{\rm def} p_c + \varepsilon \|p_c-c\|$ for some $\varepsilon > 0$.
Note that $p'_c \notin \cB(V)$ and that $p'_c$ lies ``beyond'' the facet $F_c$, relatively to $c$. If $\varepsilon$ is
small enough $(\varepsilon \le 0.1$ will certainly suffice), then $S(p'_c) \cap {\rm skel}_1 \cB(V)$ (where $S(p'_c)$ is the
unit sphere centered at $p'_c$) consists of $2k-1$ points lying on the $2k-1$ edges of $\cB(V)$ incident with
$c$, respectively, all of which are very close to $c$ (certainly of distance not more than $3 \varepsilon$ apart from $c$
for $\varepsilon \le 0.1$). Denote this set of $2k-1$ ``new'' points by $W_{2k-1}$ and define
$W=_{\rm def} \{p'_c\} \cup {\rm vert} P_{2k-1} \cup W_{2k-1}$
or, equivalently, $W =_{\rm def} \{p'_c\} \cup ({\rm vert}\cB(V) \setminus \{c\}) \cup W_{2k-1}$.
Clearly, $\# W = 2n+1$, where $n=_{\rm def} 2k-1$ (odd) and diam $W=1$.

\textbf{Claim:} \emph{$W$ is an extremal configuration for the V\'{a}zsonyi problem}.

\textbf{Proof:} There are $n$ edges of $D(W)$ in $D({\rm vert} P_{2k-1})$, $n$ edges in $D(p'_c \cup W_{2k-1})$,
and $2n$ edges with one endpoint in vert$P_{2k-1}$ and one endpoint in $W_{2k-1}$ (since every vertex of $P_{2k-1}$
is incident with two points in $W_{2k-1}$). Hence $e(W) =n+n+2n=4n=2(2n+1)-2=e(2n+1,3)$, i.e., $W$ is extremal.
This proves the claim.

It follows that $\cS \cF (\cB(W))$ is strongly self-dual (via the canonical self-duality $\varphi$). Moreover,
$\cS \cF (\cB(V))$ is isomorphic to the face complex $\cF(P)$ of an apexed $3$-prism $P$ whose basis is of order
$n\,(=2k-1)$. We leave this check to the reader with just one hint: $p'_c$ plays the role of the apex of $P$,
and $F_{p'_c}$ plays the role of its basis.
Clearly, the canonical self-duality $\varphi$ of $\cS \cF(\cB(V))$ corresponds to the strong self-duality of $\cF(P)$ whose
existence and uniqueness is assured by Proposition 9.1~b) above.

\textbf{Remarks 9.1:}
\begin{enumerate}
\item[1)] As far as we know, the configuration $W$ just now described above
is a completely new type of extremal configuration for the V\'{a}zsonyi problem. Moreover, although
the foregoing construction of $W$ was rigidly described, it allows a
considerable degree of freedom. First, the initial set $V$ of the suspended $(2k-1)$-gon allows some perturbations of
vert$P_{2k-1}$ on $S_c$, as hinted in Example 1.2 above. Second, the point $p'_c$, i.e., the center of the unit ball
$\cB(p'_c)$ by which $c$ was ``truncated'' from $\cB(V)$ to obtain $\cB(W)$, allows some perturbations near
$p_c$ as well.
\item[2)] In fact, this idea of ball truncation is the germ of a vastly general new construction of extremal configurations $W$ given
an extremal configuration $V$: choose a vertex $c \in {\rm vert} \cB(V)$ (= $V$ by Theorem 7.1~(iii)), choose a point
$p_c \in {\rm relint} F_c$, move $p_c$ slightly out of $\cB(V)$ to obtain a point $p'_c$ lying beyond the facet $F_c$
relatively to $c$, and put $W=_{\rm def} \{p'_c\} \cup (\cB(V) \cap S(p'_c)) \cup (V \setminus \{c\}))$. The following
is easy to prove.
\end{enumerate}

\textbf{Claim:} \emph{If $p'_c$ is close enough the $p_c$, then} diam$W=1$ \emph{and $W$ is extremal}.

In a subsequent paper we will exploit the idea of ball truncation and its variants, e.g., that $p'_c$
may lie not only beyond a facet but also
beyond an edge of $\cB(V)$, and even beyond a vertex. Meanwhile we refer to our paper \cite{Ku-Ma-Pe},
where this construction is explained in a more detailed way.

\textbf{Remarks 9.2:} \textbf{Applications to strongly self-dual \boldmath$3$\unboldmath-polytopes}.
\begin{enumerate}
\item[1)] The ``ball truncation'' described above has direct and obvious analogues in the
construction of a strongly
self-dual $3$-polytope $Q$ (abbr. $SSD$ $3$-polytope) given a $SSD$ $3$-polytope $Q$. ($SSD$ means
``Strongly Self-Dual'' in the
sense of \textbf{Theorem -- Definition 9.1} above.) This rich family of $3$-polytopes seems to be an uncharted territory,
and we hope that the present work will prompt a study of it. One aim of such a study might be to construct many (if not all) of the  $SSD$ $3$-polytopes
by using analogues of the ball truncation constructions described above. The model we have in mind for such a
study is Steinitz's proof that every planar $3$-connected graph is polytopal (''Steinitz's theorem'') as it appears
in \cite{Gru}, pp. 235-243. This proof shows in fact that every $3$-polytope $P$ can be obtained from the tetrahedron by a
series of either truncations of $3$-valent vertices or edges or triangular facets (the opposite operations of the ``$\eta_i$-reductions''
in Fig. 13.1.2 in p. 238 there) or addition of a point beyond a triangular facet (the opposite $\omega_i$-reductions
in Fig. 13.1.1 in p. 237 there).
\item[2)] Sallee asks in \cite{Sa}, p. 322, two questions which in the language of the present work can be
formulated as follows.
\end{enumerate}

\textbf{Question 1:} Is the face structure $\cS \cF (\cB(V))$ of every Ball polytope $\cB(V)$
of an extremal set $V \subset \bR^3 \, (4 \le \# V < \infty)$ combinatorially equivalent to a $3$-polytope?

\textbf{Question 2:} Which self-dual $3$-polytopes are combinatorially equivalent to $\cB(V)$ for some extremal $V$?

Here we insert a third question related to Question 1.

\textbf{Question 3:} For which extremal sets $V \subset \bR^3$ is
$\cB(V)$ $3$-polytopal, i.e., has a face structure combinatorially equivalent to that of a $3$-polytope?

The answer to the first question is ``no'' merely because $\cB(V)$ may have a digonal facet $F_x$ for each dangling vertex $x$. So let us modify the
question, by assuming that $\cB(V)$ has no dangling vertices.

\textbf{Question 1':} Is every Ball polytope $\cB(V)$ of a \emph{critical} (extremal) set $V \subset \bR^3 \, (4 \le \# V < \infty)$
(i.e., $V$ is extremal having no dangling vertices; cf. Proposition -- Definition 2.1 above) combinatorially
equivalent to an ordinary $3$-polytope?

The answer to this modified question is again ``no'': as mentioned in \S~1 above, there is a critical configuration $V \subset \bR^3\,,\,
\# V = 8$, for which $\cB(V)$ is not polytopal (in fact, skel$_1 (\cS \cF(\cB(V)))$ is not even $3$-connected).
This configuration will be described in the announced subsequent paper.

An obvious partial answer to Question 2 is: a necessary condition for a $3$-polytope $P$ to be combinatorially equivalent to
$\cB(V)$ for some extremal set $V \subset \bR^3 \, (4 \le \# V < \infty)$ is that $P$ is $SSD$. We conjecture that this goes vice versa,
formulating

\textbf{Conjecture 9.1:} A (necessary and) sufficient condition for the situation that the face structure $\cF(P)$ of
a $3$-polytope $P$ be isomorphic
to $\cS \cF(\cB(V))$ for some extremal (in fact: critical) configuration $V \subset \bR^3$ is that $P$ is $SSD$.

The following conjecture answers Question 3 (provided it is true).

\textbf{Conjecture 9.2:} An extremal set $V \subset \bR^3$ has a polytopal Ball polytope $\cB(V)$ (i.e., $\cB(V)$
is combinatorially equivalent to an ordinary $3$-polytope) if and only if $V$ is strongly critical, i.e., it has no proper subset which is extremal;
cf. Definition 2.2 above.

Clearly, a strongly critical set has no dangling vertex, hence ``strongly critical'' $\Rightarrow$ ``critical''. Good
reasons for why to believe in the truth of Conjecture 9.2 (perhaps even a proof of it) will be given in the mentioned subsequent paper.

For the following conjecture we remind that the \emph{polar body} $P^\triangle$ of a body $P \subset \bR^d$ with $o \in
{\rm int} P$ is $P^\triangle =_{\rm def} \{y \in \bR^d: \langle y,x\rangle \le 1 \, \mbox{ for } \, x \in P\}$,
and that if $P$ is a $d$-polytope with $o \in {\rm int} P$, then the face complexes $\cF(P^\triangle)$ of $P^\triangle$ and
$\cF(P)$ of $P$ are dual.

\textbf{Conjecture 9.3:} Let $Q$ be an $SSD$ $3$-polytope. Then there is a $3$-polytope $P$ combinatorially equivalent to
$Q$ such that $P^\triangle = - P$.

\textbf{Acknowledgement:} The authors wish to thank Dr. Margarita Spirova (Chemnitz) for preparing the figures.

\vspace{0.5cm}

\begin{tabular}{ll}
Yaakov S. Kupitz & Horst Martini\\
Micha A. Perles & Faculty of Mathematics\\
Institute of Mathematics & University of Technology Chemnitz\\
The Hebrew University of Jerusalem & 09107 Chemnitz\\
Jerusalem & GERMANY\\
ISRAEL & \emph{martini@mathematik.tu-chemnitz.de}\\
\emph{kupitz@math.huji.ac.il} &\\
\emph{perles@math.huji.ac.il} &
\end{tabular}

\end{document}